\documentclass[11pt]{article}
\usepackage{graphicx}
\usepackage{amssymb,amsfonts,amsthm}
\usepackage{amsmath,amsthm,amssymb}
\usepackage{amscd}
\usepackage{amsfonts}
\usepackage{color}

\newcommand{\Stab}{{\mathrm{Stab}}}
\newcommand{\lf}{\lfloor}
\newcommand{\rf}{\rfloor}
\newcommand{\la}{\langle}
\newcommand{\ra}{\rangle}

\newtheorem{theorem}{Theorem}[section]
\newtheorem{lemma}[theorem]{Lemma}

\newtheorem{cor}[theorem]{Corollary}
\newtheorem{proposition}[theorem]{Proposition}
\theoremstyle{definition}
\newtheorem{definition}[theorem]{Definition}

\newtheorem{remarks}[theorem]{Remarks}
\newtheorem{remark}[theorem]{Remark}

\newcommand{\hpi}{\widehat{\pi}}
\newcommand{\tpi}{\widetilde{\pi}}

\newcommand{\Cutp}{{\mathrm{Cutp}}\, }
\newcommand{\card}{{\mathrm{card}}\, }
\newcommand{\Int}{{\mathrm{Int}}\, }

\newcommand{\cg}{{\mathcal{G}}}

\newcommand{\qq}{{\mathcal Q}}
\newcommand{\qqq}{{\mathcal Q}}

\newcommand{\Con}{{\mathrm{Con}}}
\newcommand{\lm}{{\lim}}
\newcommand{\dist}{{\mathrm{dist}}}

\newcommand{\Cay}{\mathrm{Cayley}}
\newcommand{\sss}{{\mathcal S}}
\newcommand{\pp}{{\mathcal P}}

\newcommand{\Sat}{{\mathrm{Sat}}\, }
\newcommand{\Sato}{{\mathrm{Sat}}_0\, }
\newcommand{\nn}{{\mathcal N}}
\newcommand{\onn}{\overline{{\mathcal N}}}
\newcommand{\aaa}{{\mathcal A}}
\newcommand{\diam}{{\mathrm{diam}}}

\def\mh{{\mathcal{H}}} %% hiperpl. drapel
\def\calr{\mathcal{R}}   %% R caligrafic
\def\calt{\mathcal{T}}   %% T caligrafic
\def\caltr{{\mathcal{T}}_{\tilde R}}  %% T caligrafic indez tildeR
\def\calc{\mathcal{C}}   %% C caligrafic
\def\ck{\mathcal{K}}   %% K caligrafic

\newcommand {\N}{\mathbb{N}} %% positive integers
\newcommand{\Out}{\mathrm{Out}}
\newcommand{\Aut}{\mathrm{Aut}}
\newcommand{\Inn}{\mathrm{Inn}}
\newcommand {\Z}{\mathbb{Z}}            %% integers
\newcommand {\R}{\mathbb{R}} %% reals
\newcommand {\free}{\mathbb{F}} %% free
\newcommand {\q}{\mathfrak q} %% quasi-geodesic
\newcommand {\tq}{\widetilde{\mathfrak q}} %% quasi-geodesic
\newcommand {\g}{\mathfrak g} %% geodesic
\newcommand {\pgot}{\mathfrak p}
\newcommand {\cf}{\mathfrak c}
\newcommand {\me}{\medskip}

\newcommand {\iv}{^{-1}}
\newcommand{\lio}[1]{\lm^\omega\left(#1\right)}
\newcommand{\co}[1]{\Con^\omega(#1)}

\newcommand {\fn}{\footnote}
\newcommand {\Notat}{\noindent {\it{Notation}}:} %% notatie
\begin{document}

\title{Groups acting on tree-graded spaces and splittings of relatively hyperbolic groups}
\author{Cornelia Dru\c{t}u$^{a}$\thanks{The
authors are grateful to Institut des Hautes \'Etudes Scientifiques
in Paris and to Max-Planck-Institut-f\"ur-Mathematik in Bonn for
their hospitality during the work on this paper.} and Mark
Sapir$^{b,}$\thanks{The work of the second author was supported in
part by the NSF grants DMS 0245600 and DMS-0455881.}}
\date{}
\maketitle

\begin{abstract}
Tree-graded spaces are generalizations of R-trees. They appear as
asymptotic cones of groups (when the cones have cut points). Since
many questions about endomorphisms and automorphisms of groups,
solving equations over groups, studying embeddings of a group into
another group, etc. lead to actions of groups on the asymptotic
cones, it is natural to consider actions of groups on tree-graded
spaces. We develop a theory of such actions which generalizes the
well known theory of groups acting on R-trees. As applications of
our theory, we describe, in particular, relatively hyperbolic groups
with infinite groups of outer automorphisms, and co-Hopfian
relatively hyperbolic groups.
\end{abstract}

\bigskip

\begin{flushright}
Dedicated to Alexander Yurievich Olshanskii's 60th birthday
\end{flushright}

\bigskip
\tableofcontents

\section{Introduction}

\subsection{Actions of groups on real trees}
An $\R$--tree is a geodesic metric space in which every two points
are connected by a unique arc. The theory of actions of groups on
$\R$-trees, started by Bruhat-Tits \cite{Ti} and Bass-Serre
\cite{Ser} in the case of simplicial trees, and continued in the
general situation by Morgan-Shalen and Rips, and then by Sela,
Bestvina-Feighn, Gaboriau-Levitt-Paulin, Chiswell, Dunwoody and
others, turned out to be an important part of group theory and
topology. See Bestvina \cite{Best} for a survey of the theory.

Most applications of the theory are based on the following
remarkable fact first observed by Bestvina \cite{Bes1} and Paulin
\cite{Pau}. If a finitely generated group has ``many" actions by
isometries on a Gromov-hyperbolic metric space, then it acts
non-trivially (i.e. without a global fixed point) by isometries on
the asymptotic cone of that space which is an $\R$-tree. (The word
``many" is explained in Lemma \ref{Pau} below. Basically it means
that there are infinitely many actions that are pairwise
non-conjugate by isometries of the space.) Using theorems of Rips,
Bestvina-Feighn and others, this often allows to write the group
as fundamental group of a non-degenerate graph of groups (see
Section \ref{gat} for precise formulations).

The Bestvina-Paulin's observation applies, for example, given a
hyperbolic group $G$ and infinitely many pairwise non-conjugate
homomorphisms from a finitely generated group $\Lambda$ into $G$.
Then there are ``many" actions of $\Lambda$ on the Cayley graph of
$G$: every homomorphism $\phi\colon \Lambda\to G$ defines an action
$g\cdot x=\phi(g)x$. In particular, this situation occurs if one of
the following holds (see \cite{Best,RS,Se1, Se2, Se3}):

\begin{enumerate}
\item An equation $w(a,b,...,x,y,...)=1$ has infinitely many
pairwise non-conjugate solutions in a group $G$ (here $a,b,...\in
G$, $x,y,...$ are variables, $w$ is a word); then the group
$\Lambda=G*\la x,y,...\ra/(w=1)$ has infinitely many pairwise
non-conjugate homomorphisms into $G$;

\item $\Out(G)$ is infinite;

\item $G$ is not co-Hopfian, i.e. it has a non-surjective but injective endomorphism
$\phi$; then we can consider powers of $\phi$;

\item $G$ is not Hopfian, i.e. it has a non-injective but surjective
endomorphism $\phi$; then we can consider powers of $\phi$.
\end{enumerate}

In cases 2, 3, 4 above $\Lambda =G$. Note that in all these cases,
the group $\Lambda$ acts non-trivially on the asymptotic cone of the
group $G$ even if $G$ is not hyperbolic.

\subsection{Tree-graded spaces}

The asymptotic cones of non-hyperbolic groups need not be trees.
Still, in several important cases they are {\em tree-graded
spaces} in the sense of \cite{DS}. Recall the definition.

\begin{definition}\label{deftgr}
Let $\free$ be a complete geodesic metric space and let $\pp$ be a
collection of closed geodesic subsets (called {\it{pieces}}).
Suppose that the following two properties are satisfied:
\begin{enumerate}

\item[($T_1$)] Every two different pieces have at most one common
point.

\item[($T_2$)] Every simple geodesic triangle (a simple loop
composed of three geodesics) in $\free$ is contained in one piece.
\end{enumerate}

Then we say that the space $\free$ is {\em tree-graded with respect
to }$\pp$.
\end{definition}

In the property ($T_2$), we allow for trivial geodesic triangles,
consequently $\pp$ covers $\free$. In order to avoid extra
singleton pieces, we make the convention that a piece cannot
contain another piece.

Tree-graded spaces have many properties similar to the properties
of $\R$--trees (see our paper \cite{DS} and Section \ref{stg}
below).

As we mentioned in \cite[Lemma 2.30]{DS}, any complete geodesic
metric space with cut-points has a non-trivial tree-graded
structure. Namely, it is tree-graded with respect to maximal
connected subsets without cut-points. Having cut-points in all
asymptotic cones is a weak form of hyperbolicity. In \cite{KK}, it
is proved that super-linear divergence of geodesics in a Cayley
graph of a group implies the existence of cut-points in all
asymptotic cones of the group. Note that Gromov hyperbolicity is
equivalent to superlinear divergence
 of any pair of geodesic rays with common origin \cite{Short}. Here are examples
of groups and other metric spaces whose asymptotic cones have
cut-points:

\begin{itemize}
\item (strongly) relatively hyperbolic groups\fn{In this paper,
when speaking about relative hyperbolicity we shall always mean
strong relative hyperbolicity.} and metrically (strongly)
relatively hyperbolic spaces \cite{DS};

\item Mapping class groups of
punctured surfaces such that $3$ times the genus plus the number of
punctures is at least 5 \cite{Behrstock};

\item Teichm\"uller spaces with Weil-Petersson metric if
$3$ times the genus of the corresponding surface plus the number of
punctures is at least 6 \cite{Behrstock};

\item Fundamental groups of graph-manifolds which are not Sol nor Nil manifolds
\cite{KK};

\item Many right angled Artin groups \cite{BDM}.

\end{itemize}

Note that mapping class groups, Teichm\"uller spaces and fundamental
groups of graph-manifolds are not relatively hyperbolic \cite{BDM}
neither as groups nor as metric spaces (that is, there exists no
finite family of finitely generated subgroups with respect to which
they are hyperbolic, and no family of subsets with respect to which
they are metrically relatively hyperbolic in the sense of
\cite{DS}).

\subsection{Main results}

The above examples give a motivation for the study of actions of
groups on tree-graded spaces. In this paper, we show that a group
acting ``nicely" on a tree-graded space also acts ``nicely" on an
$\R$--tree, and so the group can be represented as the fundamental
group of a graph of groups with ``reasonable'' edge and vertex
groups. In order to explain the last phrase, we recall the
following definitions.

In \cite{DS}, we proved that for every point $x$ in a tree-graded
space $(\free, \pp)$, the union of geodesics $[x,y]$ intersecting
every piece by at most one point is an $\R$--tree called a {\em
transversal} tree of $\free$. A geodesic $[x,y]$ from a
transversal tree is called {\em transversal geodesic}.

\medskip

{\em Notation:} For every group $G$ acting on a tree-graded space
$(\free,\pp)$,

\begin{itemize}

\item $\calc_1(G)$ is the set of stabilizers of subsets of $\free$ all
of whose finitely generated subgroups stabilize pairs of distinct
pieces in $\pp$;

\item $\calc_2(G)$ is the set of stabilizers of pairs of points of
$\free$ not from the same piece;

\item $\calc_3(G)$ is the set of stabilizers of triples of points of
$\free$, neither from the same piece nor on the same transversal
geodesic.

\end{itemize}

Here is our main result about groups acting on tree-graded spaces.

\begin{theorem}[Theorem \ref{1}]\label{11} Let $G$ be a finitely generated group
acting by isometries on a tree-graded space $(\free,\pp)$.
 Suppose that the following hold:
 \begin{itemize}
 \item[\textbf{(i)}] every isometry $g\in G$ permutes the pieces;
 \item[\textbf{(ii)}] no piece in $\pp$ is stabilized by the whole group $G$;
 likewise no point in $\free$ is fixed by the whole group $G$.
 \end{itemize}
 Then one of the following four situations occurs:
\begin{itemize}
\item[(\textbf{I})] the group $G$ acts by isometries on an $\R$--tree non-trivially, with stabilizers of non-trivial
 arcs in $\calc_2(G)$, and with stabilizers of non-trivial tripods in $\calc_3
 (G)$;
 \item[(\textbf{II})] there exists a point $x\in \free$ such that for
 any $g\in G$ any
 geodesic $[x,g\cdot x]$ is covered by finitely many pieces: in this case
 the group $G$ acts non-trivially on a simplicial tree with
 stabilizers of pieces or points of $\free$ as  vertex stabilizers, and stabilizers of
 pairs (a piece, a point inside the piece) as edge stabilizers;
\item[(\textbf{III})] the group $G$ acts non-trivially on a
simplicial tree with edge stabilizers from $\calc_1(G)$;
 \item[(\textbf{IV})] the group $G$ acts on a complete
{$\R$--}tree by isometries, non-trivially, such that stabilizers
of non-trivial arcs are locally inside $\calc_1(G)$-by-Abelian
subgroups, and stabilizers of tripods are locally inside subgroups
in $\calc_1(G)$; moreover if $G$ is finitely presented then the
stabilizers of non-trivial arcs are in $\calc_1(G)$.
\end{itemize}
\end{theorem}

Using results of Rips, Bestvina-Feighn, Sela and Guirardel (see
Section \ref{gat}), we apply Theorem \ref{11} to deduce from actions
on tree-graded spaces splittings of the group.

If the group $G$ acting on a tree-graded space is finitely presented
then the following theorem holds.

\begin{theorem}[see Theorem \ref{split1}]\label{split3}
Let $G$ be a finitely presented group acting by isometries on a
tree-graded space $(\free ,\pp)$ such that
\begin{itemize}
 \item[\textbf{(i)}] every isometry from $G$ permutes the pieces;
 \item[\textbf{(ii)}] no piece in $\pp$ is stabilized by the whole group $G$;
 likewise no point in $\free$ is fixed by the whole group $G$;
 \item[\textbf{(iii)}] the collection of subgroups
 $\calc(G)=\calc_1(G)
 \cup \calc_2 (G)$  satisfies the ascending chain condition.
 \end{itemize}
 Then one of the following three cases occurs:
  \begin{itemize}
  \item[(1)] $G$ splits over a
  $\calc(G)$-by-cyclic group;
  \item[(2)] $G$ can be represented as the fundamental group of a
  graph of groups whose vertex groups are stabilizers of pieces and
  stabilizers of points in $\free$, and edge groups are stabilizers
  of pairs (a piece, a point in the piece);
  \item[(3)] the group $G$ has a $\calc_1(G)$-by-(free Abelian) subgroup
  of index at most 2.
  \end{itemize}
\end{theorem}

With the help of a yet unpublished result of Guirardel, Theorem
\ref{gui}, the following weaker version of Theorem \ref{split3} for
all finitely generated groups is proved. It  is still sufficient for
our applications to relatively hyperbolic groups.

\begin{theorem}[Theorem \ref{splitg}]\label{splitg1}
Let $G$ be a finitely generated group acting by isometries on a
tree-graded space $(\free ,\pp)$ such that properties \textbf{(i)}
and \textbf{(ii)} from Theorem \ref{split3} hold, moreover
\begin{itemize}
\item[\textbf{(iii)}] the subgroups in $\calc_2(G)$ are
 (finite of uniformly bounded size)-by-Abelian
 and the subgroups in $\calc_1(G)\cup \calc_3(G)$
 have uniformly bounded size.
 \end{itemize}
 Then one of the following three cases occurs:
  \begin{itemize}
  \item[(1)] $G$ splits over a
  [(finite of uniformly bounded size)-by-Abelian]-by-(virtually cyclic)
  subgroup;
  \item[(2)] same as case (2) from Theorem \ref{split3};
\item[(3)] the group $G$ has a subgroup
of index at most 2 which is a [(finite of uniformly bounded
size)-by-Abelian]-by-(free Abelian) subgroup.
 \end{itemize}
\end{theorem}

Note that if $G$ is torsion-free then one can weaken the
assumption on $\calc_2(G)$ from Theorem \ref{splitg1} by using
Sela's theorem from \cite{Se} instead of Guirardel's Theorem
\ref{gui}. The conclusion in this case is stronger:

\begin{theorem}[see Theorem \ref{split5}]\label{split5g}
Let $G$ be a finitely generated group acting by isometries on a
tree-graded space $(\free ,\pp)$ such that properties \textbf{(i)}
and \textbf{(ii)} from Theorem \ref{split1} hold, and in addition
\begin{itemize}
\item[\textbf{(iii)}] the collection of subgroups
 $\calc_2(G)$ satisfies the ascending chain condition
 and every subgroup in $\calc_1(G)\cup \calc_3(G)$ is trivial.
 \end{itemize}
 Then one of the following three cases occurs:
  \begin{itemize}
  \item[(1)] $G$ splits over a
  $\calc_2(G)$-by-cyclic group or an Abelian-by-cyclic group;
  \item[(2)] same as case (2) from  Theorem
\ref{split3};
  \item[(3)] the group $G$ has a metabelian subgroup
  of index at most 2.
  \end{itemize}
\end{theorem}

\subsection{The case of relatively hyperbolic groups}

As the authors together with D. Osin showed in \cite{DS}, relatively
hyperbolic groups provide  natural examples of groups with
tree-graded asymptotic cones. We describe now some applications of
Theorems \ref{split3} and \ref{splitg1} to relatively hyperbolic
groups.

Questions related to homomorphisms into relatively hyperbolic groups
have attracted considerable attention especially in the case of
hyperbolic groups, Kleinian groups, limit (in another terminology
fully residually free) groups and more generally relatively
hyperbolic groups with Abelian peripheral subgroups.

{Recall that Z. Sela has proved in \cite{Se1} that a hyperbolic
group which is non-elementary and torsion free is co-Hopfian if and
only if it is freely indecomposable.} Co-Hopf geometrically finite
Kleinian groups (these are hyperbolic relative to Abelian subgroups)
have been described by T. Delzant and L. Potyagailo in \cite{DP}. A
description of all homomorphisms of a finitely presented group into
a relatively hyperbolic group via Makanin-Razborov diagrams has been
provided in the case of limit groups \cite{Ali}, and more generally
in the case of torsion-free groups that are hyperbolic relative to
Abelian subgroups \cite{Gro3}. The structure of $\Out(G)$ has been
clarified in the case of limit groups in \cite{BKM}, and in the case
of relatively hyperbolic groups in \cite{Gro2}.

One of the strongest known results about homomorphisms into
relatively hyperbolic groups is due to Dahmani.

\begin{definition}\label{dah1} Following Dahmani \cite{Dah}, we say that a homomorphism
$\phi$ from a group $\Lambda$ into a relatively hyperbolic group
$G$ has an {\em accidental parabolic} if either $\phi(\Lambda)$ is
parabolic {(in which case $\phi$ is called a {\em parabolic
homomorphism})} or $\Lambda$ splits over a subgroup $C$ such that
$\phi(C)$ is either parabolic or finite.\fn{It is easy to see that
this definition is equivalent to the definition in \cite{Dah}.}
\end{definition}

Dahmani proved in \cite{Dah} that if $\Lambda$ is finitely
presented, and $G$ is relatively hyperbolic then there are
finitely many subgroups of $G$, up to conjugacy, that are images
of $\Lambda$ in $G$ by homomorphisms without accidental
parabolics.

Instead of homomorphic images, we consider the set of
homomorphisms. {Recall that a group $\Lambda$ satisfies {\em
property FA} of Serre \cite{Ser} if every action of $\Lambda$ on a
simplicial tree has a global fixed point. The class of groups with
property FA includes all finite groups, finite index subgroups in
$\mathrm{SL}_n(\Z)$, $n\ge 3$, and, more generally all groups with
Kazhdan property T, $\Aut(F_n)$, $n\ge 3$, mapping class groups of
surfaces with positive genus or of spheres with at least 5
punctures \cite{Bo,CV}, subgroups of finite index in Chevalley
groups of rank $\ge 2$ over rings of integers  \cite{Ti}, etc.}

\begin{theorem}[See Corollary \ref{corcor}]\label{corcor1} If a
finitely generated group $\Lambda$ satisfies property FA then for
every relatively hyperbolic group $G$ there are only finitely many
pairwise non-conjugate non-parabolic homomorphisms $\Lambda\to G$.
\end{theorem}

If the group $\Lambda$ {is the fundamental group of a non-degenerate
graph of groups} (i.e. it does not have property FA), then the
structure of homomorphisms into relatively hyperbolic groups is more
complicated. Note that if a group $G$ splits over an Abelian
subgroup $C$, say $G=A*_C B$, then it typically has many outer
automorphisms that are identity on $A$ and conjugate $B$ by elements
of $C$. Hence in order to get a finiteness result for homomorphisms
up to conjugacy, we need to modify the definition of accidental
parabolics as follows.

\begin{definition}\label{dahdef1} We say that a homomorphism
$\phi\colon \Lambda\to G$ has a {\em weakly accidental parabolic} if
one of the following two cases occurs:
\begin{itemize}
    \item $\phi(\Lambda)$ is either parabolic or virtually cyclic;
    \item $\Lambda$ splits over a subgroup $C$ such that $\phi(C)$ is
    either parabolic or virtually cyclic.
\end{itemize}
\end{definition}

Here we formulate some of our results not in the strongest
possible form. To simplify the formulation, we impose the
restriction that peripheral subgroups do not contain free
non-Abelian subgroups (see Section \ref{homs} where this condition
is removed). From now on till the end of the section, $G$ is a
(strongly) relatively hyperbolic group.

\begin{theorem}[An immediate corollary of Theorem \ref{prop888} and Remark
\ref{rem888}.]\label{1Gdah4}
Suppose that the peripheral subgroups
of $G$ do not contain free non-Abelian subgroups. Let $\Lambda$ be
a finitely generated subgroup of $G$ which is neither virtually
cyclic nor parabolic. Assume moreover that $\Lambda$ does not
split over a parabolic subgroup nor over a virtually cyclic
subgroup.

Then there are finitely many conjugacy classes in $G$ of injective
homomorphisms $\Lambda\to G$ whose image is not parabolic.
\end{theorem}

Dahmani's theorem from \cite{Dah} and Theorem \ref{1Gdah4}
immediately imply the following corollary.

\begin{cor}\label{Gdah10} Let $\Lambda$ be a finitely presented group
and assume that the peripheral subgroups of $G$ do not contain
non-Abelian free subgroups. Then there are finitely many conjugacy
classes in $G$ of homomorphisms $\psi\colon \Lambda\to G$ without
weakly accidental parabolics and with Hopfian images.
\end{cor}

\proof  By Dahmani \cite{Dah}, up to conjugacy in $G$, there are
finitely many subgroups that are images of $\Lambda$ by
homomorphisms $\psi\colon \Lambda\to G$ without accidental
parabolics. So we may assume that the image of $\psi$ is fixed.
Then $\psi\colon \Lambda\to G$ is the composition of a fixed
surjective homomorphism $\psi_1\colon \Lambda\to \psi_1 (\Lambda)$
and a homomorphism $\psi_2$ from $\Lambda_1=\psi_1 (\Lambda)$ onto
itself. If $\Lambda_1$ is Hopfian then $\psi_2$ is necessarily
injective and we can apply Theorem \ref{1Gdah4}.
\endproof

We do not know if a non-parabolic finitely generated subgroup of a
relatively hyperbolic group can be non-Hopfian provided peripheral
subgroups are Hopfian. Hence it might be possible to strengthen
Corollary \ref{Gdah10} at least in the case of Hopfian peripheral
subgroups. If peripheral subgroups are not Hopfian, the whole group
can be non-Hopfian as well: consider a free product $A*B$ where $A$
is non-Hopfian.

Theorem \ref{1Gdah4} also immediately implies the following
corollary. For every subgroup $\Lambda<G$ let $N_G(\Lambda)$ (resp.
$C_G(\Lambda)$ ) be the normalizer (resp. the centralizer) of
$\Lambda$ in $G$. Clearly there exists a natural embedding
$\varepsilon$ of $N_G(\Lambda)/C_G(\Lambda)$ into the group of
automorphisms $\Aut(\Lambda)$.

\begin{cor}[See Corollary \ref{cor646} and Remark \ref{rem888}]
\label{Gcor646}
Suppose that $\Lambda\le G$ is neither virtually cyclic nor
parabolic, and that it does not split over a parabolic or virtually
cyclic subgroup. Suppose that peripheral subgroups of $G$ do not
contain free non-Abelian subgroups.

 Then $\varepsilon(N_G(\Lambda)/C_G(\Lambda))$ has finite
index in $\Aut(\Lambda)$. In particular, if $\Out(\Lambda)$ is
infinite then $\Lambda$ has infinite index in its normalizer.
\end{cor}

In the case when $\Out(G)$ is infinite or $G$ is not co-Hopfian, we
need weaker or no extra assumptions on the peripheral subgroups to
obtain splittings of $G$.

\begin{theorem}[Theorem \ref{thout}, Remark \ref{remark459}] \label{thout1}
Suppose that the peripheral subgroups of $G$ are not relatively
hyperbolic with respect to proper subgroups. If $\Out(G)$ is
infinite then one of the followings cases occurs:
\begin{itemize}
    \item[(1)] $G$ splits over a virtually cyclic subgroup;
    \item[(2)] $G$ splits over a parabolic
  (finite of uniformly bounded size)-by-Abelian-by-(virtually cyclic)
  subgroup;
    \item[(3)] $G$ can be represented as a non-trivial amalgamated product or
    HNN extension with one of the vertex groups a maximal
    parabolic subgroup of $G$.
\end{itemize}
\end{theorem}

\begin{remark}
If a peripheral subgroup is relatively hyperbolic then, according to
Lemma \ref{subgroups}, it can be replaced in the list of peripheral
subgroups by its own peripheral subgroups. In an arbitrary
relatively hyperbolic group this process of getting smaller and
smaller peripheral subgroups might not terminate. A example is that
of a finitely generated non-Abelian free group hyperbolic relative
to a finitely generated non-Abelian free subgroup. The free subgroup
in its turn can be replaced by a proper free subgroup, and this
process can continue indefinitely. Still, for this example there
exists a terminal point of the process, as the free group is
hyperbolic relative to the trivial subgroup.

One may then ask if every relatively hyperbolic group has such a
terminal structure, that is if it has a list of peripheral subgroups
that are not relatively hyperbolic. It turns out that the answer is
negative. An example is the inaccessible group of Dunwoody
\cite{Dun1} (the argument showing it is in \cite{BDM}). It is not
known if torsion free or finitely presented examples of this sort
exist.
\end{remark}

\begin{theorem}[Theorem \ref{cohopf}]\label{cohopf1}
Suppose that a relatively hyperbolic group $G$ is not co-Hopfian.
Let $\phi$ be an injective but not surjective homomorphism $G\to G$.
Then one of the following holds:

\begin{itemize}
\item $\phi^k(G)$ is parabolic for some $k$.

\item $G$ splits over a parabolic or virtually cyclic subgroup.
\end{itemize}
\end{theorem}

Various weaker versions of Theorems \ref{1Gdah4}, \ref{thout1} and
\ref{cohopf1} (with various restrictions on parabolic subgroups)
have been proved before by K. Ohshika and L. Potyagailo \cite{OP},
D. Groves (\cite{Gro1, Gro2, Gro3}), I. Belegradek and A.
Szczepa\'{n}ski \cite{BS}. Instead of the asymptotic cones of the
relatively hyperbolic group itself, Belegradek and Szczepa\'{n}ski
used actions of relatively hyperbolic groups on locally compact
hyperbolic spaces. This led them in, say, the case when $\Out(G)$ is
infinite or $G$ is not co-Hopfian, to an action of $G$ on an
$\R$--tree (the asymptotic cone of the hyperbolic space) with
parabolic stabilizers. But since they do not have control on the
stabilizers, they have to assume that all subgroups of peripheral
subgroups are finitely generated (i.e. the peripheral subgroups are
slender).

\subsection{Possible further applications}

Relatively hyperbolic groups form a small subclass of the class of
groups whose asymptotic cones have cut points (we call such groups
{\em constricted} in \cite{DS}). For example, as we mentioned above,
mapping class groups, fundamental groups of graph manifolds are also
constricted. It would be very tempting to apply our results about
groups acting on tree-graded spaces to obtain statements similar to
Theorem \ref{prop888} for mapping class groups and other constricted
groups. By Lemma \ref{Pau}, if $G$ is constricted, and $\Lambda$ has
infinitely many pairwise non-conjugate homomorphisms into $G$ then
$\Lambda$ acts on an asymptotic cone $\cal C$ of $G$ with no
globally fixed points. In order to apply Theorem \ref{11}, one would
need to know that the action also does not stabilize any piece in
the natural tree-graded structure of $\cal C$, i.e. maximal
connected subsets of $\cal C$ without cut-points. Then in order to
apply Theorem \ref{split3}, one also needs to know the stabilizers
of the pieces in $\Lambda$. In Section \ref{srelhyp}, we show how
this strategy works in the case of relatively hyperbolic groups.

\subsection{Organization of the paper}
The plan of the paper is the following.

In Section 2, we present some general facts about tree-graded
spaces. In particular, we define an $\R$--tree appearing as a
natural quotient of a tree-graded space. It can be described as the
factor-space $\free/\approx$, where $x\approx y$ if in any geodesic
joining $x$ and $y$ non-trivial subarcs contained in pieces compose
a dense subset. We also define a natural construction that can be
applied to any tree-graded space $(\free,\pp)$ and gives a
tree-graded space $(\free, \pp')$ with ``bigger" pieces. In the same
section we also recall results about groups acting on $\R$--trees
(theorems of Bestvina-Feighn, Levitt, Sela and Guirardel). We prove
a version of Levitt's theorem about actions of finitely generated
groups on trees by homeomorphisms (Theorem \ref{lev1}).

In Section \ref{gatg}, we develop our theory of actions on
tree-graded spaces and prove Theorem \ref{11}. The proof proceeds
as follows. Let $G$ be a group acting on a tree-graded space
$(\free, \pp)$ such that conditions (i) and (ii) from Theorem
\ref{11} hold. First we analyze the situation (Case A) when $G$
acts on the $\R$--tree $\free/\approx$ non-trivially. Then we show
that the stabilizers of arcs of that $\R$--tree are from $\calc_2
(G)$. This gives Case (I) from Theorem \ref{11}.

In Case B, $G$ acts on $\free/\approx$ with a global fixed point.
The corresponding $\approx$--equivalence class $R$ stabilized by
$G$ is also a tree-graded space with a collection of pieces $\calr
\subset \pp$. We define an increasing transfinite sequence of
tree-graded structures on $R$,
$\pp_0=\calr<\pp_1=\pp_0'<...<\pp_\alpha$, using the construction
from Section \ref{gat} and corresponding equivalence relations
$\sim_0$, $\sim_1$,... on $R$ (two points $a, b$ of $R$ are
$\sim_\beta$--equivalent if a geodesic $[a,b]$ is covered by
finitely many pieces from $\pp_\beta$). This sequence of
tree-graded structures must stabilize for some cardinal $\alpha$.
We consider two subcases: Case B1 is when $G$ stabilizes a piece
in $\pp_\alpha$, and Case B2 when $G$ does not stabilize a piece
there.

In Case B1, we take the minimal cardinal $\delta\le \alpha$ such
that $G$ fixes a piece $A$ in $\pp_\delta$. We show that $\delta$ is
not a limit cardinal, so $\delta-1$ exists. Then we construct a
simplicial tree and an action of $G$ on it. The vertices of the tree
are pieces of $\pp_{\delta-1}$ intersecting non-trivially a copy in
$\free$ of the Cayley graph of $G$, and points of intersection of
pieces. Edges connect a piece and an intersection point contained in
that piece. If $\delta=1$ then we get Case (II) of Theorem \ref{11},
if $\delta>1$ then we get Case (III).

In Case B2, the pieces of $\pp_\alpha$ do not intersect. Then the
$\pp_\alpha$--pieces have a natural structure of pre-tree in the
sense of Bowditch \cite{Bow} and the group $G$ acts on this
pre-tree by automorphisms. The pre-tree is $G$-isomorphic (as a
pre-tree with a $G$-action) to an $\R$--tree. We prove this by
first embedding it $G$-equivariantly into an $\R$--tree, and then
showing that the embedding is surjective. (For this embedding we
might have used a result of Bowditch and Crisp \cite{BC}, but we
found an easy proof of the existence of such an embedding before
we were aware of the reference \cite{BC}).) Then we use our
version of Levitt's theorem (Theorem \ref{lev1}), and obtain an
action of $G$ on an $\R$--tree by isometries with ``good" arc
stabilizers (so that Case (IV) of the theorem holds).

In Section \ref{srelhyp}, we consider applications of our results to
relatively hyperbolic groups. First we recall results from \cite{DS}
about relatively hyperbolic groups and asymptotically tree-graded
spaces and prove some modifications of these results. Then we
consider the isometry groups of asymptotic cones of relatively
hyperbolic groups and describe stabilizers of pieces, pairs of
pieces and triples of points of the cone. This allows us to apply
Theorems \ref{11} and \ref{splitg1} and deduce Theorems \ref{thout1}
and \ref{cohopf1}.

{\bf Acknowledgement.} We are grateful to I. Belegradek, F. Dahmani,
V. Guirardel, I. Kapovich, M. Kapovich, G. Levitt, F. Paulin and Z.
Sela for many helpful conversations.

\section{Preliminaries}\label{prel}

\subsection{Notation and definitions}\label{ndef}

Throughout this paper, we shall use the following notation.

For every path $\pgot$, we denote the start of $\pgot$ by
$\pgot_-$ and the end of $\pgot $ by $\pgot_+$. \me

Let $X$ be a metric space. Given two quasi-geodesics in $X$,
$\pgot\colon [0,a]\to X$ and $\q\colon [0,b]\to X$, such that
$\pgot_+=\q_-$, we denote by $\pgot\sqcup \q$ the map from
$[0,a+b]$ to $X$ that is equal to $\pgot$ on $[0,a]$ and to
$\q[t-a]$ on $[a,a+b]$.

If $a,b$ are two points in $X$ then $[a,b]$ denotes any geodesic
$\q$ with $\q_-=a, \q_+=b$.

If $x$ is a point in $X$ and $r\ge 0$ then $B(x,r)$
 denotes the ball of radius $r$ around $x$ in $X$.

For every $Y\subseteq X$, $r\ge 0$, $\nn_r(Y)$ and $\onn_r(Y)$
denote the open and respectively the closed $r$-tubular neighborhood
of $Y$ in $X$.

Recall that a group is said to have some property \textit{locally}
if any finitely generated subgroup of it has this property.

\me

\subsection{Tree-graded metric spaces}\label{stg}

\begin{lemma}[\cite{DS}, Proposition 2.17]\label{pt2}
The property ($T_2$) in the definition of tree-graded spaces (see
Definition \ref{deftgr}) can be replaced by the assumption that
$\pp$ covers $\free$ {together with} the following property (which
can be viewed as an extreme version of the bounded coset
penetration property of \cite{Fa}):

\begin{quotation}
($T_2'$)\quad for every topological arc $\cf:[0,d]\to \free$ and
$t\in [0,d]$, let $\cf[t-a,t+b]$ be a maximal sub-arc of $\cf$
containing $\cf (t)$ and contained in one piece. Then every other
topological arc with the same endpoints as $\cf$ must contain the
points $\cf (t-a)$ and $\cf (t+b)$.
\end{quotation}
\end{lemma}

Throughout the rest of the section, $(\free, \pp)$ is a tree-graded
space.

The following statement is an immediate consequence of
\cite[Corollary 2.11]{DS}.

\begin{lemma}\label{ab}
Let $A$ and $B$ be two pieces in $\pp$. There exist a unique pair of
points $a\in A$ and $b\in B$ such that any topological arc joining
$A$ and $B$ contains $a$ and $b$. In particular $\dist (A,B)=\dist
(a,b)$.
\end{lemma}

\begin{definition}\label{sats}
For every topological arc $\g $ in $\free$ we define its
{\textit{strict saturation}}, denoted by $\Sato \g$, as
 the union of $\g$ with all the pieces intersecting $\g$ non-trivially.
\end{definition}

 \Notat\quad For every arc $\g$ in $\free$, we denote by
$I(\g)$ the collection of non-trivial sub-arcs which appear as
intersections of $\g$ with pieces from $\Sato \g$.

\medskip

The following statement immediately follows from property $(T_1)$.

\begin{lemma}\label{arcs} Every subarc in
$I(\g)$ is a maximal subarc of $\g$ contained in a piece.
\end{lemma}

\begin{definition}
The  {\textit{saturation}} of $\g$, denoted by $\Sat \g$, is the
union of $\Sato \g$ and all the pieces intersecting $\g$ by single
points outside the arcs from $I(\g)$.
\end{definition}

Obviously, $\Sato \g \subseteq \Sat \g$.

\begin{definition}\label{ctpset} Let $\g$ be a
geodesic segment, ray or line in $\free$. We denote by $\Cutp(\g)$
 the complementary set in $\g$ of the union of all the interiors of subarcs from $I(\g)$.
 We call it the
 {\textit{set of cut-points
 on}} $\g$.
\end{definition}

\begin{lemma}[Lemma
  2.23, (3) and Corollary 2.10 in \cite{DS}]\label{stronglyconvex}
If $\cf_1\sqcup \cf_2\sqcup....\sqcup \cf_k$ is a polygonal line
then $\Sato \cf_1\cup \Sato \cf_2\cup....\cup \Sato \cf_k$ is
strongly convex, i.e. it contains all topological arcs with
endpoints in it.
 \end{lemma}

Property $(T_2')$ and {Lemma \ref{stronglyconvex}} immediately imply
the following statement.

\begin{cor}\label{convsat}
Two topological arcs with the same endpoints have the same strict
saturation, the same saturation and the same set of cut-points.

In particular a topological arc joining two points in a piece is
contained in the piece.
\end{cor}

Corollary \ref{convsat} implies that we can define the following
notions.

\begin{definition}
Let $x,y$ be any two distinct points in $\free$. We define the
{\textit{strict saturation}}  of {the pair of} points $x,y$, which
we denote by $\Sato \{x,y\}$, as the common strict
 saturation of all the topological arcs joining $x$ and $y$.
 The {\textit{saturation}} of the {pair of } points $x,y$,
 $\Sat \{x,y\}$, is defined similarly.

We likewise define the {\textit{set of cut-points separating $x$ and
$y$}}, which we denote by $\Cutp\{x,y\}$, as the set of cut-points
of some (any) topological arc joining $x$ and $y$.
\end{definition}

\begin{remark} \label{211}
If an isometry fixes two points $x$ and $y$ in $\free$ then it
stabilizes $\Sato \{x,y\}$ and $\Sat \{x,y\}$, and it fixes
$\Cutp\{x,y\}$ pointwise.
\end{remark}

\begin{remark}
Note that $I(\g)$ together with the singletons not included in any
sub-arc of $I(\g)$ define  a structure of tree-graded space on $\g$
induced by the tree-graded structure of $\free$.
\end{remark}

\begin{lemma}\label{t2} For every $\epsilon>0$ let $x,y,x',y'$ be points in $X$ such that
  $\dist(x,x')<\epsilon$, $\dist(y,y')<\epsilon${. L}et $\g$ and
 $\g'$ be two geodesics connecting $x$ {with} $y$ and $x'$ {with} $y'$
 respectively, and
 let $\g_\epsilon =\g \setminus (B(x,\epsilon )  \cup B(y,\epsilon
 ))$.

 Then $\g_\epsilon\subset \Sato\g'$. In particular there exists an injective map $\iota: I(\g_\epsilon)\to
 I(\g')$ preserving the order of the arcs and such that:

\begin{itemize}
\item[(1)] $\mathfrak{a}$ and
 $\iota(\mathfrak{a})$ are in the same piece for each
$\mathfrak{a}\in I(\g_\epsilon)$;
\item[(2)] $\mathfrak{a}$ and $\iota(\mathfrak{a})$ have the same endpoints
  for all but at most two extremal intervals $\mathfrak{a}\in I(\g_\epsilon)$;
\item[(3)]  the sum of lengths of the
 sub-arcs in $I(\g)$ differs from the sum of lengths of sub-arcs in
 $I(\g')$ by at most $2\epsilon$.
\end{itemize}
 \end{lemma}

\proof Let $[x,x']$ and $[y,y']$ be two arbitrary geodesics. Let
$\bar{x}$ be the farthest from $x$ point in $\g \cap [x,x']$. The
point $\bar{y}$ is defined similarly for $y$ and $\g \cap [y,y']$.
Denote by $\bar{\g}$ the sub-arc of $\g$ between $\bar{x}$ and
$\bar{y}$.  Then $[x', \bar{x}]\cup \bar{\g} \cup [\bar{y},y']$ is a
topological arc joining $x'$ and $y'$ and containing $\g_\epsilon$.
The statements of the lemma now follow from $(T_2')$.
\endproof

\begin{lemma}\label{union} Let $U$ be a union of pieces of a
tree-graded space $(\free,\pp)$. Suppose that for every two points
$x,y$ in $U$, $\Sato \{x,y\}\subseteq U$. Then every geodesic {in
$\free$ } connecting two points from the closure $\overline{U}$ is
contained in $\overline{U}$, moreover its interior is contained in
$U$.
\end{lemma}

\proof  Let $x,y$ be two points in $\overline{U}$. Then $x$ is the
limit point of a sequence $x_n\in U$ and $y$ is the limit point of a
sequence $y_n\in U$. Let $[x_n,y_n]$ be a geodesic {in $\free$. Then
$[x_n,y_n]$ is inside $\Sato\{x_n,y_n\}\subseteq U$.} Let $[x,y]$ be
an arbitrary geodesic joining $x$ and $y$ in $\free$. By Lemma
\ref{t2}, for any $\epsilon$, $[x,y]\setminus
\left(B(x,\epsilon)\cup B(y,\epsilon)\right)$ is contained in $\Sato
\{ x_n,y_n\}\subset U$. Therefore $[x,y]$ is contained in the
closure of $U$ and its interior is contained in $U$.
\endproof

\subsection{{$\R$--tree quotients of tree-graded spaces}}

\Notat\quad For every tree-graded space $(\free,\pp)$ and every two
points $x,y$ in $\free$ let $\widetilde\dist(x,y)$ be  $\dist(x,y)$
minus the sum of lengths of sub-arcs from $I([x,y])$.

\begin{lemma}\label{well} 1. The number $\widetilde\dist(x,y)$ is well
defined (i.e. it does not depend on the choice of a geodesic
$[x,y]$).

2. The function $\widetilde\dist(x,y)$ is symmetric and satisfies
the triangular inequality.
\end{lemma}

\proof Statement 1 immediately follows from Corollary \ref{convsat}.
Statement 2 follows from Lemma \ref{stronglyconvex}.\endproof

Let us define the equivalence relation $\approx$ by $$x\approx y
\hbox{ if and only if } \widetilde\dist(x,y)=0.$$

\begin{lemma}\label{equiv}
\begin{itemize}
\item[(1)] If $x\approx y$ then
  for every geodesic $[x,y]$, its saturation is contained in the same
  $\approx$--equivalence class as $x$ and $y$. The same
  holds for every topological arc joining $x$
  and $y$. In particular,
   every piece intersecting an $\approx$--class is contained in it.
\item[(2)] The equivalence relation $\approx$ is closed.
\end{itemize}
\end{lemma}

\proof The statement in part (1) for geodesics follows immediately
from the definition of $\approx$. Together with property $(T_2')$ it
then implies the same statement for topological arcs.

(2) Let $x=\lim x_n, y=\lim y_n$ where $x_n\approx y_n$, that is
$\widetilde\dist(x_n,y_n)=0$. We need to show that
$\widetilde\dist(x,y)=0$. For every $\epsilon>0$ consider $x_n, y_n$
such that $\dist(x_n,x)<\epsilon, \dist(y_n,y)<\epsilon$. Since
$\widetilde\dist(x_n,y_n)=0$, the union of non-trivial intersections
of pieces with a geodesic $[x_n,y_n]$ is dense in the geodesic. By
Lemma \ref{t2}, the same is true for $[x,y]\setminus
(B(x,\epsilon)\cup B(y,\epsilon))$. Hence $\dist(x,y)$ minus the sum
of the lengths of sub-arcs from $I([x,y])$ cannot be bigger than
$2\epsilon$. Thus $\widetilde\dist(x,y)<2\epsilon$ for every
$\epsilon>0$, hence $\widetilde\dist(x,y)=0$.\endproof

\begin{remark} \label{rem7} Note that $\approx$ is in general not the smallest equivalence relation
satisfying the properties in Lemma \ref{equiv}. For example,
consider the collection of closures of subintervals of the interval
$[0,1]$ used {in } creating the Cantor set (the middle thirds). This
collection together with the singletons not contained in any
interval form a tree-graded structure $\pp$ on the unit interval. It
is easy to see that in this case $\approx$ has only one equivalence
class. On the other hand the smallest equivalence relation with the
properties listed in Lemma \ref{equiv} has as equivalence classes
all the pieces $\pp$.
\end{remark}

 \begin{lemma}\label{equiv1} Every $\approx$--class is a
connected  union of pieces.
\end{lemma}

\proof  It follows immediately from Lemma \ref{equiv}.
\endproof

Recall \cite{DS} that for every point $x$ in a tree-graded space
$\free$, the {\em transversal tree at $x$}  consists of all
$y\in\free$ such that any geodesic $[x,y]$ has only trivial
intersections with pieces. We proved in \cite{DS} that transversal
trees either coincide or do not intersect. {Any two points in the
same transversal tree are joined by a unique geodesic in $\free$.
These geodesics are called {\em transversal geodesics} in $\free$.}
The following lemma immediately follows from the definition of
$\approx$.

\begin{lemma}\label{transv} If $x\ne y$ are in the same transversal tree of $\free$
then $\dist(x,y)={\widetilde\dist}(x,y)$. In particular,
$x\not\approx y$. Thus every transversal tree projects into
$\free/{\approx}$ isometrically.
\end{lemma}

Let $T$ be the quotient $\free/{\approx}$.

\begin{lemma}\label{Ttree} $T$ is an $\R$--tree with respect to the metric induced by
 $\widetilde\dist$. Every geodesic in $\free$ projects onto a geodesic in $T$.
\end{lemma}

\proof We first show that $T$ is a geodesic metric space. Consider
two points $\bar{x}$ and $\bar{y}$ in $T$. We shall construct one
geodesic $\g_{\bar{x},\bar{y}}$ joining them. Let $x, y\in \free$ be
representatives of $\bar x$ and $\bar y$, and consider a geodesic
$[x,y]$. Let $\mu$ be the Lebesgue measure on $[x,y]$ (here the
geodesic is identified with an interval). Define a new measure
$\mu_0$ on all Borel sets by
$$
\mu_0 (B)= \mu \left(B\setminus \bigcup_{\mathfrak{a}\in
  I([x,y])}\mathfrak{a}\right)\, .
$$

Note that $\mu_0$ is absolutely continuous with respect to $\mu$,
hence there exists a measurable function $f:[x,y]\to \R$ such that
$$
\mu_0(B)=\int_B f\, d\mu\, .
$$

Let $D=\dist (x,y)$ and let $\delta
=\widetilde\dist(\bar{x},\bar{y})$. For every $t\in [0,D]$ let $x_t$
be the point on $[x,y]$ at distance $t$ from $x$ and let $[x,x_t]$
be the sub-arc of $[x,y]$ between $x$ and $x_t$. The function
$F:[0,D]\to [0,\delta]$, $F(t)=\int_{[x,x_t]} f\, d\mu$ is monotone
non-decreasing and continuous. Define the map $\g_{\bar{x},\bar{y}}
: [0, \delta ] \to T$ such that $\g_{\bar{x},\bar{y}}(s)$ is the
projection $\bar{x}_{t(s)}$, where $x_{t(s)}$ is a point in $F\iv
(s)$.

 Let $s<r$ be two numbers in $[0, \delta ]$. In
order to compute $\widetilde\dist (\g_{\bar{x},\bar{y}} (s),
\g_{\bar{x},\bar{y}}(r))$ consider the sub-arc
 $\left[x_{t(s)},x_{t(r)} \right]$ of $[x,y]$. We have the following
equalities
$$
\widetilde\dist (\g_{\bar{x},\bar{y}} (s),
\g_{\bar{x},\bar{y}}(r))=\mu_0 \left(\left[x_{t(s)},x_{t(r)} \right]
\right)=\int_{\left[x_{t(s)},x_{t(r)} \right]}f\, d\mu F\left(x_{t(r)} \right)-F\left(x_{t(s)}\right)=r-s\, .
$$

Note that for every $m\in [x,y]$, its projection $\bar{m}$ onto $T$
is in $\g_{\bar{x},\bar{y}}$. Indeed, let $d=\widetilde\dist
(\bar{x},\bar{m})$ and let $x_{t(d)}$ be chosen as before. Then
$d=\int_{[x,m]} f\, d\mu=\int_{[x,x_{t(d)}]} f\, d\mu$, hence
$$\int_{[m,x_{t(d)}]} f\, d\mu
    = \widetilde{\dist} (\bar{m}\, ,\, \bar{x}_{t(d)})=0
$$ and
      $\g_{\bar{x},\bar{y}} (d)=\bar{x}_{t(d)}=\bar{m}$.

Now we show that for every two points $\bar{x},\bar{y}$ there exists
a unique geodesic joining them in $T$. Let $\pgot$ be an arbitrary
geodesic joining $\bar{x}$ and $\bar{y}$ in $T$. Let
$\g_{\bar{x},\bar{y}}$ be the geodesic constructed above between
$\bar{x},\bar{y}$. Let $\bar{z}$ be an arbitrary point on $\pgot$
and let $z$ be a representative of $\bar{z}$ in $\free$. For two
arbitrary geodesics $[x,z],[z,y]$ let $z'$ be the farthest from $z$
point in $[x,z]\cap [z,y]$. Then $[x,z']\cup [z',y]$ is a
topological arc. Consider the maximal sub-arc $[a,b]$ in it
containing $z'$ and contained in a piece. Note that this sub-arc can
be a point. Property $(T_2')$ implies that $\{ a,b\} \subset [x,y]$.
By definition $\bar{a}=\bar{b}=\bar{z'}$, which is a point on
$\g_{\bar{x},\bar{y}}$. Consequently, $\widetilde\dist
(\bar{x},\bar{z'})+\widetilde\dist
(\bar{y},\bar{z'})=\widetilde\dist (\bar{x},\bar{y}).$ On the other
hand, since $z'\in [x,z]\cap [z,y]$, $\widetilde\dist
(\bar{x},\bar{z'})\leq \widetilde\dist (\bar{x},\bar{z})$ and
$\widetilde\dist (\bar{y},\bar{z'})\leq \widetilde\dist
(\bar{y},\bar{z})$. It follows that the previous two inequalities
are in fact equalities, and $\bar{z'}=\bar{z}$.\endproof

\begin{lemma} \label{inv} Let $\free$ be a tree-graded space.
\begin{itemize}
\item[(1)] An isometry $\phi$ of $\free$ permuting the pieces induces an isometry
$\tilde \phi$ of the real tree $T$.
\item[(2)] For every non-trivial geodesic $\g$ in $T$ there exists a non-trivial geodesic $\pgot$
 in $\free$ such that its projection on $T$ is $\g${.
 Moreover}
 the isometry $\tilde \phi$ fixes $\g$ pointwise if and only
 if $\phi$ fixes $\Cutp(\pgot)$ pointwise.

In the particular case when $\g$ is the projection of a geodesic
$\g_0$ in a transversal tree, $\pgot$ can be taken equal to $\g_0$.
\end{itemize}
\end{lemma}

\proof  (1) Since for every piece $A$ and every geodesic $\g$, the
lengths of $\phi(A)\cap \phi(\g)$ and $A\cap \g$ are the same,
$\widetilde\dist(\g_-,\g_+)=\widetilde\dist(\phi(\g_-),\phi(\g_+))$.

(2) Consider the endpoints $a\ne b$ of $\g$ in $T$. By Lemma
\ref{Ttree} if $x,y\in \free$ are representatives of $a$ and $b$
respectively, and $[x,y]$ is a geodesic joining them, the latter
 projects onto $\g$.  The intersection of
the $\approx$--equivalence class containing $x$ with $[x,y]$ is
closed and connected, by Lemma \ref{equiv}. Hence it is a geodesic
segment $[x,x']\subsetneq [x,y]$. By replacing if necessary $x$ by
$x'$, we may therefore assume that $[x,y]$ contains only one point
in the $\approx$--equivalence class of $x$. A similar argument
%%%%%%12/12/05%%%
allows to assume that the intersection of $[x,y]$ and the
$\approx$--class containing $y$ is $\{y\}$. The same can then be
said about $\phi[x,y]$ and its endpoints. We take $\pgot =[x,y]$.
Since the projection of $\Cutp(\pgot)$ is $\g$, if $\phi$ fixes
$\Cutp(\pgot)$ pointwise then $\tilde \phi$ fixes $\g$ pointwise.
Now we show that the converse also holds.

As $\tilde \phi (a)=a$ and $\tilde \phi (b)=b$, it follows that
$\phi(x)\approx x$ and $\phi(y)\approx y$. Any geodesic $[x, \phi
(x)]$ intersects $\phi(\pgot )$ only in
  $\phi(x)$, likewise any geodesic $[y,\phi(y)]$ intersects $\phi(\pgot )$ in
  $\phi(y)$. Then
  $\g_1=[x, \phi (x)]\sqcup \phi(\pgot )\sqcup [\phi (y), y]$ is a topological arc and
  $\phi(x),
  \phi (y)$ belong to $\Cutp(\g_1)$.
It follows from $(T_2')$ that the geodesic $\pgot $ must contain
$\phi(x) ,
  \phi (y)$. On the other hand $\dist (x,y)=\dist (\phi(x) ,
  \phi (y))$. Consequently $x=\phi (x)$ and $y=\phi (y)$. This and Remark \ref{211} imply that
  $\Cutp(\phi(\pgot))=\phi(\Cutp(\pgot))=\Cutp(\pgot)$. Moreover, since $\phi$ is an
   isometry it must fix $\Cutp(\pgot)$ pointwise.\endproof

\subsection{Nested tree-graded structures}\label{nested}

If $\pp$ and $\pp'$ are two  collections of subsets of a set
$\free$, we write $\pp \prec \pp'$ if for every set $A\in \pp$
there exists $A'\in \pp'$ such that $A\subset A'$. The relation
$\prec$ induces a partial order on the set of tree-graded
structures of a tree-graded space $\free$ because of our
convention that in a tree-graded structure, pieces cannot contain
each other.

For every tree-graded space $(\free,\pp)$, consider the following
equivalence relation $\sim$. Two points $x$ and $y$ are
$\sim$--equivalent if one (hence any by Corollary \ref{convsat})
geodesic $[x,y]$ is inside the union of a finite number of pieces.
This relation is transitive, by Lemma \ref{stronglyconvex}.

\medskip

Let $\kappa$ be an equivalence class for $\sim$. The following lemma
is obvious.

\begin{lemma}\label{geod}
For any two points $x,y$ in $\kappa$, the saturation $\Sat\{ x,y\}$
is contained in $\kappa$.
\end{lemma}

\begin{lemma}\label{k}
The union of strict saturations of pairs of points in $\kappa$ is
equal to $\kappa$.
\end{lemma}

\proof By Lemma \ref{stronglyconvex}, the union of strict
saturations of pairs of points from $\kappa$ is contained in
$\kappa$. On the other hand $\kappa$ is the union of all geodesics
with endpoints in $\kappa$. Since every geodesics is inside its
strict saturation, we conclude that $\kappa$ is inside the union of
strict saturations of geodesics with endpoints in $\kappa$.
\endproof

Let $\pp'$ be the collection of closures of $\sim$--equivalence
classes. For every $A$ in $\pp'$ let $\kappa(A)$ be the
$\sim$--equivalence class such that $A=\overline{\kappa(A)}$.

\begin{lemma}\label{ppp} $\free$ is tree-graded with respect to $\pp'$.
\end{lemma}

\proof By construction all pieces in $\pp'$ are closed. By Lemmas
\ref{k} and \ref{union}, all pieces are geodesic subspaces.

Clearly every piece of $\pp$ is inside a piece of $\pp'$. Therefore
pieces of $\pp'$ cover $\free$ and $(T_2)$ is satisfied.

Assume that the intersection of two sets $A$ and $B$ from $\pp'$
contains two points $a,b$. The set $A$ is the closure of a
$\sim$--equivalence class $\kappa(A)$ and $B$ is the closure of
$\kappa(B)$.

By Lemma \ref{union}, the interior of a geodesic $[a,b]$ is
contained in $\kappa(A)\cap \kappa(B)$, a contradiction (distinct
equivalence classes do not intersect).
\endproof

\begin{lemma}\label{seqcut}
Let $A$ be a piece in $\pp'$ which is the closure of an
equivalence class $\kappa$ of $\sim$, let $x\in \kappa$ and let
$p\in A\setminus \kappa $. On every geodesic $[x,p]$, there exists
an infinite sequence of pairwise distinct points
$x_0=x,x_1,x_2,...., x_n,...$ appearing in this order from $x$ to
$p$, such that $x_n$ converges to $p$ and $[x_n,x_{n+1}]$ is the
intersection of $[x,p]$ with a piece in $\pp$.

Moreover all geodesics joining $x$ and $p$ contain this ordered
countable set of points.
\end{lemma}

\proof Let $(y_n)$ be a sequence of points from $\kappa$ converging
to $p$. By the definition of $\sim$, any geodesic $[x,y_n]$ is
covered by finitely many pieces from $\pp$. By replacing, if
necessary, $y_n$ by the farthest point from it in $[x,y_n]\cap
[y_n,p]$ one can assume that $[x,y_n]$ is contained in a topological
arc joining $x$ to $p$. Property $(T_2')$ implies that all endpoints
of intersections of $[x,y_n]$ with pieces of $\pp$ are also on
$[x,p]$. Let $c_n$ be the nearest to $y_n$ such a point. Then $c_n$
also converges to $p$. Indeed, otherwise, all but finitely many
$c_n$ are equal to a point $c\in [x,p]$, and $y_n$ must be in the
same piece $P_n$ as $c$. If all but finitely many pieces $P_n$
intersect $[c,p]$ nontrivially, then these pieces are the same, and
$p\in P_n$, a contradiction. If infinitely many pieces $P_n$
intersect $[c,p]$ by a point (which must by $c$) then $c$ is the
projection of $p$ onto $P_n$ for infinitely many $n$'s. Thus $\dist
(y_n,p)\geq \dist (c,p)$ for infinitely many $n$'s, therefore
$c=p\in \kappa$, a contradiction.

We have, therefore, found a sequence $c_n\in \kappa\cap [x,p]$ of
endpoints of non-trivial intersections of $\pp$--pieces with
$[x,p]$ that converges to $p$, as required. The ``moreover"
statement of the lemma follows from $(T_2')$.
\endproof

\begin{lemma}\label{stabk}
Let $\psi$ be an isometry of $\free$ permuting pieces in $\pp$ and
such that $\psi(A)=B$, where $A,B$ are pieces in $\pp'$. Then
$\psi(\kappa(A))=\kappa(B)$.
\end{lemma}

\proof If $\kappa(A)$ is a singleton then the result is obvious.
Thus we suppose that $\kappa(A)$ is not a singleton. Let $x$ be an
arbitrary point in $\kappa(A)$, let $y$ be a point in $\kappa(A)
\setminus \{x\}$ and let $\g$ be a geodesic joining $x$ with $y$.
By the definition of $\sim$, this geodesic is covered by finitely
many non-trivial intersections with pieces contained in
$\kappa(A)$. Then $\psi(\g)$ is a geodesic joining $\psi(x)$ with
$\psi(y)$, and it is also covered by finitely many non-trivial
intersections with pieces from $\pp$. Since $\psi(x)$ and
$\psi(y)$ are two points in $B=\overline{\kappa(B)}$, by Lemma
\ref{union} the interior of $\psi(\g)$ is contained in
$\kappa(B)$. It follows that all pieces intersecting $\psi(\g)$
 non-trivially are contained in $\kappa(B)$,
hence $\psi(\g) \subset \kappa(B)$. In particular $\psi(x),
\psi(y) \in \kappa(B)$. Since $x$ was arbitrary in $\kappa(A)$, we
conclude that $\psi(\kappa(A) )\subseteq \kappa(B)$. A similar
argument applied to $\psi\iv$ implies that $\psi(\kappa(B)
)\subseteq \kappa(B)$.\endproof

\begin{lemma}\label{stabkp}
Let $\psi$ be an isometry of $\free$ permuting pieces in $\pp$ and
let $A$ be a piece in $\pp'$ such that $\psi(A)=A$. Let $x$ be an
arbitrary point in $\kappa=\kappa(A)$, $p$ a point in $
A\setminus\kappa$ and $x_1,x_2,...$ the sequence of points on
$[x,p)$ converging to $p$ as in Lemma \ref{seqcut}. If $\psi(p)=p$
then $\psi$ fixes all but finitely many of the $x_i$'s. In
particular $\psi$ stabilizes all but finitely many
  pieces in $\pp$ intersecting $[x,p]$ non-trivially.
\end{lemma}

\proof By Lemma \ref{stabk}, $\psi(x) \in \kappa$, and
$\psi(x_0),\psi(x_1),..., \psi(x_n),... $ is a sequence of points
on $[\psi(x),p]$ which converges to $p$ and consists of endpoints
of arcs from $I([\psi(x),p])$.
 A geodesic $\g=[x,\psi(x)]$ is covered by finitely many pieces
 since $x\sim \psi(x)$.
 Let $x'$ be the farthest from $\psi(x)$ point in $\g \cap [\psi(x), p]$. It
 is different from $p$, and $[x,x']\sqcup [x',p]$ is a topological arc. Property
 $(T_2')$ implies that this arc contains all $x_i$, and since $[x,x']$
  is covered by finitely many pieces, for some $n_0$, all $x_n$ with $n\geq n_0$
  are in $[x',p]$. By its definition, this sequence coincides with
 the intersection of
  the sequence $(\psi(x_n))$ with $[x',p]$. Thus we have that
   $x_n=\psi(x_{n+k})$ for all $n\geq n_0$ and some fixed $k\ge 0$. On the
   other hand, $\dist(x_n,p)=\dist(\psi(x_{n+k}),p)=\dist(x_{n+k},p)$ implies
   that $k=0$. Therefore $\psi$ fixes all $x_n$ with $n\geq n_0$.
 \endproof

\begin{lemma}\label{union2} Consider a sequence $\pp_1\prec\pp_2\prec...$
  of tree-graded structures on $\free$ and an ascending sequence
  of pieces $A_1\subseteq A_2\subseteq...$ where
$A_i\in \pp_i$. The closure $\widehat{A}$ of the union $\bigcup
A_i$ contains together with any two points any geodesic joining
them. Moreover the interior of such a geodesic is contained in
$\bigcup A_i$ and a non-trivial sub-arc of it is contained in all
but finitely many $A_i$.
\end{lemma}

\proof Let $x,y$ be two arbitrary points in $\widehat{A}$ and let
$\g$ be an arbitrary geodesic with endpoints $x,y$. We have
$x=\lim_{n\to \infty }x_n$ and $y=\lim_{n\to \infty }y_n$ where
$x_n, y_n \in \bigcup A_i$. Without loss of generality we may
suppose that both $x_n$ and $y_n$ are in some piece $A_{i_n}$.
Lemma \ref{t2} implies that $\g \setminus \nn_\epsilon
(\{x,y\})\subset A_{i_n}$ for $n$ large enough. Thus the interior
of $\g$ is contained in $\bigcup A_i$ and $\g$ is contained in
$\widehat{A}$.\endproof

\begin{lemma}\label{up} For every sequence $\pp_1\prec\pp_2\prec...$
  of tree-graded structures on $\free$ there exists the smallest
  tree-graded structure $\bigcup\pp_i$ such that $\pp_j\prec \bigcup \pp_i$
  for every $j$.

  The pieces in $\bigcup \pp_i$ are closures of unions
of ascending sequences $A_1\subseteq A_2\subseteq...$ where $A_i\in
\pp_i$.
\end{lemma}

\proof Consider all possible sequences $A_1\subseteq
A_2\subseteq...$ where $A_i\in \pp_i$. Let a collection
$\widehat\pp$ of subsets of $\free$ consist of closures of unions of
all these sequences of sets. It is clear that
$\pp_j\prec\widehat\pp$ for every $j$.

Let us prove that $\free$ is tree-graded with respect to
$\widehat\pp$. It is obvious that all pieces in $\widehat\pp$ are
closed and that $(T_2)$ is satisfied.

The fact that every piece in $\widehat\pp$ is a geodesic subspace
follows from Lemma \ref{union2}.

Let us prove $(T_1)$. Let $A=\overline{\bigcup A_i}$,
$B=\overline{\bigcup B_i}$ be two pieces in $\widehat\pp$ and $x\ne
y\in A\cap B$. Then every geodesic $[x,y]$ is in $A\cap B$. By Lemma
\ref{union2}, a non-trivial portion $[x',y']$ of a geodesic $[x,y]$
belongs to all but finitely many $A_i$ and to all but finitely many
$B_i$. Since $\pp_i$ satisfy $(T_1)$, $A_i=B_i$ for all but finitely
many $i$. Hence $A=B$.

It remains to show that every tree-graded structure
$\widetilde\pp\succ\pp_i$ for all $i$ satisfies
$\widetilde\pp\succ\widehat\pp$. Indeed, let $A\in \widehat\pp$.
Then $A=\overline{\bigcup A_i}$. If all $A_i$ are points then $A$ is
a point and $A$ is contained in a piece from $\widetilde\pp$.
Otherwise, all but finitely many $A_i$'s are not points and by
$(T_1)$ applied to $\widetilde\pp$, all $A_i$ are in the same piece
$\tilde A\in \widetilde\pp$. Since $\tilde A$ is closed, $A\subseteq
\tilde A$.
\endproof

The following lemma is obvious.

\begin{lemma} \label{perm} 1. Let $(\free,\pp)$ be a tree-graded space
and let $\phi$ be an isometry of $\free$ permuting pieces of $\pp$.
Then $\phi$
  permutes pieces of $\pp'$.

2. If $(\free,\pp_i)$ is a sequence of tree-graded structures on
$\free$,
   $\pp_i\prec\pp_{i+1}$ and an isometry $\phi$ permutes pieces of
   each $\pp_i$ then $\phi$ permutes pieces of $\bigcup\pp_i$.
\end{lemma}

\begin{lemma}\label{stabu}
Let $\pp_1\prec\pp_2\prec...$ be a sequence of tree-graded
structures on a complete geodesic metric space $\free$.

Let $\psi$ be an isometry of $\free$ which permutes the pieces in
$\pp_i$ for any $i$, and let $A,B$ be pieces in $\bigcup \pp_i$ such
that $\psi(A)=B$. Suppose that A and B are closures of unions
$\bigcup A_i$ and respectively $\bigcup B_i$, with $A_i, B_i\in
\pp_i$. Then $\psi(A_i )=B_i$ for $i$ large enough.
\end{lemma}

\proof Suppose that $\bigcup A_i$ is not a singleton, otherwise the
statement clearly holds. Let $x$ be a point in $\bigcup A_i$. For
some $i$ large enough, $x$ is in $A_i$ together with a non-trivial
geodesic $[x,y]$. The image of this geodesic $\psi([x,y])$ is in
$\psi(A_i)\in \pp_i$ and it is also in $B$. Lemma \ref{union2}
implies that the interior of $\psi([x,y])$ is in $\bigcup B_i$ and
that a non-trivial subarc $\g$ of $\psi([x,y])$ is contained in some
$B_j$.

Let $k=\max (i,j)$. Then $\psi(A_i)\subset \psi(A_k)$ and
$B_j\subset B_k$. Hence $\psi(A_k)$ and $B_k$ are two pieces in
$\pp_k$ which have in common $\g$. Property $(T_1)$ implies that
$\psi(A_k)=B_k$ and the same holds for every $m\geq k$.\endproof

\subsection{Groups acting on $\R$--trees}
\label{gat}

An action of a group $G$ on an $\R$--tree $T$ is called {\em
stable} if the set of stabilizers of arcs of $T$ satisfies the
ascending chain condition (ACC). An action is called {\em minimal}
if $T$ does not have a proper invariant subtree.

We begin by recalling two known result of Rips, Bestvina and Feighn,
and Sela on non-trivial stable actions on trees.

\begin{theorem}[Rips-Bestvina-Feighn {\cite[Theorem 9.5]{BF}}]\label{bf}
Let $\Lambda$ be a finitely presented group with a
  non-trivial stable minimal action by isometries on an
  $\R$--tree $T$.
  Then one of the following two cases occurs:
  \begin{itemize}
  \item[(1)] $\Lambda$ splits over an extension $E$-by-cyclic,
  where $E$ is the stabilizer of a non-trivial arc of $T$;
  \item[(2)] $T$ is a line and $\Lambda$ has a subgroup of index at most 2
  that is the extension of the kernel of the action of $\Lambda$ on $T$ by a finitely
  generated free Abelian group.
  \end{itemize}
\end{theorem}

\begin{theorem}[Sela {\cite[Section 3]{Se}}] \label{sela} Let $\Lambda$ be a finitely
generated group with a
  non-trivial stable minimal action by isometries on an
  $\R$--tree $T$ and assume that the stabilizers of all tripods in
  $T$ are trivial. Then the conclusion of Theorem \ref{bf} holds.
\end{theorem}

The following version of Theorem \ref{sela} is proved by V.
Guirardel in \cite{Gui}. Since the proof is still not published, we
present this stronger version along with Theorem \ref{sela}.

\begin{definition} \label{defgui} The \textit{height} of an arc in an $\R$--tree with respect to
the action of some
 group $G$ on it is the maximal length of a decreasing chain of sub-arcs with distinct stabilizers.
 If the height of an arc is zero then it follows that all sub-arcs of it have the same stabilizer.
  In this case the arc is called \textit{stable}.

  The tree $T$ is \textit{of finite height} if any arc of
  it can be covered by finitely many arcs with finite height.
  If the action is minimal and $G$ is finitely generated
  then this condition is equivalent to the
  fact that there exists a finite collection of arcs
  $\mathcal{I}$ of finite height such that any arc is covered by finitely many translates of arcs in $\mathcal{I}$ \cite{Gui}.
\end{definition}

\begin{theorem}[{Guirardel} \cite{Gui}]\label{gui}
Let $\Lambda$ be a finitely generated group and let $T$ be a real
tree on which $\Lambda$ acts minimally and with finite height.
Suppose that the stabilizer of any non-stable arc in $T$ is finitely
generated.

Then one of the following three situations occurs:

\begin{itemize}
    \item[(1)] $\Lambda$ splits over the stabilizer of a
non-stable arc or over the stabilizer of a tripod;

    \item[(2)] $\Lambda$ splits over a virtually
cyclic extension of the stabilizer of a stable arc;

\item[(3)] $T$ is a line and $\Lambda$ has a subgroup of index at most
2 that is the extension of the kernel of that action by a finitely
  generated free Abelian group.
\end{itemize}
\end{theorem}

  \medskip

Stability of the action on an $\R$--tree is a necessary condition
and cannot be removed from Theorem \ref{bf} or \ref{sela} (see
Dunwoody \cite{Dun}). The next lemma shows that in some cases
stability and finite height follow from the algebraic structure of
stabilizers of arcs.

\begin{lemma}\label{stable} Let $G$ be a finitely generated group
acting on an $\R$--tree $T$ with finite of size at most $D$ tripod
stabilizers, and (finite of size at most $D$)-by-Abelian arc
stabilizers, for some constant $D$. Then
\begin{itemize}
    \item[(1)] an arc with stabilizer of size $>(D+1)!$ is stable;
    \item[(2)] every arc of $T$ is of finite height (and so the action is
of finite height and stable).
\end{itemize}
\end{lemma}

\proof (1) Let $H$ be the stabilizer of an arc $\g$ in $T$,
$|H|>(D+1)!$ and let $\g_1$ be a sub-arc in $\g$ with stabilizer
$H_1>H$. Then $H_1$ is an extension of a subgroup $U$ of size at
most $D$ by an Abelian group, and the centralizer $C(U)$ of $U$ in
$H_1$ has index at most $D!$ since $H_1/C(U)$ is embedded into
$\Aut(U)$. Therefore $|C(U)\cap H|>D$.

For every $h\in H_1\setminus H$, $hHh\iv$ fixes $h\g$. Since
$h\not\in H$, $h\g\ne\g$ but $\g_1\subseteq \g\cap h\g$. Hence $\g
\cup h\g$ contains at least one tripod.  Thus the group $H \cap hH
h\iv$, which stabilizes $\g \cup h\g$, is of size at most $D$. If
$h\in U$ then $H \cap hH h\iv$ contains $C(U)\cap H$. This gives
$U\le H$. Since $H_1/U$ is Abelian, $H$ contains the derived
subgroup of $H_1$. Hence $H$ is normal in $H_1$. Therefore $H \cap
hH h\iv=H$, and $|H|\le D$, a contradiction.

(2) Indeed, from (1) it immediately follows that the height of every
arc in $T$ cannot exceed $(D+1)!+1$.\endproof

We need a pretree version of Levitt's theorem \cite[Theorem 1]{Lev}.

\begin{definition}[see \cite{Bow}, page 10] \label{pretreed} A {\em pretree} is a set
equipped with a ternary {\em betweenness} relation $xyz$ satisfying
the following conditions:
\begin{itemize}
\item{(PT0)} $(\forall x,y) (\neg xyx)$.

\item{(PT1)} $xzy \Leftrightarrow yzx$.

\item{(PT2)} $(\forall x,y,z) (\neg (xyz\wedge xzy))$.

\item{(PT3)} $xzy$ and $z\neq w$ then $(xzw \vee yzw)$.
\end{itemize}
\end{definition}

An \textit{interval} in a pretree is a set $\lfloor x,y\rfloor$
composed of all $z$ such that $xzy$ holds.

\begin{definition}\label{nn}
An automorphism $g$ of a pretree $T$ is called {\em
  non-nesting} if for every interval $I$, $g\cdot I\subseteq I$
  implies $g\cdot I=I$.
\end{definition}

\begin{theorem}\label{lev} If a finitely presented group $G$ admits a
non-trivial non-nesting action by pretree automorphisms on an
$\R$--tree $T$, then it admits a non-trivial isometric action on
some complete $\R$--tree $T'$ with the following properties:

\begin{itemize}
\item[(1)] the stabilizer  of an arc in $T'$ is also the
stabilizer of an arc in $T$;

\item[(2)] the stabilizer of a tripod in $T'$ is also the
stabilizer of a tripod in $T$;

\item[(3)] if $G$ stabilizes a line in $T'$ then it stabilizes a
line in $T$.
\end{itemize}
\end{theorem}

\proof Property (1) was proved in \cite{Lev} under the assumption
that $G$ acts on $T$ by homeomorphisms. On the other hand, it is
proved in Mayer and Oversteegen \cite{MO} that for every action of
a finitely generated group $G$ on an $\R$--tree $T$ by pretree
automorphisms, one can modify the metric on $T$ (preserving the
pretree structure) so that $G$ acts on $T$ by homeomorphisms. Thus
we can apply the results of \cite{MO}, and then Levitt's theorem.

Another way of proving property (1) is the following. It is easy
to see that every pretree automorphism of an $\R$--tree preserves
geodesic intervals, hence restricted to a finite subtree it
becomes a homeomorphism. Since Levitt's proof only uses
restrictions of the homeomorphisms to finite subtrees, it carries
without any change.

Property (2) follows from property (1).

Property (3) easily follows from the proof of Levitt's theorem in
\cite{Lev}.
\endproof

We are going to prove a version of Theorem \ref{lev} for finitely
generated groups. We start with the following Lemma. For the
definition of the notion of $\omega$--limit used in it see Section
\ref{ac}.

\begin{lemma} \label{indlim} Let $G=\la S\ra$ be an inductive
limit of groups $G_n=\la S\ra$ and surjective homomorphisms
$G_1\to G_2\to...$ that are identical on $S$. For every $n$, let
$(T_n, dist_n)$ be a complete $\R$--tree upon which $G_n$ acts
non-trivially by isometries. For every $x\in T_n$ let
$$d_n(x)=\sup_{a\in S}\dist(ax,x),$$ and let $x_n$ be a point in
$T_n$ such that $d_n=d_n(x_n)\le \inf_{x\in T_n}d_n(x)+1$. Let
$\omega$ be any ultrafilter. Then

\begin{itemize}
\item[(1)] $G$ acts non-trivially by isometries on the $\omega$--limit $T$ of
$(T_n,\dist_n/d_n, (x_n))$ by $$(g_n)\lim^\omega (y_n)=\lim^\omega
(g_ny_n)\; ;$$

\item[(2)] for every arc $l$ in $T$ with stabilizer $\Stab_G(l)$
there exists a sequence of arcs $l_n$ in $T_n$
 such that $\lio{l_n}\subset l$ and such that any finitely generated subgroup $K$ in
$[\Stab_G(l),\Stab_G(l)]$ is inside the inductive limit of
stabilizers $\Stab_{G_n} (l_n)$ for $n\in I_K\subseteq \mathbb{N}$,
where $\omega(I_K)=1$;

\item[(3)] for every tripod  $abc$ in $T$ there exists a sequence of tripods
$\alpha_n\beta_n\gamma_n$ in $T_n$ such that
$\lio{\alpha_n\beta_n\gamma_n}\subset abc$ and such that the
following holds. Any finitely generated subgroup $L$ stabilizing
$abc$ in $G$ is inside the inductive limit of stabilizers of
$\alpha_n \beta_n \gamma_n$ in $G_n$ for $n\in I_L\subseteq
\mathbb{N}$, where $\omega(I_L)=1$.
\end{itemize}
\end{lemma}

\proof (1) We need to prove that for every $g=(g_n)\in G$ and
every $\lim^\omega(y_n)\in T$ the point $\lim^\omega(g_ny_n)$ is
defined (i.e. the distance from it to $(x_n)$ is not infinity),
and that the action $(g_n)\lim^\omega(y_n)=\lim^\omega(g_ny_n)$ is
well defined (i.e. it does not depend on the sequence representing
$g$). The first statement follows from the choice of $d_n$ and
$x_n$, as $\dist (g_ny_n,x_n)\leq \dist (y_n,x_n) + |g|_S d_n$.
The second statement follows from the fact that $\omega$ is a
non-principal ultrafilter and so every set of natural numbers with
finite complement has $\omega$--measure $1$.

\me

(2) Since $T$ is a tree and so any pair of points is connected by
a unique arc, every arc $l$ in $T$ is the $\omega$--limit of arcs
$\ell^n\subseteq T_n$. Since the arc $l$ has non-zero length,
$|\ell^n|=O(d_n)$. Suppose that $(g_n)$ fixes $l$. Then
\begin{equation}\label{tl}
\lim_\omega \frac{\dist(g_n\ell^n_-, \ell^n_-)}{d_n}=\lim_\omega
\frac{\dist(g_n\ell^n_+, \ell^n_+)}{d_n}=0\, .
\end{equation}

Since $T_n$ is a tree, (\ref{tl}) can only happen if $\omega$-a.s.
there exists a number $r_n=o(d_n)$ and a number $\epsilon_n=o(d_n)$
such that for every point $x_n\in \ell^n\setminus
\nn_{\epsilon_n}(\{\ell^n_-,\ell^n_+\})$,
$$\dist(g_nx_n,x_n)=r_n.$$ This immediately implies that for any fixed $\varepsilon >0$,
any two elements $(g_n)$ and $(h_n)$ from $G$ stabilizing $l$,
$[g_n,h_n]$ fixes the arc $l_n=\ell^n\setminus \nn_{\varepsilon
d_n}(\{\ell^n_-,\ell^n_+\})$ for $n\in I_{g,h}\subseteq \N$ where
$\omega(I_{g,h})=1$. Hence $[g,h]$ is in the induction limit of
stabilizers $\Stab_{G_n}(l_n)$, $n\in I_{g,h}$. This implies Part
(2).

\me

(3) The proof is similar to the proof of (2). A tripod $abc$ in
$T$ is the $\omega$--limit of tripods $a_nb_nc_n$ in $T_n$. If
$(g_n)$ fixes the tripod $abc$, then $g_n$ must move the ends
$a_n, b_n, c_n$ of the tripod $a_nb_nc_n$ by distance $o(d_n)$
$\omega$--a.s. This implies that $g_n$ must fix the center mass
$m_n$ of $a_nb_nc_n$ $\omega$--a.s. Therefore for any
$\varepsilon>0$, the element $g_n$ must fix the tripod
$\alpha_n\beta_n\gamma_n$ where $\alpha_n\in [a_n,m_n], \beta_n\in
[b_n,m_n], \gamma_n\in [c_n,m_n]$, and
$\dist(a_n,\alpha_n)=\dist(b_n,\beta_n)= \dist(c_n,
\gamma_n)=\varepsilon d_n$ $\omega$--a.s.
\endproof

\begin{theorem}\label{lev1} If a finitely generated group $G$ admits a
non-trivial non-nesting action by pretree automorphisms on an
$\R$--tree $T$, then it admits a non-trivial isometric action on
some complete $\R$--tree $T'$ and

\begin{itemize}
\item[(1)] the derived subgroup of a stabilizer of an arc in $T'$
is locally inside the stabilizer of an arc in $T$;

\item[(2)] the stabilizer of a tripod in $T'$ is locally inside the
stabilizer of a tripod in $T$.

\end{itemize}
\end{theorem}
%45362043

\proof If $G$ is finitely presented then we can apply Theorem
\ref{lev}. If $G$ is not finitely presented, then we can represent
$G$ as the inductive limit of a sequence of finitely presented
groups and surjective homomorphisms $G_1\to G_2\to...$. Each $G_n$
acts on $T$ by pretree automorphisms, and we can apply Theorem
\ref{lev} to this action. Hence each $G_n$ acts by isometries on a
complete $\R$--tree $T_n$ and the stabilizers of arcs (tripods) of
$T_n$ in $G_n$ are stabilizers of arcs (tripods) of $T$ in $G_n$.
By Lemma \ref{indlim}, $G$ acts on a complete $\R$--tree $T'$ by
isometries and properties (1), (2), (3) of the lemma hold.

Let $l$ be a non-trivial arc in $T'$. Consider  the stabilizer
$\sss$ of $l$ in $G$. Then by Lemma \ref{indlim} any finitely
generated subgroup of $[\sss,\sss]$ is inside the inductive limit of
stabilizers $\sss_n, n\in I\subseteq \mathbb{N}$, of arcs $l^n$ in
$T_n$, $\omega(I)=1$. Each $\sss_n$ stabilizes an arc $\ell^n$ in
$T$. Notice that if $\sss_n$ stabilizes a non-trivial arc in $T$
then the image $\sss_n'$ of $\sss_n$ in $G$ also stabilizes this
arc. Hence any finitely generated subgroup of $[\sss,\sss]$ is
inside a stabilizer of a nontrivial arc in $T$. This proves (1).

Statement (2) is proved similar to (1).
\endproof

\section{Groups acting on tree-graded spaces}
\label{gatg}

Let $G$ be a finitely generated group acting by isometries on a
 tree-graded space $(\free , \pp)$. Let $S=S\iv$ be a finite generating set of $G$.
 Given a piece $B\in \pp$ and a
 point $p\in B$ we denote by $\Stab(B)$ (resp. $\Stab(p)$ and $\Stab(B,p)$) the
 stabilizer of $B$ (resp. stabilizer of $p$ and the intersection
 $\Stab(B)\cap \Stab(p)$).

\medskip

Let us fix the following notation:
\begin{itemize}

\item $\calc_1(G)$ is the set of stabilizers of subsets of $\free$
all of whose finitely generated subgroups stabilize pairs of
distinct pieces in $\pp$.

\item $\calc_2(G)$ is the set of stabilizers of pairs of
points of $\free$ not from the same piece.

\item $\calc_3(G)$ is the set of stabilizers of triples of
points of $\free$ neither from the same piece nor from the same
transversal geodesic.
\end{itemize}

\medskip

By Lemma \ref{inv}, the action of $G$ on $\free $ induces an
action by isometries of $G$ on the $\R$--tree $T=\free/{\approx}$.
The main result of this section is the following theorem.

\begin{theorem}\label{1} Let $G$ be a finitely generated group
acting on a tree-graded space $(\free,\pp)$.
 Suppose that the following hold:
 \begin{itemize}
 \item[\textbf{(i)}] every isometry $g\in G$ permutes the pieces;
 \item[\textbf{(ii)}] no piece in $\pp$ is stabilized by the whole group $G$;
 likewise no point in $\free$ is fixed by the whole group $G$.
 \end{itemize}
 Then one of the following four situations occurs:
\begin{itemize}
\item[(\textbf{I})] the group $G$ acts by isometries on the real
tree $T=\free/{\approx}$
 non-trivially, with stabilizers of non-trivial
 arcs in $\calc_2(G)$, and with stabilizers of non-trivial tripods in $\calc_3
 (G)$;
 \item[(\textbf{II})] there exists a point $x\in \free$ such that for
 any $g\in G$ any
 geodesic $[x,g\cdot x]$ is covered by finitely many pieces: in this case
 the group $G$ acts non-trivially on a simplicial tree with
 stabilizers of vertices of the form $\Stab(B)$, $B\in \pp$, or of
 the form $\Stab(p)$, $p\in \free$, and  stabilizers of edges
 of the form $\Stab(B,p)$;
\item[(\textbf{III})] the group $G$ acts non-trivially on a
simplicial tree with edge stabilizers from $\calc_1(G)$;
 \item[(\textbf{IV})] the group $G$ acts on a complete
{$\R$--}tree by isometries, non-trivially, stabilizers of
non-trivial arcs are locally inside $\calc_1(G)$-by-Abelian
subgroups, and stabilizers of tripods are locally inside subgroups
in $\calc_1(G)$; moreover if $G$ is finitely presented then the
stabilizers of non-trivial arcs are in $\calc_1(G)$.
\end{itemize}
\end{theorem}

{\bf Case A}. Suppose that the action of  $G$ on $T=\free/{\approx}$
does not have a global fixed point. Then the action of $G$ on $T$
has all the other required properties from (I). Indeed, Lemma
\ref{inv}, (2), and Remark \ref{211} imply that any stabilizer in
$G$ of a non-trivial arc of $T$ coincides with the stabilizer of two
distinct points in $\free$, not contained in the same piece,
therefore it is an element of $\calc_2(G)$. The same results imply
that the stabilizer of a tripod in $T$ coincides with the stabilizer
of a triple of points in $\free$ projecting onto the vertices of the
tripod in $T$. These three points  are not in the same piece nor on
the same transversal geodesic (otherwise their images in $T$ would
coincide or they would be on a geodesic).

{\bf Case B.} Suppose that $G$ fixes a point in $T$. Let $t\in T$
be this point. Let $R$ be the $\approx$--equivalence class
projecting onto $t$. By Lemma \ref{equiv1}, $R$ is tree-graded
with respect to the pieces from $\pp$ contained in $R$. Let $\cal
R$ be this set of pieces.

Lemmas \ref{ppp} and \ref{up}  allow us to define a transfinite
sequence of tree-graded structures on $R$. Set $\pp_0=\cal R$, for
every non-limit cardinal $\alpha+1$, we define $\pp_{\alpha+1}$ as
$\pp_\alpha'$, and for every limit cardinal $\alpha$ we set
$\pp_\alpha=\bigcup_{\beta<\alpha} \pp_\beta$.

\medskip

\begin{definition}\label{int}
Let $B$ be an arbitrary piece in $\pp_\alpha$, $\alpha \geq 1$.
Assume that $\alpha$ is not a limit cardinal, and let
$\sim_{\alpha-1}$ be the $\sim$--equivalence relation defined in
Section \ref{nested} corresponding to $(\free ,\pp_{\alpha-1})$.
Then $B$ is the closure of a $\sim_{\alpha-1}$--equivalence class.
We denote this equivalence class by $\Int (B)$ and we call it the
\textit{interior} of $B$.

Assume that $\alpha$ is a limit cardinal. Then $B$ is, by Lemma
\ref{up}, the closure of an increasing union of pieces $B_\beta$
from $\pp_\beta$, $\beta<\alpha$. In this case we denote
$\bigcup_{\beta<\alpha} B_\beta$ by $\Int (B)$ and also call it the
\textit{interior} of $B$.
\end{definition}

\medskip

Note that Lemma \ref{geod} in the first case and Corollary
\ref{convsat} in the second case imply that the interior of a piece
in $\pp_\alpha$, with $\alpha \geq 1$, is always convex.

\begin{lemma}\label{claim} If $B\ne B'\in \pp_\alpha$ then
$\Int(B)\cap\Int(B')=\emptyset$.
\end{lemma}

\proof This is obviously true if $\alpha$ is not a limit cardinal or
if either $B$ or $B'$ is a singleton. Suppose that $\alpha$ is a
limit cardinal, and that $B=\overline{\cup_{\beta<\alpha} B_\beta},
B'=\overline{\cup_{\beta<\alpha}B_\beta'}$ are not singletons. If
there exists $p\in \Int(B)\cap\Int(B')$ then $p\in B_\beta\cap
B_\beta'$ for some $\beta<\alpha$. Therefore $B_\beta$ and
$B_\beta'$ are inside the same piece of $\pp_{\beta+1}$. Moreover
$B_\beta$ and $B_\beta'$ are not singletons. Property $(T_1)$ of
tree-graded spaces then implies that $B_\xi=B'_\xi$ for every
$\xi>\beta$. Hence $B=B'$, a contradiction. \endproof

 It is obvious that the
sequence $\pp_\alpha$ stabilizes (for example when the cardinality
of $\alpha$ is bigger than the cardinality of the set of collections
of subsets of $R$). It is also clear that if $\pp$ contains two
pieces that intersect then $\pp'\ne \pp$. Hence for some $\alpha$,
$\pp_\alpha$ consists of pairwise non-intersecting pieces. Then the
equivalence classes of the relation $\sim_\alpha$ corresponding to
$\pp_\alpha$ are just the pieces from $\pp_\alpha$. We shall assume
that $\alpha$ is minimal with this property. By Lemma \ref{perm}, we
have an induced action of $G$ on $\tilde R=R/\sim_\alpha$. Consider
two cases.

{\bf Case B.1.} Assume the group $G$ fixes a point  in $\tilde R$,
that is $G$ stabilizes a piece in $\pp_\alpha$. Let $\delta\leq
\alpha$ be the minimal cardinal such that $G$ fixes a piece $A$ in
$\pp_\delta$. By (ii),  $\delta\geq 1$.

Pick a point $x\in \Int A$. By Lemmas \ref{stabk} and \ref{stabu},
$G\cdot x \subset \Int A$.

Let us modify the set of generators $S$ as follows. We know that
there exists $s\in S$ such that $s\cdot x\ne x$. Let $s_1,...,s_k$
be all elements of $S$ such that $s_i\cdot x=x$. Then let us replace
each $s_i$ by $ss_i$, and $s_i\iv$ by $s_i\iv s\iv$. The set $S'$ of
generators thus obtained is closed under taking inverses, and no
element of $S'$ fixes $x$. Without loss of generality we can assume
that $S$ itself satisfies this property.

Since $\Int (A)$ is convex, geodesics connecting pairs of points
from $G\cdot x$ are in $\Int (A)$. Let us define an image of the
Cayley graph $\Cay(G,S)$ in $\Int (A)$. For every $s\in S$ choose a
geodesic $\cf_s$ connecting $x$ and $s\cdot x$ such that
$s\cf_{s^{-1}}=\cf_s$. Now for every $g\in G$ let the image of the
edge $(g, gs)$ be $g\cdot \cf_s$. These geodesics will be also
called {\em edges}. Note that by our assumption about $S$, none of
the edges can be of length 0. This will be needed in the proof of
Lemma \ref{notalimit} below. Thus $\Int (A)$ contains an image
$\mathfrak{C}$ of the Cayley graph $\Cay(G,S)$ such that all edges
of $\Cay(G,S)$ map to geodesic intervals in $\Int (A)$.

\begin{lemma}\label{notalimit} $\delta$ is not a limit cardinal.
\end{lemma}

\proof Suppose that $\delta$ is a limit cardinal. Then
$\pp_\delta=\cup_{\beta<\delta} \pp_\beta$. Therefore $\Int A$ is an
increasing union of pieces $A_\beta\in \pp_\beta$, $\beta<\delta$.
Then there exists $\beta<\delta$ such that $\bigcup_{s\in S}\cf_s
\subset A_\beta$. By Lemma \ref{perm}, for every $g\in G$ the union
$\bigcup_{s\in S} g\cdot \cf_s$ belongs to the piece $g\cdot
A_\beta\in \pp_\beta$. Using induction on the length $|g|_S$, we
prove that $g\cdot A_\beta=A_\beta$.

For $|g|_S=0$ this is trivial. Suppose that the statement is true
for $|g|=n$, consider an arbitrary element $|g|$ of length $n+1$.
Then $g=g_1s$, $|g_1|_S=n$. By induction the edge $g_1\cf_s$ is in
$A_\beta$. On the other hand, $g_1\cf_s$ is the same geodesic as
$g\cf_{s^{-1}}$ in $g\cdot A_\beta$. So the pieces $A_\beta$ and
$g\cdot A_\beta$ have a non-trivial arc in common, hence
$A_\beta=g\cdot A_\beta$ by property $(T_1)$ of tree-graded spaces.

This implies that $G\cdot x\subset A_\beta$, hence that $G\cdot
A_\beta =A_\beta$, a contradiction with the minimality of $\delta$.
\endproof

According to Lemma \ref{notalimit} and (ii), there exists
$\delta-1$. The collection $\pp_\delta$ is thus equal to
$\pp_{\delta-1}'$. We have that $\mathfrak{C}\subset \Int (A)$,
hence all edges $g\cdot \cf_s$ are covered by finitely many pieces
in $\pp_{\delta-1}$. On the other hand, no piece in $\pp_{\delta-1}$
contains the whole graph $\mathfrak{C}$.

Consider the following graph $\Gamma$. The vertices of $\Gamma$ are
of two types. Vertices of the first type are the pieces in
$\pp_{\delta-1}$ intersecting non-trivially the edges of
$\mathfrak{C}$. Vertices of the second type are the intersection
points of the pieces representing vertices of the first type. We
connect a vertex of the first type $B$ with a vertex of the second
type $p$ if and only if $p\in B$. It is easy to show using Corollary
\ref{convsat} and Lemma \ref{stronglyconvex} that $\Gamma$ is a
simplicial tree. The action of $G$ on $\free$ induces a simplicial
action of $G$ on $\Gamma$ (by Lemma \ref{perm}). This action does
not fix a point. Indeed, by the minimality of $\delta$, no vertex of
the first type is fixed by the whole $G$. The fact that no vertex of
the second type is fixed follows from (ii). The group $G$ cannot fix
a midpoint of any edge either, because it would have to fix its
endpoints.

Suppose that $\delta-1=0$. Then $\pp_{\delta -1}=\mathcal{R}$. Let
$K$ be the stabilizer of an edge $(B,p)$. Then, by definition, $K$
is of the form $\Stab(B,p)$, $B\in \pp, p\in \free$. Note that in
this case and by Lemma \ref{stronglyconvex}, any geodesic $[x,g\cdot
x],\, g\in G,$ is covered by a finite number of pieces from
$\mathcal{R}$. Thus case (II) of the theorem occurs.

Therefore we can assume that $\delta>1$. Consider an edge $(B,p)$ of
$\Gamma$. Let $K\le G$ be the stabilizer of this edge. We have that
$p\in B$, $B=\overline{\Int(B)}$.

\begin{lemma}\label{claim1} Suppose that $p$ belongs to $\Int(B)$. Then the
stabilizer of the edge $(B,p)$ coincides with the stabilizer of the
vertex $p$.
\end{lemma}

\proof We can assume that $B$ is not a singleton. Let $\xi$ be the
smallest cardinal such that a $\pp_\xi$--piece containing $p$ is
not a singleton. Then $\xi<\delta-1$ since $\Int(B)$ is not a
singleton either. Now suppose that $g\in G$ fixes $p\in \Int(B)$.
Then $g$ permutes the pieces of $\pp_\xi$ containing $p$. The
union of these pieces is inside a piece from
$\pp_{\xi+1}=\pp_\xi'$. Note that $\xi+1\le \delta-1$. The element
$g$ stabilizes this union. Since different pieces of
$\pp_{\delta-1}$ intersect by at most a point, $g$ must stabilize
$B$. Thus we proved that $\Stab(B,p)=\Stab(p)$.
\endproof

Edges $(p,B)$ of $\Gamma$ satisfying the condition $p\in \Int (B)$
will be called {\em redundant}. Lemma \ref{claim} shows that for
every vertex $p$ of $\Gamma$ there can be at most one redundant edge
of the form $(B,p)$.

Let $\Gamma'$ be the simplicial tree obtained by collapsing the
redundant edges of $\Gamma$ (i.e. removing the interior of each
redundant edge and identifying its ends). One can describe $\Gamma'$
directly as follows. The vertices of $\Gamma'$ are of two types.
Vertices of the first type are pieces from $\pp_{\delta-1}$
intersecting non-trivially the edges of $\mathfrak{C}$. Vertices of
the second type are points $p=B_1\cap B_2$ where $B_1, B_2$ are
vertices of the first type, such that $p\not\in \Int(B)$ for any
vertex of the first type {$B$}. Vertices $B$ and $B'$ of the first
type are connected if and only if $\Int(B)\cap B'$ or $B\cap
\Int(B')$ is not empty. A vertex $B$ of the first type is connected
to the vertex $p$ of the second type if $p\in B$.

According to Lemmas \ref{stabk} and \ref{stabu} the group $G$
permutes redundant edges of $\Gamma$, hence it acts on $\Gamma'$ by
simplicial automorphisms.

Let us prove that the stabilizer $K$ of any edge in $\Gamma'$ is in
$\calc_1(G)$. If $K$ stabilizes an edge of type $(B,B')$, i.e. it
stabilizes two pieces $B,B'\in \pp_{\delta-1}$ with the intersection
point $p$ in $\Int (B)\sqcup \Int (B')$, then $K$ stabilizes $p$.
Suppose that $p\in \Int(B')$. Then $K\subset \mathrm{Stab}(B,p)$,
moreover by Lemma \ref{claim1}, $\mathrm{Stab}(B,p)\subset
\mathrm{Stab} (B,B')=K$. Thus $K=\mathrm{Stab}(B,p)$.

Similarly, if $(B,p)$ is an edge of $\Gamma'$ where $B\in
\pp_{\delta-1}$, $p\in\free$, then $p\not\in\Int(B)$, and
$K=\mathrm{Stab}(B,p)$. Hence it is enough to prove the following
statement.

\begin{lemma} \label{claim2} For every $B\in \pp_{\delta-1}$ and
$p\in B\setminus\Int(B)$ the stabilizer $K$ of the pair $(B,p)$ is
in $\calc_1(G)$.
\end{lemma}

\proof Suppose first that $\delta-1$ is a limit cardinal. Then
$B=\overline{\bigcup_{\beta <\delta-1}B_\beta }$.

According to Lemma \ref{stabu}, for every $g\in K$ there exists
$\beta(g)$ such that $gB_\beta =B_\beta$ for every $\beta
>\beta(g) $.

Let $K_1$ be a finitely generated subgroup of $K$, and let
$k_1,...,k_m$ be a set of generators of $K_1$. There exists
$\beta_0$ such that for every $\beta>\beta_0$, $k_i B_\beta
=B_\beta$ for all $i$. Therefore $K_1$ also stabilizes the piece
$B_\beta$, and it fixes the point $p$, hence it fixes the projection
$y$ of $p$ onto $B_\beta$. Since $p,y$ are in $R$, $p\approx y$. On
the other hand $p \nsim_\beta y $ (otherwise $p\in B_{\beta+1}$),
hence $p \nsim_0 y $ and $\Sato \{p,y\}$ contains at least two
pieces. Then $K_1$ must also fix these pieces. Hence $K_1$
stabilizes a pair of distinct pieces.

\medskip

Now suppose that $\xi=\delta-2$ exists. Let again $K_1$ be a
finitely generated subgroup of $K$ and let $S$ be a finite set
generating $K_1$. Let $y$ be a point in $\Int (B)$ and let
$y_0,y_1,..., y_n,...$ be a sequence of points converging to $p$
on the geodesic $[y,p]$ as in Lemma \ref{seqcut}. For every $s\in
S$, according to Lemma \ref{stabkp}, there exists an $n_s$ such
that $s$ fixes all $y_n$ with $n\geq n_s$. It follows that the
whole group $K_1$ fixes all $y_n$'s for $n\ge m=\max_{s\in S}
n_s$. Let $z=y_m$ and $z'=y_{m+2}$. They are by construction in
the same $\approx$--class $R$, and $\Sato \{ z,z'\}$ must contain
at least two distinct pieces in $\mathcal{R}$. Since $K_1$
stabilizes $\{ z,z'\}$, it must stabilize these two pieces, by
Remark \ref{211}. We conclude that $K\in \calc_1(G)$.
\endproof

Summarizing, we obtain the following proposition.

\begin{proposition} \label{ii1} In Case B.1,
either $G$ acts non-trivially on a simplicial tree with edge
stabilizers in $\calc_1(G)$, or property (II) of Theorem \ref{1}
holds.
\end{proposition}
\medskip

{\bf Case B.2.} Suppose $G$ does not stabilize a point in $\tilde
R=R/\sim_\alpha$.

Take a point $p$ in $R$, and consider the corresponding image
${\mathfrak C}$ of the Cayley graph $\Cay(G,S)$  as in Case B.1.
Recall that it consists of vertices $g\cdot p,\, g\in G$, connected
by edges $g\cdot \cf_s , s\in S,$, where $\cf_s$ is a geodesic
joining $p$ and $sp$, with the assumption that
$s\cf_{s^{-1}}=\cf_s$.

Let us define a ternary {\em betweenness} relation on $\tilde R$ as
follows. For every three points $x,y,z$ in $\tilde R$ we set $xyz$
if there exists a geodesic connecting a point from $x$ with a point
from $z$ and containing a point from $y$. Recall that pieces from
$\pp_\alpha$ are pairwise disjoint, thus a point from $y$ is neither
in $x$ nor in $z$.

\begin{lemma}\label{intervals} For every $x,y,z\in \tilde R$, $xyz$ if
 and only if any geodesic connecting a point in $x$ with a point in $z$
intersects $y$.
\end{lemma}

\proof It immediately follows from property $(T_2')$.
\endproof

\begin{lemma}\label{bow} The set $\tilde R$ with the betweenness
  relation is a pretree in the sense of Definition \ref{pretreed}.
\end{lemma}

\proof (PT0) $(\forall x,y) (\neg xyx)$. Indeed, suppose that
$a,b,c$ are points in $x,y,x$ respectively and a geodesic $[a,c]$
contains $b$. Since $x$ is convex (by Lemma \ref{union}), $[a,c]$ is
inside $x$. Therefore $x$ and $y$ must intersect, a contradiction.

(PT1) $xzy \Leftrightarrow yzx$. That is obvious.

(PT2) $(\forall x,y,z) (\neg (xyz \wedge xzy))$. Indeed, suppose
that $a,a'\in x, b,b'\in y, c,c'\in z$ and $[a,c]$ contains $b$
while $[a',b']$ contains $c'$. Let $b''$ be the entry point of
$[a',b']$ in $y$, and let $b'''$ be the exit point of $[a,c]$ from
$y$. Then the union $[c', b'']\sqcup [b'',b''']\sqcup [b''',c]$ is
an arc by \cite[Lemma 2.28]{DS}. By Property $(T_2)$ we have that
any geodesic connecting $c'$ and  $c$ must pass through $b''$.
Therefore we have that $z$ is between $y$ and $y$ which contradicts
(PT0).

(PT3) $xzy$ and $z\neq w$ then $(xzw \vee yzw)$. Indeed, suppose
that $\neg xzw$. Since $xzy$, we can find  $a\in x, b\in y, c\in z$
such that $[a,b]$ contains $c$. Since $\neg xzw$, any geodesic
connecting $a$ and a point $e$ in $w$ does not intersect $z$. Let
$a'$ be the farthest from $a$ intersection point of $[a,b]$ and
$[a,e]$. Then the union $[e,a']\sqcup [a',b]$ is an arc. Since it
passes through $z$, every geodesic connecting $b$ and $e$ passes
through $z$, and we have $yzw$. \endproof

\medskip

For every two points $x,y$ in the pretree $\tilde R$ the set
$\{z\mid xzy\}$, denoted by $\lfloor x,y\rfloor$, is called the
{\textit{interval}} between $x$ and $y$.

By $\rfloor x,y\lfloor$ we denote $\lfloor x,y\rfloor$ without
$x,y$, and we call it
 {\textit{the open interval}} between $x$ and $y$.

Every interval $\lfloor x,y\rfloor$ has a natural order: $a<b$ if
$a\in \lf x,b\rf$.

\begin{lemma} \label{median} {$\tilde R$ is a median pretree in the terminology of \cite{Bow}, i.e.
}for every three points $x,y,z\in \tilde
  R$ the intersection of intervals $\lf x,y\rf$, $\lf x,z\rf$ and $\lf
  y,z\rf$ is a singleton in $\tilde R$.
\end{lemma}

\proof Consider a geodesic triangle $abc$ in $R$ where $a\in x, b\in
y, c\in z$. By Lemma \ref{intervals}, $\lf x,y\rf$ (resp. $\lf
y,z\rf$, $\lf x,z\rf$) consist of all pieces intersecting the side
$[a,b]$ (resp. $[b,c]$, $[a,c]$). Let $a'$ be the farthest from $a$
common point of $[a,b]$, $[a,c]$, and $b'$, $c'$ defined likewise.
Then $a'b'c'$ is a simple geodesic triangle in $R$. By Property
$(T_2)$, $a'b'c'$ is contained in a piece. Since pieces of
$\pp_\alpha$ do not intersect, this piece is unique. It is the
intersection of the intervals $\lf x,y\rf$, $\lf x,z\rf$ and $\lf
  y,z\rf$.
\endproof

\medskip

\noindent {\textit{ Notation:}} For every three points $x,y,z\in
\tilde
  R$ we denote the unique common point of $\lf x,y\rf$, $\lf x,z\rf$ and $\lf
  y,z\rf$  by $m(xyz)$.

\medskip

\begin{definition}\label{ptdense}
We call a subset $U$ of a set $V$ {\em dense} in $V$  if between any
two points $x,y$ of $V$, there is a point from $U$ distinct from
$x,y$.
\end{definition}

The following lemma is well known (see, for example, \cite{Haus}).

\begin{lemma}\label{order}
A countable ordered set $L$ is order isomorphic to the set of
rational numbers from the open unit interval of the real line if and
only if $L$ is dense in itself and every element of $L$ is between
two other elements of $L$ (i.e. there are no terminal elements in
$L$).
\end{lemma}

\begin{lemma}\label{separated} For every two points $x,y$ in $\tilde
  R$ the interval $\lfloor x,y\rfloor$ in the  pretree $\tilde R$ contains a dense subset
that is order isomorphic to the set of rational numbers in the unit
interval.
\end{lemma}

\proof Indeed, let $x,y\in \tilde R$. Let $\pgot$ be a geodesic
connecting a point $a\in x$ to a point $b\in y$. We assume that
$a,b$ are the exit and entry point of $\pgot$ from $x$ and to $y$
respectively. Since different $\pp_\alpha$--pieces do not
intersect, and the transversal trees of $R$ are trivial by Lemma
\ref{transv}, there are countably many $\pp_\alpha$--pieces that
intersect $\pgot$ nontrivially, and between any two intersections
of pieces with $\pgot$, there is a non-trivial intersection of
$\pgot$ with a piece.  By Lemma \ref{order}, the set of pieces
that intersect $\pgot$ non-trivially form an ordered set that is
order isomorphic to ${\mathbb Q}\cap (0,1)$.
\endproof

 Note that every $\R$--tree can be considered as a dense pretree
with the natural betweenness relation: $xyz$ if and only if $y$
belongs to the geodesic $[x,z]$.

\begin{proposition}\label{Prtree}
The pretree $\tilde R$ is isomorphic to an $\R$--tree $\Theta$.
\end{proposition}

The proof of Proposition \ref{Prtree} is done in several steps.

\begin{definition}
A subset $S$ in a pretree $R$ is called \textit{connected} if for
every two points $s_1,s_2$ in $S$, $\lfloor s_1, s_2 \rfloor$ is
contained in $S$.
\end{definition}

\begin{lemma}\label{rtree} The pretree $\tilde R$ can be embedded into
  an $\R$--tree $\Theta$ such that the following properties hold:

\begin{enumerate}
\item [(R1)] $\tilde R$ is dense in $\Theta$ in the sense of
  Definition \ref{ptdense};

\item[(R2)] every element in $\Theta$ is between two elements in $\tilde
R$;

\item[(R3)] the topology induced by the metric topology of $\Theta$ on every
  interval $\lf x,y\rf$ of $\tilde R$ coincides with the subinterval topology on
  this interval.
\end{enumerate}
\end{lemma}

\proof Consider the set $ \qq$ of triples $(S, \phi,\Theta_S)$ where

\begin{enumerate}
\item[(S1)] $S$ is a connected sub-pretree of $\tilde R$, i.e. for every
  $s_1,s_2\in S$, $\lf s_1,s_2\rf\subseteq S$;

\item[(S2)] $\Theta_S$ is an $\R$--tree, and $\phi$ is an embedding of
  pretrees with dense image;

\item[(S3)] the topology induced by the metric on $\Theta_S$
 on any interval $\lf s_1,s_2\rf$ of $S$  coincides with the subinterval topology;

\item[(S4)] every element of $\Theta_S$ is between two elements of
  $\phi(S)$.
\end{enumerate}

Note that the set $\qq$ is not empty because it contains the triple
$(\{s\},\iota,\{s\})$ where $s\in \tilde R$ and $\iota$ is the
identity map.

Define a partial order $\prec$ on $\qq$ as follows. We say that
$(S,\phi,\Theta_S)$ is smaller than $(S',\phi',\Theta_{S'})$ if
$S\subseteq S'$, $\Theta_S\subseteq \Theta_{S'}$, the metric of
$\Theta_S$ is the restriction of the metric of $\Theta_{S'}$, and
$\phi$ is the restriction of $\phi'$ onto $S$.

Let $(S_1,\phi_1,\Theta_{S_1})\prec
(S_2,\phi_2,\Theta_{S_2})\prec...$ be an increasing sequence of
elements of $\qqq$. Let $S=\bigcup S_i, \Theta_S=\bigcup
\Theta_{S_i}$ and $\phi\colon S\to \Theta_S$ be defined in the
natural way. It is obvious that $(S,\phi,\Theta_S)\succ
(S_i,\phi_i,\Theta_{S_i})$ for every $i$ and that
$(S,\phi,\Theta_S)$ satisfies (S1)-(S4). Thus $(\qq,\prec)$
satisfies the condition of the Zorn lemma.

Let $(S,\phi,\Theta_S)$  be a maximal element of $\qqq$. Let us show
(by contradiction) that $S=\tilde R$. Suppose that $S\ne \tilde R$.

Consider the union $\cup S$ of $\pp_\alpha$--pieces contained in
$S$. Then $\cup S\ne R$. Note that $\cup S$ is connected. Moreover
for any two points $a,b$ in $\cup S$ the saturation $\Sat \{
a,b\}$ is in $\cup S$. Indeed, consider a geodesic $\pgot$ in $R$
connecting $a$ and $b$. Since $R$ is convex, $\pgot$ is covered by
its intersections with $\pp_\alpha$--pieces from $R$. By (S1), all
these pieces are in $S$. Therefore $\Sat \pgot\subseteq \cup S$.

\medskip

{\bf Case 1.} Suppose that $\cup S$ is not closed in $R$. Take $p$
on the boundary $\partial(\cup S)\setminus \cup S$, and the unique
$\pp_\alpha$--piece $P$ containing $p$. Let $S'=S\cup\{P\}$.

Consider a geodesic $\g$ connecting a point $a\in \cup S$ with
$p$. By Lemma \ref{union}, $\g \setminus \{p\}$ is contained in
$\cup S$. Let $a_n$, $n\ge 0$, be a sequence of elements of
$\g\setminus \{p\}$ converging to $p$ such that $a_0=a$ and
$a_{n}\in [a,a_{n+1}]$. Let $A_n$ be the $\pp_\alpha$--piece
containing $a_n$, $n=0,1,2,...$. Then $\lf A_0,A_n\rf\subset\lf
A_0, A_{n+1}\rf$  in the pretree $\tilde R$. Note that every piece
in $\lf A,P\lf$ is inside one of $\lf A_0, A_n\rf$. Therefore
every piece in $\tilde R$ between $P$ and $A$ is in $S'$. Hence
$S'$ is connected.

Consider the increasing union $U=\bigcup [\phi(A_0), \phi(A_n)]$
of geodesic intervals in the $\R$--tree $\Theta_S$. Then $U$ is
either a geodesic ray or a half-open geodesic interval in
$\Theta_S$.

Suppose that $U$ is a half-open geodesic interval in $\Theta_S$. We
define $\Theta_{S'}$ as the tree obtained from $\Theta_S$ by adding,
if necessary, an endpoint $t$ to the half-open geodesic interval and
extending the metric in the natural way. Then consider $\phi_{S'}$
as the extension of $\phi$ to $S'$ defined by $\phi_{S'} (P)=t$.

We show that $(S',\phi_{S'},\Theta_{S'} )$ satisfies the
conditions (S1)-(S4). That would contradict the maximality of
$(S,\phi_S, \Theta_S)$. Property (S1) is already proved, (S3) and
(S4) are obvious.

In order to show (S2) it suffices to prove that if $B,B'\in S$ are
such that $B'$ is between $B$ and $P$, then $\phi_{S'}(B')$ is
between $\phi_{S'}(B)$ and $\phi_{S'}(P)$.  Let $b\in B$. Lemma
\ref{t2} implies that for $n$ large enough the piece $B'$
intersects the geodesic $[b,a_n]$. Consequently $\phi_{S'}(B')$ is
in the geodesic $[\phi_{S'}(B),\phi_{S'}(A_n)]$ for $n\geq n_0$.
Thus, it is in $[\phi_{S'}(B),\phi_{S'}(P)]$.

Suppose that $U$ is a geodesic ray in $\Theta_S$. Then we redefine
the metric on $\Theta_S$ so that the topology stays the same but
the ray $U$ becomes isometric to $[0,1)$, for instance by taking
the stereographic projection of $\mathbb R$ onto the unit circle
and restricting it to $[0,\infty)$. Note that after this
modification of the metric on the tree $\Theta_S$ properties
(S1)-(S4) still hold. This and the previous case will end the
argument.

\medskip

{\bf Case 2.} Now suppose that $\cup S$ is closed. Let $p\in
R\setminus \cup S$.

For every point $a$ in $\cup S$ a geodesic $\g_a$ joining $a$ with
$p$ contains a point $p_a$ which is the nearest to $p$ point in
$\cup S$. If two geodesics $\g_a , \g_b$ are such that $p_a\neq
p_b$ then it follows that $p_a, p_b$ and a point $c$ in $\g_a \cap
\g_b$ are the vertices of a simple geodesic triangle. Then
$p_a,p_b$ and $c$ are in the same $\pp_\alpha$--piece, hence this
piece is also in $\cup S$, and so is $c$. This contradicts the
choice of $p_a$ and $p_b$. Thus $p_a$ is a point $p'$ that does
not depend on the point $a$ or on the geodesic $\g_a$, $p'$ is in
fact the nearest point projection of $p$ onto $\cup S$.

Let $P$ and $P'$ be the $\pp_\alpha$--pieces containing $p$ and
$p'$, respectively. Let $\g$ be a geodesic connecting $p'$ and
$p$, and let $\Sat\g$ be its saturation. Let $S'$ be the union of
$S$ and all pieces in $\Sat\g$. Clearly $S'$ is a connected
sub-pretree of $\tilde R$. By Lemma \ref{separated}, there exists
a countable dense subset $Y$ of the interval $\lf P', P\rf$,
containing $P$ and $P'$, and order isomorphic to ${\mathbb Q}\cap
[0,1]$. Let $\phi'$ be an order isomorphism from $Y$ onto
${\mathbb Q}\cap [0,1]$, such that $\phi' (P')=0$ and $\phi'
(P)=1$. The map $\phi'$ can be uniquely extended to an
order-preserving map from $\lf P', P\rf$ to a dense subset of the
interval $[0,1]$ of the real line. We denote that extension by
$\phi'$ also. The union $\Theta_{S'}$ of $\Theta_S$ and $[0,1]$,
with $0$ identified to the point $\phi (P')$, is an $\R$--tree
with the natural metric extending the metric on $\Theta_S$ and the
interval metric on $[0,1]$. The map $\phi_{S'}$ is defined as
$\phi$ on $S$ and as $\phi'$ on $\lf P', P\rf$. It is easy to see
that $\phi_{S'}$ is an embedding of the pretree $S'$ into
$\Theta_{S'}$ and properties (S3) and (S4) are satisfied.
Properties (S1) and (S2) have been established before. This
contradicts the maximality of the triple $(S,\phi,\Theta_S)$.
\endproof

\begin{remark}
After completing the proof of Lemma \ref{rtree}, we discovered a
paper by Bowditch and Crisp \cite{BC} where for an arbitrary
pretree $\tilde P$ satisfying some mild conditions an embedding
into an $\R$--tree $\Phi$ is constructed, so that every
automorphism of $\tilde P$ extends to an automorphism of the
$\R$--tree $\Phi$ (considered as a pretree), and non-nesting
automorphisms are extended to non-nesting automorphisms (for the
definition of ``non-nesting" see Definition \ref{nn})

Thus, it would suffice to prove that our pretree $\tilde R$
satisfies their conditions and then apply \cite{BC}. Still we
decided to include our original proof. The reason is that the
embedding in \cite{BC} is based on a relatively complicated
construction from Bowditch \cite{Bow}, and our embedding is
straightforward and the proof is much shorter. Note that our
construction in Lemma \ref{rtree} can shorten the arguments in
\cite{Bow} too. Note also that other embeddings of pretrees into
$\R$--trees have been considered in Chiswell \cite{C}.
\end{remark}

\medskip

\noindent {\textit{Conventions:}} In what follows, we identify
$\tilde R$ with its image in $\Theta$. Eventually we shall show,
proving Proposition \ref{Prtree}, that $\tilde R=\Theta$.

For two points $x,y\in \Theta$ we denote the geodesic joining them
by $[x,y]$.

\medskip

\begin{lemma}\label{inested}
Every point $p$ in $\Theta\setminus \tilde R$ is the unique
intersection of a nested sequence of intervals $[x_n,y_n]$ with
$x_n, y_n\in \tilde R$.
\end{lemma}

\proof Property (R2) of Lemma \ref{rtree} implies that $p\in
[x_0,y_0]$ with $x_0,y_0\in\tilde R$. By Lemma \ref{separated}, the
interval $\lf x_0,y_0\rf$ in $\tilde R$ contains a countable dense
subset $U$. Let $U\cap [x_0,p]=\{u_0,u_1,..... \}$ and $U\cap [p,
  y_0]=\{v_0,v_1,..... \}$. Then $\{ p\} =\bigcap [u_n,v_n]$. Indeed,
suppose that $\bigcap [u_n,v_n]$ is an interval $[a,b]$ containing
$p$. Property (R1) applied twice implies that there exist $u,v\in
]a,p[\cap \tilde R$. Density of $U$ implies that there exists
$u_i\in ]u,v[\subset ]a,p[$. This yields a contradiction.

Now for every $n\in {\mathbb N}$ take $x_n$ to be the first $u_k$
at distance less than $1/n$ from $p$. Likewise we define $y_n$
using $v_k$. Obviously the sequence of intervals $[x_n,y_n]$ is
nested and its intersection is~$\{ p\}$.\endproof

\noindent \textit{Proof of Proposition \ref{Prtree}.}\quad  It
suffices to prove that ${\tilde R} = \Theta$, that is the
embedding in Lemma \ref{rtree} is surjective. Suppose that there
exists a point $p$ in $\Theta \setminus \tilde R$. By Lemma
\ref{inested}, $p$ is the unique intersection of a nested sequence
of intervals $[x_n,y_n]$ with $x_n, y_n\in \tilde R$. For every
$n\in \N$ let $a_n\in x_n$ and $b_n\in y_n$ be the pair of points
minimizing the distance in $R$. Since $x_n,y_n \in [x_{n-1},
y_{n-1}]$ it follows that $x_n,y_n$ intersect any geodesic
$[a_{n-1}, b_{n-1}]$, therefore by Lemma \ref{ab} the points
$a_n,b_n$ are contained in $[a_{n-1}, b_{n-1}]$. Thus one can
produce a nested sequence of geodesics $[a_n,b_n]$ in $R$ whose
intersection is a geodesic $[a,b]$. Let $x$ be the unique
$\pp_\alpha$--piece in $\tilde R$ containing $a$ and let $y$ be
the unique piece containing $b$. The interval $\lfloor x, y
\rfloor$ is contained in any $\lfloor x_n , y_n \rfloor$, hence in
$\Theta$ the arc $[x,y]$ is contained in $\bigcap [x_n,y_n] \{p\}$.
It follows that $x=y=p\in \tilde R$, a
contradiction.\hspace*{\fill}$\Box$

\begin{proposition}\label{non-nesting} The group $G$ acts on the
  $\R$--tree $\Theta$ by non-nesting pretree automorphisms,
  non-trivially, and with stabilizers of
  non-trivial arcs in $\calc_1(G)$.
\end{proposition}

\proof It suffices to prove that the action of $G$ on $\tilde R$
satisfies all the required properties.

{\bf Step 1.} The group $G$ acts on the pretree $\tilde R$ by
automorphisms. This follows from Lemma \ref{intervals} and the fact
that an isometry from $G$ takes geodesics in $R$ to geodesics in $R$
and permutes $\pp_\alpha$--pieces (by Lemma \ref{perm}).

\medskip

{\bf Step 2.} The action of $G$ on $\tilde R$ is non-nesting.
Indeed, if $I=\lfloor a,b \rfloor$ is an interval of $\tilde R$,
then by Lemma \ref{intervals}, $\lfloor a,b \rfloor$ consists of
all $\pp_\alpha$--pieces intersecting a certain geodesic $\g$ in
$R$. We can assume that $\g$ is the shortest geodesic joining the
pieces $a$ and $b$. Recall that by Lemma \ref{ab} the endpoints of
$\g$ are uniquely defined, and that any geodesic joining $a$ and
$b$ must contain them.

Suppose that $h\cdot I\subsetneq I$, $h\in G$. Note that $h\cdot I$
consist of all pieces intersecting $h\cdot \g$, and that $h\cdot \g$
is the shortest geodesic joining $h\cdot a$ and $h\cdot b$. Since
$h\cdot I\subset I$, the pieces $h\cdot a$ and $h\cdot b$ must
intersect $\g$ therefore the endpoints $h\cdot \g_-$ and $h\cdot
\g_+$ must be contained in $\g$. On the other hand, $h\cdot I\neq I$
implies that either $h\cdot \g_-\ne \g_-$ or $h\cdot \g_+\ne \g_+$.
In particular $\dist\left( h\cdot \g_-\, ,\, h\cdot \g_+ \right)$ is
smaller than the length of $\g$, a contradiction.

\medskip

{\bf Step 3.} The arc stabilizers of the action of $G$ on $\tilde R$
are in $\calc_1(G)$.

Indeed, let $[a,b]$ be a non-trivial interval in $\tilde R$ and
let $K$ be its stabilizer.

If $g\in G$ fixes $a$ and $b$ then in $R$ it fixes the two points
$x\in a$ and $y\in b$ such that $\dist(x,y)=\dist(a,b)$. According
to the proof of Lemma \ref{separated}, the strict saturation of
$\{ x,y\}$ in $(R, \pp_\alpha )$ contains countably many disjoint
pieces. Consequently the same is true for the strict saturation
$\Sato\{x,y\}$ in $(R, \calr).$ Thus every element $g$ fixing
$x,y$ must stabilize countably many distinct pieces in $\calr$.

%%%%%%%%%%%%%%end%%%%%%%%%%%%%%%%%%%%%%%%%

We conclude that $K$ is the stabilizer of the set $x\cup y$ in
$\free$, and that it stabilizes countably many distinct pieces in
$\calr$, hence $K$ is in $\calc_1(G)$.
\endproof

Theorem \ref{lev} in case when $G$ is finitely presented and
Theorem \ref{lev1} in the general case, and Proposition
\ref{non-nesting} imply the following statement.

\begin{proposition} \label{B2} In Case B.2, the group $G$ acts on an $\R$--tree
non-trivially by isometries so that
\begin{itemize}
\item[(a)] the stabilizers of non-trivial arcs are locally inside
$\calc_1(G)$-by-Abelian subgroups of $G$, and stabilizers of
tripods are locally inside subgroups in $\calc_1(G)$;

\item[(b)] if $G$ is finitely presented then the stabilizers of
non-trivial arcs are in $\calc_1(G)$.
\end{itemize}
\end{proposition}

This ends the proof of Theorem \ref{1}.

\begin{proposition}\label{prop1}
Let $G$ be a finitely  presented group satisfying the conditions
of Theorem \ref{1}. Suppose that $G$ stabilizes a (bi-infinite)
line $\calt$ in the $\R$--tree it acts upon according to the
theorem. Then $G$ stabilizes the strict saturation and the set of
cut-points of a bi-infinite geodesic $\g$ in $\free$ that is not
contained in a piece of $\pp$. In addition, any element from $G$
fixes the line $\calt$ pointwise if and only if it fixes the set
of cut-points from $\Cutp(\g)$ pointwise.
\end{proposition}

\proof
 We shall consider the different cases in the proof of Theorem \ref{1},
 which yield different real trees with a good action of $G$ on them.

\medskip

{\bf Case A}. Suppose that $\calt$ is a line in the $\R$--tree $T
=\free/{\approx}\, $. First we construct a geodesic $\mathfrak G$
whose projection onto $T$ is $\calt$.

Let $\calt_0$ be the maximal subinterval of $\calt$ that can be
obtained as the projection of a geodesic segment or ray from
$\free$. Note that we may suppose that $\calt_0$ is not a
singleton. Assume that $\calt_0\ne\calt$. Then by Lemma
\ref{Ttree} it is a geodesic segment or ray (i.e. $\calt_0$ cannot
be an open or semi-open interval), and we can find a geodesic
segment or ray $\g_0$ in $\free$ projecting onto $\calt_0$.  Let
$x$ be an endpoint of $\g_0$. The intersection of the
$\approx$--equivalence class containing $x$ with $\g_0$ is closed
and connected, by Lemma \ref{equiv}. It cannot contain the whole
$\g_0$, otherwise the projection of the latter would be a point.
Hence it is a geodesic segment $[x',x]$. By replacing if necessary
$x$ by $x'$, we may therefore assume that $\g_0$ contains only one
point in the $\approx$--equivalence class of $x$.

Let $\bar{x} \in \calt_0$ be the projection of $x$ in $T$. Then it
is an endpoint of $\calt_0$. Pick $\bar{y}\in \calt_0$, $\bar{y}
\neq \bar{x}$ and let $y\in \g_0$ be a representative of it. Let
$\bar{z}\in \calt \setminus \calt_0$ be such that $\bar{x} \in
(\bar{y},\bar{z})$ and let $z\in \free \setminus \g_0$ be a
representative of it. Then
$\widetilde\dist(y,z)=\widetilde\dist(y,x)+\widetilde\dist(x,z)$.

Consider two geodesics $[y,x]\subset \g_0$ and $[x,z]$ in $\free$.
Assume that they have a point $x'\neq x$ in common. Then
$$
\widetilde\dist(y,z)\leq
\widetilde\dist(y,x')+\widetilde\dist(x',z)\leq
\widetilde\dist(y,x)+\widetilde\dist(x,z)=\widetilde\dist(y,z)\, .
$$

It follows that $\widetilde\dist(y,x)=\widetilde\dist(y,x')$, hence
that $\widetilde\dist(x,x')=0$. This contradicts the choice of $x$
on $\g_0$.

Therefore, $[y,x]\cap [x,z]=\{x\}$, and $\cf=[y,x]\sqcup [x,z]$ is a
topological arc. By the choice of $x$, the maximal sub-arc in $\cf$
containing $x$ and appearing as intersection with a piece is either
$\{x\}$ or it has $x$ as an endpoint.

Let now $[y,z]$ be a geodesic in $\free$. By $(T_2')$ it has to
contain $x$. Therefore $\dist (y,z)=\dist (y,x) +\dist(x,z)$.

We have thus obtained that for all $y\in \g_0 \setminus \{x\}$,
 $\dist (y,z)=\dist (y,x) +\dist(x,z)$. One easily deduces from this that given
 a geodesic $[x,z]$,
  $\g_0 \sqcup [x,z]$ is a geodesic segment or ray. It projects onto a
  set strictly larger than
  $\calt_0$, contradicting its maximality.

We conclude that $\calt$ is the projection of a bi-infinite geodesic
$\g$. We now show that $G$
 stabilizes the strict saturation and the set of cut-points of this geodesic.

Let $x,y\in \g$ be two points with $\widetilde\dist(x,y){>}0$
  and such that $[x,y]\subset \g$ intersects the
  $\approx$--equivalence classes of $x$  and $y$
only in $x$ and $y$ respectively. In particular $x$ and $y$ belong
to $\Cutp\g$.

For every $g\in G$, the same can be said about $g\cdot{[x,y]}$ and
its endpoints. On the other hand, $g\cdot [x,y]$
  projects into $\calt$. Thus there
  exist $x'$ and $y'$ in $\g$ with $\widetilde\dist(x',g\cdot x)=0$ and
  $\widetilde\dist(y',g\cdot y)=0$.
  Since $\widetilde\dist(x',g\cdot x)=0$,
  any geodesic $[x', g\cdot x]$ intersects $g\cdot [x,y]$ only in
  $g\cdot x$, likewise for a geodesic $[y',g\cdot y]$. Then
  $[x', g\cdot x]\sqcup g\cdot [x,y]\sqcup [g\cdot y, y']$ is
  a topological arc, and $g\cdot x , g\cdot
  y$ are endpoints of intersections of that arc with pieces. It follows that
  the geodesic $[x',y']\subset \g$
   must contain $g\cdot x , g\cdot y$. This and Corollary
   \ref{convsat} imply that $\Sato g\cdot [x,y]= g\cdot \Sato [x,y]\subset \Sato\g$
   and that $\Cutp g\cdot [x,y] = g\cdot \Cutp[x,y]\subset \Cutp\g$.

Since there exists an
  increasing sequence of segments $[x_n,y_n]$ as above such that
 $\bigcup [x_n,y_n] = \g$, the equalities $G\cdot \Sato\g = \Sato\g$
and $G\cdot \Cutp\g=\Cutp\g$ follow immediately.

 Obviously if $g$ fixes $\Cutp\g$ pointwise then it fixes $\calt$ pointwise, because the
 projection of $\Cutp\g$ on $T$ is $\calt$. Let now
 $g$ in $G$ be such that $g\cdot \bar{t}=\bar{t}$ for every $\bar{t}\in \calt$.

 Let again $x,y\in \g$ be two points with $\widetilde\dist(x,y)>0$
  and such that $[x,y]\subset \g$
  intersects the $\approx$--equivalence classes of $x$  and $y$
only in $x$ and respectively $y$. With the argument above it follows
that $[x,y]$ must contain
  $g\cdot x , g\cdot y$.

Since $\dist (x,y)=\dist (g\cdot x , g\cdot y)$ it follows that
$g\cdot x =x$ and $g\cdot y = y$. Then $g\cdot
\Cutp\{x,y\}=\Cutp\{x,y\}$. Moreover, since $g$ is an isometry, it
fixes $\Cutp   \{x,y\}$ pointwise.

 We now take an increasing sequence of segments $[x_n,y_n]$ as above such that
 $\bigcup [x_n,y_n] = \g$. Then $\Cutp\g= \bigcup \Cutp\{x_n,y_n\}$, hence $g$
  fixes every point in $\Cutp\g$.

\medskip

{\bf Case B.1}. In this case $\calt$ is a line in the simplicial
tree $\Gamma$. Let $(B_n)_{n\in \Z }$ be the set of vertices of the
first type in $\calt$ enumerated in the order in which they appear
in $\calt$. Then each $B_n$ is a piece in $\pp_{\delta -1}$
intersecting $\mathfrak{C}$, and $B_n\cap B_{n+1}$ is a point $p_n$
for all $n\in \Z $. Let $\cf_n\subset B_n$ be a geodesic joining
$p_{n-1}$ to $p_n$. By \cite[Lemma 2.28]{DS}, $\g=\bigcup_{n\in \Z }
\cf_n$ is a geodesic line. The group $G$ stabilizes the set of
points $\{ p_n \mid n\in \Z \}$ on $\g$ hence it stabilizes
$\Sato\g$ and $\Cutp\g$. An element $g\in G$ fixes $\calt$ pointwise
if and only if it fixes pointwise $\{ p_n \mid n\in \Z \}$. Since
$\Cutp\g = \bigcup_{n\in \N } \Cutp\{ p_{-n} , p_n\}$ it follows
that $g\in G$ fixes pointwise $\{ p_n \mid n\in \Z \}$ if and only
if it fixes pointwise $\Cutp\g$.

\medskip

{\bf Case B.2}. Suppose that the isometric action of $G$ on the
$\R$--tree $T'$ from Theorem \ref{lev} stabilizes a line $L$. Then
by Theorem \ref{lev}, the action of $G$ on $\Theta$ stabilizes a
line as well. Let $\calt $ be that line. It can also be seen that
an element of $G$ fixes $L$ pointwise if and only if it fixes
$\calt$ pointwise by Lemma \ref{rtree}. Let $\caltr$ be the line
in $\tilde R$ corresponding to the line $\calt$ in $\Theta$. Let
$\lf l_n,r_n\rf$ be an increasing sequence of intervals in $\tilde
R$ such that $\bigcup \lf l_n,r_n\rf =\caltr$. Let $x_n\in l_n$
and $y_n\in r_n$ be the unique pair of points minimizing the
distance between the pieces $l_n,r_n$. The inclusion $\lf
l_n,r_n\rf \subset \lf l_{n+1},r_{n+1}\rf$ and Lemma \ref{ab}
imply that $x_n,y_n$ are contained in any geodesic joining
$x_{n+1}$ and $y_{n+1}$. It follows that if $\g_n$ is an arbitrary
geodesic of endpoints $x_{n+1}, x_n$, $\g_n'$ is an arbitrary
geodesic of endpoints $y_n, y_{n+1}$ and $\g_0$ is a geodesic
joining $x_1,y_1$, then
$$\g=...\sqcup \g_n \sqcup ... \sqcup \g_1
\sqcup \g_0 \sqcup \g_1' \sqcup ...\sqcup \g_n' \sqcup ...$$ is a
geodesic line in $\free$. Since for every $g\in G$, $g\cdot \lf
l_n,r_n\rf \subset \caltr$, it follows that $g\cdot l_n$ and
$g\cdot r_n$ intersect $\g$, hence by Lemma \ref{ab},  $g\cdot
x_n$ and $g\cdot y_n$ are contained in $\g$. Thus $G\cdot \{ x_n ,
y_n \}\subset \g$ for any $n\in \N$. An argument as in Case A
allows to deduce from this that $G$ stabilizes $\Sato \g$ and
$\Cutp\g$, and that $g\in G$ fixes $\caltr$ if and only if it
fixes $\Cutp\g$ pointwise.
\endproof

\begin{remark} \label{rem89}
In the proof of Proposition \ref{prop1}, one can replace
``finitely presented'' by ``finitely generated'' in Cases A and
B.1, since finite presentability is not used there.

In Case B.2, the finite presentability is used because we need to
apply Theorem \ref{lev}. Assume that a non-finitely presented
group $G$ stabilizes a bi-infinite line $L'$ in the $\R$--tree
$T'$. Since the group of isometries of the line is (torsion-free
Abelian)-by-$\mathbb{Z}/2\mathbb{Z}$, $G$ has an index at most 2
subgroup $G_1$ that is an extension of a subgroup $K$ fixing $L'$
pointwise by a finitely generated free Abelian group. By Theorem
\ref{lev1}, the derived subgroup $[K,K]$ is locally inside the
pointwise stabilizer of an arc in the pretree $\tilde R$ from the
Case B.2 of the proof of Theorem \ref{1}. Then as in Case B.2 of
the proof of Proposition \ref{prop1}, we can deduce that $[K,K]$
is locally inside a subgroup in $\calc_1(G)$. Thus $G_1$ is inside
a ($\calc_1(G)$-by-Abelian)-by-(finitely generated free Abelian
group).
\end{remark}

Combining Theorems \ref{1}, \ref{bf} and Proposition \ref{prop1},
 we obtain the following result.

\begin{theorem}[see Theorem \ref{split3}]\label{split1}
Let $G$ be a finitely presented group acting by isometries on a
tree-graded space $(\free ,\pp)$ such that:
\begin{itemize}
 \item[\textbf{(i)}] every isometry from $G$ permutes the pieces;
 \item[\textbf{(ii)}] no piece in $\pp$ is stabilized by the whole group $G$;
 likewise no point in $\free$ is fixed by the whole group $G$;
 \item[\textbf{(iii)}] the collection of subgroups
 $\calc(G)=\calc_1(G)
 \cup \calc_2 (G)$  satisfies the ascending chain condition.
 \end{itemize}
 Then one of the following three cases occurs:
  \begin{itemize}
  \item[(1)] $G$ splits over a
  $\calc(G)$-by-cyclic group;
  \item[(2)] $G$ can be represented as the fundamental group of a
  graph of groups whose vertex groups are of the form $\Stab(B)$
  or $\Stab(p)$ and
  edge groups are of the form $\Stab(B,p)$, $B\in \pp, p\in \free$;
  \item[(3)] the group $G$ has a $\calc_1(G)$-by-(free Abelian) subgroup
  of index at most 2.
  \end{itemize}
\end{theorem}

The following statement is a version of Theorem \ref{split1} for
finitely generated groups, using Theorem \ref{sela}.

\begin{theorem}\label{split5}
Let $G$ be a finitely generated group acting by isometries on a
tree-graded space $(\free ,\pp)$ such that properties \textbf{(i)}
and \textbf{(ii)} from Theorem \ref{split1} hold, and in addition
\begin{itemize}
\item[\textbf{(iii)}] the collection of subgroups
 $\calc_2(G)$ satisfies the ascending chain condition
 and every subgroup in $\calc_1(G)\cup \calc_3(G)$ is trivial.
 \end{itemize}
 Then one of the following three cases occurs:
  \begin{itemize}
  \item[(1)] $G$ splits over a
  $\calc_2(G)$-by-cyclic group or over an Abelian-by-cyclic group;
  \item[(2)] same as case (2) in Theorem
\ref{split1};
  \item[(3)] the group $G$ has a metabelian subgroup
  of index at most 2.
  \end{itemize}
\end{theorem}

\proof The statement would follow from Theorems \ref{1},
\ref{sela}, Proposition \ref{prop1}, and Remark \ref{rem89} if we
prove that in Case B.2 of the proof of Theorem \ref{1}, the action
is stable. But in that case, by Proposition \ref{B2}, the action
is with Abelian arc stabilizers and with trivial tripod stabilizers,
 since all subgroups in
$\calc_1(G)$ are trivial by our assumption. Hence we can apply
Lemma \ref{stable}.
\endproof

If instead of Theorem \ref{sela} we use Theorem \ref{gui}, then by
Theorem \ref{1}, Lemma \ref{stable} and Remark \ref{rem89}, we
obtain the following version of our theorem:

\begin{theorem}[Theorem \ref{splitg1}]\label{splitg}
Let $G$ be a finitely generated group acting by isometries on a
tree-graded space $(\free ,\pp)$ such that properties \textbf{(i)}
and \textbf{(ii)} from Theorem \ref{split1} hold, and in addition
\begin{itemize}
\item[\textbf{(iii)}] the subgroups in $\calc_2(G)$ are
 (finite of uniformly bounded size)-by-Abelian
 and the subgroups in $\calc_1(G)\cup \calc_3(G)$
 have uniformly bounded cardinality.
 \end{itemize}
 Then one of the following three cases occurs:
  \begin{itemize}
  \item[(1)] $G$ splits over a
  [(finite of uniformly bounded size)-by-Abelian]-by-(virtually cyclic)
  subgroup
  \item[(2)] same as case (2) in Theorem
\ref{split1};
\item[(3)] the group $G$ has a subgroup
of index at most 2 which is a [(finite of uniformly bounded
size)-by Abelian]-by-(free Abelian) subgroup.
 \end{itemize}
\end{theorem}

\begin{remark}\label{remark459} The proofs of Theorems
\ref{split1}, \ref{split5}, \ref{splitg} show that in case (2) of
these theorems, the group splits as a non-trivial amalgamated
product or HNN extension with vertex subgroup of the form
$\Stab(B)$.
\end{remark}

\section{Applications: relatively hyperbolic groups}\label{srelhyp}

\subsection{Asymptotic cones}\label{ac}

Let $I$ be an arbitrary countable set. Recall that a
\textit{non-principal ultrafilter} $\omega$ over $I$ is a finitely
additive measure on the class $\mathcal{P}(I)$ of subsets of $I$
such that each subset has measure either $0$ or $1$ and all finite
sets have measure 0. Since we only use non-principal ultrafilters,
the word non-principal will be omitted throughout the paper.

If a statement $P(i)$ holds for all $i$ from a set $J$ such that
$\omega (J)=1$, then we say that $P(i)$ holds {\em $\omega$--a.s.}.

\begin{remark}\label{udisj} The definition of an ultrafilter
immediately implies the following.  Let $P_1(i)$, $P_2(i)\dots $,
$P_m(i)$, $i\in I$, be statements such that for any $i\in I$ no two
of them can be true simultaneously. If the disjunction of these
statements holds $\omega$--a.s. then there exists $k\in \{ 1,2,\dots
, m\}$ such that $\omega$--a.s. $P_k(i)$ holds and all $P_j(i)$ with
$j\neq k$ do not hold.
\end{remark}

Given a sequence of sets $(X_n)_{n\in I}$ and an ultrafilter
$\omega$, the {\em ultraproduct corresponding to $\omega$}, $\Pi
X_n/\omega$, consists of equivalence classes of sequences
$(x_n)_{n\in I}$, $x_n\in X_n$, where two sequences $(x_n)$ and
$(y_n)$ are identified if $x_n=y_n$ $\omega$--a.s. The equivalence
class of a sequence $x=(x_n)$ in $\Pi X_n/\omega$ is denoted {either
by $x^\omega$ or by $(x_n)^\omega$}. In particular, if all $X_n$ are
equal to the same $X$, the ultraproduct is called the {\em
ultrapower} of $X$ and it is denoted  by $\Pi X/\omega$.

If $G_n$, $n\ge 1$, are groups then $\Pi G_n/\omega$ is again a
group with the multiplication
$(x_n)^\omega(y_n)^\omega=(x_ny_n)^\omega$.

\begin{lemma} \label{D} Let $\omega$ be an ultrafilter and let $(X_i)$ be a
sequence of sets of cardinality at
  most $D$ $\omega$--a.s. Then the ultraproduct $\Pi
  X_i/\omega$ contains at most $D$ elements.
\end{lemma}

\proof Consider $X_i=\{x_i^1,..., x_i^D \}$. If the cardinality of
$X_i$ is strictly less that $D$ then the last element is repeated
as many times as necessary. Let $x^j_\omega =(x_i^j)^\omega$ for
$j=1,2,...,D$. For every $y_\omega=(y_i)^\omega$ in $\Pi
  X_i/\omega$, there exists $j\in \{ 1,2,...,D\}$ such that
  $\omega$--a.s. $y_i=x_i^j$, hence $y_\omega = x^j_\omega$. Thus $\Pi
  X_i/\omega = \{x^1_\omega,..., x^D_\omega\}$.\endproof

\begin{remark} Lemma \ref{D} is a simple corollary of the
well known Los' theorem: if a first order property holds in
$\omega$--almost all $X_i$ then it holds in the ultraproduct $\Pi
X_i/\omega$: consider the formula $\exists x_1,...,x_D \forall y
(y=x_1 \vee \dots \vee y=x_D)$.
\end{remark}

For every sequence of points $(x_n)_{n\in I}$ in a topological space
$X$, its $\omega$--\textit{limit} $\lm_\omega x_n$ is a point $x$ in
$X$ such that every neighborhood $U$ of $x$ contains $x_n$
$\omega$--a.s.

\begin{itemize}
\item Suppose that the metric space $X$ is Hausdorff. For every sequence $(x_n)$ in $X$,
if the $\omega$--limit $\lm_\omega x_n$ exists, then it is unique.

\item Every sequence of elements in a compact space has an
$\omega$--limit \cite{Bou}.
\end{itemize}

\begin{definition}[$\omega$--limit of metric spaces]

Let $(X_n,\dist_n)$, $n\in I$, be a sequence of metric spaces and
let $\omega$ be an ultrafilter over $I$. Consider the ultraproduct
$\Pi X_n/\omega$. For every two points $x=(x_n)^\omega,
y=(y_n)^\omega$ in $\Pi X_n/\omega$ let
$$D(x,y)=\lm_\omega \dist_n(x_n,y_n)\, .$$

Consider an {\em observation point} $e=(e_n)^\omega$ in $\Pi
X_n/\omega$ and define $\Pi_e X_n/\omega$ to be the subset of $\Pi
X_n/\omega$ consisting of elements which are finite distance from
$e$ with respect to $D$. The function $D$ is a pseudo-metric on
$\Pi_e X_n/\omega$, that is, it satisfies the triangle inequality
and the property $D(x,x)=0$, but for some $x\ne y$ the number
$D(x,y)$ can be $0$.

The {\em $\omega$--limit} $\lio{X_n}_e$ {\em of the metric spaces}
$(X_n,\dist_n)$ {\em relative to the observation point} $e$ is the
metric space obtained from $\Pi_e X_n/\omega$ by identifying all
pairs of points $x,y$ with $D(x,y)=0$. The equivalence class of a
sequence $(x_n)$ in $\lio{X_n}_e$ is denoted by $\lio{x_n}$.
\end{definition}

Note that if $e,e'\in \Pi X_n/\omega$ and $D(e,e')<\infty$ then
$\lio{X_n}_e=\lio{X_n}_{e'}$.

\begin{definition}[asymptotic cone] Let $(X,\dist)$ be a metric
space, $\omega$ be an ultrafilter over a set $I$, $e=(e_n)^\omega$
be an observation point. Consider a sequence of numbers
$d=(d_n)_{n\in I}$ called {\em scaling constants} satisfying
$\lm_\omega d_n=\infty$.

The $\omega$--limit $\lio{X,\frac{1}{d_n}\dist}_e$ is called an
{\em asymptotic cone of $X$.} It is denoted by $\co{X;e,d}$.

Note that if $X$ is a group $G$ endowed with a word metric then
$\Pi_1 G/\omega$ is a subgroup of the ultrapower of $G$.
\end{definition}

\begin{definition}
For a sequence $(A_n), n\in I,$ of subsets of $(X,\dist)$ we denote
by $\lio{A_n}$ the subset of $\co{X; e, d}$ that consists of all the
elements $\lio{x_n}$ such that $x_n\in A_n$ $\omega$--a.s. Notice
that if $\lim_\omega \frac{\dist(e_n,A_n)}{d_n}=\infty $ then the
set $\lio{A_n}$ is empty.
\end{definition}

\medskip

\textit{Properties of asymptotic cones:}

\begin{enumerate}
\item Any asymptotic cone of a metric space is a complete metric space \cite{VDW}.
The same proof gives that $\lio{A_n}$ is always a closed subset of
the asymptotic cone $\co{X;e,d}$.

\item {Let $G$ be a finitely generated group endowed with a word metric.
The group $\Pi_1 G/\omega$ acts on $\co{G;1,d}$ transitively by
isometries:}
$$(g_n)^\omega \lio{x_n}=\lio{g_nx_n}.$$

{Given an arbitrary sequence of observation points $x$, the group
$x^\omega (\Pi_1 G/\omega) (x^\omega)\iv$ acts transitively by
isometries on the asymptotic cone $\co{G;x,d}$. In particular, every
asymptotic cone of $G$ is homogeneous.}

More generally if a group $G$ acts by isometries on a metric space
$(X,\dist)$ and there exists a bounded subset $B \subset X$ such
that $X=G\cdot B$, then all asymptotic cones of $X$ are
homogeneous metric spaces.
\end{enumerate}

\me

\textit{Convention:} When we consider an asymptotic cone of a
finitely generated group{, unless otherwise stated, we shall assume}
that the observation point $e$ is $(1)^\omega$.

\subsection{Asymptotically tree-graded metric
spaces}\label{satg}

\begin{definition}[asymptotically tree-graded spaces]\label{asco}

Let $(X,\dist)$ be a metric space and let $\mathcal{A}=\{ A_i \mid
i\in I \}$ be a collection of subsets of $\free$. In every
asymptotic cone $\co{\free;e,d}$, we consider the collection of
subsets
$$
\mathcal{A}_\omega = \left\{ {\lio{ A_{i_n}}}\mid (i_n)^\omega\in
\Pi I/\omega \hbox{ such that the sequence } \left(
\frac{\dist(e_{n},A_{i_n})}{d_{n}}\right) \hbox{ is
bounded}\right\}\, .
$$

We say that $X$ is \textit{asymptotically tree-graded with respect
to} $\aaa$ if every asymptotic cone $\co{X; e,d}$ is tree-graded
with respect to $\aaa_\omega $.

This notion is a generalization, in the setting of metric spaces, of
the usual notion of (strongly) relatively hyperbolic group.
\end{definition}

There is no need to vary the ultrafilter in Definition \ref{asco}:
if a space is tree-graded with respect to a collection of subsets
for one ultrafilter, it is tree-graded for any other with respect to
the same collection of subsets \cite[Corollary 4.30]{DS}.

We need the following facts from \cite{DS}.

\begin{lemma}[Theorem 4.1, Remark 4.2, (2), in \cite{DS}] \label{alphas}
Let $(X,\dist)$ be a geodesic metric space and let $\aaa=\{ A_i \mid
i\in I \}$ be a collection of subsets of $X$. If the metric space
$X$ is asymptotically tree-graded with respect to $\aaa$ then the
following properties are satisfied:
\begin{itemize}
  \item[$(\alpha_1)$] for every $\delta >0$ the diameters
of the intersections $\nn_{\delta}(A_i)\cap \nn_{\delta}(A_j)$ are
uniformly bounded for all $i\ne j$;

\item[$(\alpha_2)$] for every $L\geq 1$, $C\ge 0$, and $\theta
\in \left[ 0, \frac{1}{2} \right)$ there exists $M>0$ such that
for every $(L,C)$--quasi-geodesic $\q$ defined on $[0,\ell ]$ and
every $A\in \aaa$ such that
  $\q(0),\q(\ell )\in \nn_{\theta\ell/L}
  (A)$ we have $\q([0, \ell ])\cap \nn_{M}
  (A)\neq \emptyset$.
\end{itemize}
\end{lemma}

\medskip

\noindent \textit{Notation:}  Let $M(L,C)$ be the maximum between
the constant defined in \cite[Definition 4.20]{DS} and the
constant given by $(\alpha_2)$ for $\theta =\frac{1}{3}$.

\medskip

\begin{definition}\label{defsat} Let $\q$ be an
$(L,C)$--quasi-geodesic in $X$. The \textit{saturation of }$\g$,
denoted by $\Sat(\g)$, is the union of $\g$ and all the sets from
$\aaa$ whose $M(L,C)$--tubular neighborhood crosses~$\g$.
\end{definition}

\begin{lemma}[Lemma 4.15 in \cite{DS}]
\label{815} Let $X$ be a geodesic metric space which is
asymptotically tree-graded with respect to a collection of subsets
$\aaa$. For every $L\geq 1$ and $C\geq 0$, there exists $t\geq 1$
such that for every $d\geq 1$ and for every $A\in \aaa$, every
$(L,C)$--quasi-geodesic joining two points in $\nn_d(A)$ is
contained in $\nn_{td}(A)$.
\end{lemma}

\begin{lemma}\label{813}
Let $X$ be a metric space asymptotically tree-graded with respect
to a collection of subsets $\aaa$. For every $L\geq 1,\, C\geq
0,\, M\geq M(L,C)$ and $\delta >0$, there exists $D_0>0$ such that
the following holds. Let $A\in \aaa$ and let $\q_i : [0,\ell_i]
\to X\, ,\, i=1,2,$ be two $(L,C)$--quasi-geodesics with the
 two respective start points $a_i$ in $\nn_M (A)$, and such that the
 diameter of $\q_i \cap \overline{\nn}_M(A)$ does not exceed $\delta $
 for $i=1,2$.

If $\dist (a_1,a_2)\geq D_0$ then $\q_1\sqcup [a_1,a_2]\sqcup
\q_2$ is an $(L+C+1,C_1)$--quasi-geodesic, where $C_1=C_1(D_0,
\delta , C)$.
\end{lemma}

\proof {According to \cite[Lemma 4.19]{DS}, $\q_1\sqcup [a_1,a_2]$
and $[a_1,a_2]\sqcup \q_2$ are $(L_1,C_1)$--quasi-geodesics. It
remains to prove that for every $t\in [0,\ell_1]$ and $s\in
[0,\ell_2]$,
$$
\dist (\q_1(t),\q_2(s))+O(1)\geq \frac{1}{L}(t+s+\dist (a_1,a_2))\,
.$$}

Lemma 4.28 and Corollary 8.14 in \cite{DS} imply that a geodesic
$[\q_1(t),\q_2(s)]$ contains two points $a_1'$ and $a_2'$ at
distance $\varkappa$ of $a_1$ and $a_2$, respectively, with
$\varkappa =\varkappa (X, \aaa)$. Then $\dist
(\q_1(t),\q_2(s))\geq \dist (\q_1(t),a_1)+\dist (a_1,a_2)+\dist
(a_2,\q_2(s))-4\varkappa\geq \frac{1}{L}(t+s)-2C+\dist
(a_1,a_2)-4\varkappa$.\endproof

\begin{lemma}[Lemmas 4.26 and 4.28 in \cite{DS}]\label{pl} Let
$\bigcup_{i=1}^m \q_i$ be a polygonal line composed of
$(L,C)$--quasi-geodesics.
    \begin{enumerate}
    \item \textbf{(uniform property $(\alpha_2)$ for saturations of polygonal
    lines)} For every $\lambda \geq 1$, $\kappa \geq 0$ and $\theta \in \left[
0, \frac{1}{2} \right) $ there exists $R$ such that for every
$(\lambda ,\kappa )$--quasi-geodesic $\cf:[0,\ell ]\to X $ joining
two points in $\nn_{\theta \ell /\lambda} \left(\bigcup_{i=1}^m
\Sat\left(\q_i\right) \right)$, we have $\cf([0,\ell ])\cap
\nn_{R} (\bigcup_{i=1}^m \Sat\left(\q_i\right))\neq \emptyset$ (in
particular, the constant $R$ does not depend on $\q_i$, only on
$m$).
    \item \textbf{(uniform quasi-convexity of saturations of polygonal
    lines)}  For every $\lambda
\geq 1$, $\kappa \geq 0$, there exists $\tau $ such that for every
$R\ge1$, the union of saturations $\bigcup_{i=1}^m
\Sat\left(\q_i\right)$ has the property that every $(\lambda ,
\kappa )$--quasi-geodesic $\cf$ joining two points in its
$R$--tubular neighborhood is entirely contained in its $\tau
R$--tubular neighborhood.
    \item  For every $\delta >0$ and every $A\in \aaa$ such that $A\not\subset\bigcup_{i=1}^m \Sat
\left( \q_i \right)$, the intersection $\nn_\delta(A)\cap \nn_\delta
\left( \bigcup_{i=1}^m \Sat \left( \q_i \right) \right)$ has
diameter bounded uniformly in $A$, $\q_1,\ldots,\q_m$.

\item For every $R>0$ and $ \delta
>0$ there exists $\varkappa
>0$ such that if $A,B\in \aaa, A\cup B\subset \bigcup_{i=1}^m
\Sat(\q_i), A\neq B$, the following holds. Let $a\in \nn_R(A)$ and
$b\in \nn_R(B)$ be two points that can be joined by a quasi-geodesic
${\mathfrak p}$ such that ${\mathfrak p} \cap \nn_R(A)$ and
${\mathfrak p} \cap \nn_R(B)$ has diameter at most $\delta$. Then
$\{ a,b \}\subset \nn_{\varkappa} \left( \bigcup_{i=1}^m \q_i
\right) $.
\end{enumerate}
\end{lemma}

\begin{lemma}[Corollary 8.14 in \cite{DS}]\label{ends} For
every $L\geq 1$, $C\geq 0$, $M\geq M(L,C)$ and $\delta>0$ there
exists $D_1>0$ such that the following holds. Let $A\in \aaa$ and
let $\q_i\colon [0,\ell_i]\to X$, $i=1,2$, be two
$(L,C)$--quasi-geodesics with one common endpoint $b$ and the
other two respective start points $a_i\in \nn_M (A)$, such that
the diameter of $\q_i\cap \overline{\nn}_M (A)$ does not exceed
$\delta$. Then $\dist(a_1,a_2)\leq D_1$.
\end{lemma}

In what follows we consider the image of a quasi-geodesic $\q
:[0,\ell ] \to X$ endowed with the order from $[0,\ell]$.

Given a subset $A\in X$ intersecting $\q [0,\ell ]$ we call
\textit{entrance point} of $\q$ into $A$ the image $\q(t)$ of the
smallest number $t\in [0,\ell ]$ such that $\q(t)\in A$ and $\q
(t-1)$ or $\q (0)$ if $0\leq t<1$ is in the complementary of $A$.

We call \textit{exit point} of $\q$ from $A$ the image $\q(s)$ of
the largest number $s\in [0,\ell ]$ such that $\q(s)\in A$ and $\q
(s+1)$ or $\q (\ell)$ if $\ell \geq s>\ell-1$ is in the
complementary of $A$.

\begin{lemma}\label{mm}
Let $\bigcup_{i=1}^m \q_i$ be a polygonal line composed of
$(L,C)$--quasi-geodesics. Let $R$ be the constant given by Lemma
\ref{pl}, (1), for $\lambda \geq 1$, $\kappa \geq 0$ and $\theta
=\frac{1}{3}$. For every $R\leq R_1<R_2$ and every $(\lambda ,
\kappa)$ --quasi-geodesic $\q :[0, \ell ] \to X$ with $\q (0)$ in
$\mathcal{N}_1=\nn_{R_1} \left(\bigcup_{i=1}^m
\Sat\left(\q_i\right) \right)$ and $\q (\ell)$ outside
$\mathcal{N}_2=\nn_{R_2} \left(\bigcup_{i=1}^m
\Sat\left(\q_i\right) \right)$ the exit point $\q(t_1)$ from
$\mathcal{N}_1$ and the exit point $\q(t_2)$ from $\mathcal{N}_2$
satisfy $ t_2-t_1\leq 3\lambda (R_2+\lambda +\kappa )$. In
particular the two exit points are at uniformly bounded distance.
\end{lemma}

\proof If on the contrary $t_2-t_1> 3\lambda (R_2+\lambda +\kappa
)$, this together with Lemma \ref{pl}, (1), applied to $\q|_{[t_1+1,
t_2]}$ would contradict the fact that $\q (t_1)$ is an exit point
from $\mathcal{N}_1$.\endproof

\begin{lemma}\label{endsg}
The statement in Lemma \ref{ends} holds with $A$ replaced by the
saturation of a third $(L,C)$--quasi-geodesic $\pgot_0$ and $M\geq
\max (M(L,C), R)$, where $R$ is the constant given in Lemma
\ref{pl}, (1), for $(\lambda, \kappa)=(L,C)$, $\theta=\frac{1}{3}$
and $m=1$.
\end{lemma}

\begin{center}
\begin{figure}[!ht]
%TeXCAD Options
%\grade{\on}
%\emlines{\off}
%\epic{\off}
%\beziermacro{\on}
%\reduce{\on}
%\snapping{\off}
%\quality{8.00}
%\graddiff{0.01}
%\snapasp{1}
%\zoom{4.0000}
\unitlength .7mm % = 2.85pt
\linethickness{0.4pt}
\ifx\plotpoint\undefined\newsavebox{\plotpoint}\fi % GNUPLOT compatibility
\begin{picture}(120.75,95.75)(-45,0)
\qbezier(58.75,90.25)(52.63,89.13)(56,83.5)
\qbezier(56,83.5)(59.13,78.13)(53.75,77.25)
\qbezier(53.75,77.25)(47.75,76.75)(51.75,70.25)
\qbezier(51.75,70.25)(56.5,64.25)(49.25,63.25)
\qbezier(49.25,63.25)(43.75,62.75)(47.25,57.25)
\qbezier(47.25,57.25)(53,49.5)(44.75,49.75)
\qbezier(44.75,49.75)(38.25,49.75)(41.75,40.75)
\qbezier(41.75,40.75)(44.5,33)(39.25,33.25)
\qbezier(39.25,33.25)(34.13,33.38)(36.5,25)
\qbezier(58.5,90.25)(64.5,87.13)(61.5,81.5)
\qbezier(61.5,81.5)(59.75,78.63)(65,73.25)
\qbezier(65,73.25)(69.75,67.63)(67.5,65.5)
\qbezier(67.5,65.5)(64.88,63.38)(70.75,55.75)
\qbezier(70.75,55.75)(76.88,47.75)(74.5,45.75)
\qbezier(74.5,45.75)(70.25,45.25)(78,34.75)
\qbezier(78,34.75)(84,26.75)(81,26.75) \put(58.75,90){\circle*{2}}
\put(36.5,25.25){\circle*{2}} \put(81.5,26.75){\circle*{2}}
\qbezier(8,27.5)(9.38,37.13)(33.25,40.25)
\qbezier(33.25,40.25)(68.88,44.88)(94,41)
\qbezier(94,41)(116.63,38)(118.75,27)
\qbezier(118.75,27)(120.75,16.88)(57.75,7.25)
\qbezier(57.75,7.25)(19.88,1.5)(13.5,13.75)
\qbezier(13.5,13.75)(7.25,24.38)(8,27.5)
\put(61,95.75){\makebox(0,0)[cc]{$b$}}
\put(31,21.5){\makebox(0,0)[cc]{$a_1$}}
\put(85,24.25){\makebox(0,0)[cc]{$a_2$}}
\put(72.5,65.25){\makebox(0,0)[cc]{$\q_2$}}
\put(45,69){\makebox(0,0)[cc]{$\q_1$}}
\put(55,15.5){\makebox(0,0)[cc]{$\Sat(\pgot_0)$}}
\end{picture}
%\centering
\caption{Lemma \ref{endsg}.} \label{fig2}
\end{figure}
\end{center}

\proof  We have $a_i\in \nn_M \left( \Sat \pgot_0 \right)$. There
are three cases.

\medskip

 \textbf{Case 1.} Both $a_i$ are in $\nn_M \left( \pgot_0 \right)$.
Let $a_i'$ be a point in $\pgot_0$ such that $\dist (a_i,a_i')\leq
M$. Let $\pgot:[0,\ell ]\to X$ be a sub-quasi-geodesic of $\pgot_0$
of endpoints $a_1'$ and $a_2'$. Note that { the number $\ell$ is of
the order of $\dist (a_1,a_2)$. }

By Lemma \ref{pl}, (2), $\pgot [0,\ell ] \subset \nn_\tau  (\Sat
(\q_1))\cup \nn_\tau (\Sat (\q_2))$. Let $m=\pgot (\ell/2)$. If
$m\in \nn_\tau (\q_1)\cup \nn_\tau ( \q_2)$ and $\ell$ is large
enough this and Lemma \ref{pl}, (1), contradict the hypothesis
that the diameter of $\q_i\cap \overline{\nn}_M (\Sat (\pgot ))$
is at most $\delta$.

Assume that $m\in \nn_\tau (A)$ with $A\subset \Sat (\q_1)$. Note
that the entrance point $e_1$ of $\pgot$ in $\onn_\tau (A)$ and the
entrance point $e_2$ of $\q_1$ in $\nn_M (A)$ are at distance
$O(1)$, by Lemmas \ref{ends} and \ref{mm}. If the distance from
$e_2$ to $a_1$ is too large then Lemma \ref{pl}, (1), allows { us }
to find a point in $\q_1 \cap \nn_R(\Sat (\pgot))$ at distance $\gg
\delta$ from $a_1$. Thus the diameter of $\{ a_1, e_1, e_2\}$ is
$O(1)$. If $\ell$ is large enough then $\dist (m, e_1)$ is { larger
than any constant fixed in advance}, and by Lemma \ref{pl}, (3),
$A\subset \Sat (\pgot ). $ It follows that if $e_2'$ is the exit
point of $\q_1$ from $\nn_M(A)$, then it is at distance at most
$\delta$ from $a_1$. Therefore $\dist (e_2', m)$ is also { larger
than any constant fixed in advance}.

Let $e_1'$ be the exit point of $\pgot$ from $\onn_\tau (A)$. We
have that $m$ is between $e_1$ and $e_1'$, hence { we can assume
that $\dist (e_2', e_1')$ is large enough}. The sub-quasi-geodesic
of $\q_1$ between $b$ and $e_2'$, together with $[e_2',e_1']$ and
the sub-quasi-geodesic of $\pgot$ between $e_1'$ and $a_2$ compose
a quasi-geodesic, by Lemma \ref{813}. According to Lemma \ref{pl},
(2), this quasi-geodesic is contained in $\nn_{\tau'} (\Sat
(\q_2))$. In particular $\nn_\tau (A)$ intersects $\nn_{\tau'}
(\Sat (\q_2))$ in a set of diameter at least $\dist (e_1',e_2')$,
so by Lemma \ref{pl}, (3), it is contained in $\Sat (\q_2)$. The
intersection of $\q_2$ with $\nn_M(A)$ is contained in its
intersection with $\nn_M(\Sat (\pgot_0 ))$, therefore it is in
$B(a_2, \delta)$. We have thus obtained that both $a_1$ and $a_2$
are at distance $O(1)$ from the entrance points of $\q_1$ and
respectively $\q_2$ in $\nn_M(A)$. We may then use Lemma
\ref{ends}.

\medskip

 \textbf{Case 2.} Assume that $a_1\in \nn_M(A)$ for some $A\subset \Sat
 (\pgot_0 )$ while $a_2\in \nn_M(\pgot_0)$. Let $a_2'$ be a point in $\pgot_0\cap B(a_2, M)$.
 Assume that $\nn_M(A)$ intersects the sub-quasi-geodesic $\pgot_1$ of $\pgot_0$
 between $\pgot_0(0)$ and $a_2'$. Let $e$ be the exit point of
 $\pgot_1$ from $\onn_M(A)$,
 and let $\pgot$ be the sub-quasi-geodesic of $\pgot_1$ between $e$ and $a_2'$.

 The union $\q_1'=\q_1\sqcup [a_1,e]$ is a quasi-geodesic. Indeed, if $\dist
 (a_1,e)$ is large enough, this follows from Lemma \ref{813} while
 if $\dist (a_1,e)\leq D_0$ it is obvious. Also if the intersection
 of $\q_1'$ with $\Sat (\pgot )$ contains a point too far from
 $e$, then $A$ has a large intersection with $\Sat (\pgot)$,
 hence it is contained in it, and this contradicts the choice of
 $\pgot
 $. Thus, the intersection  $\q_1'\cap \Sat (\pgot )$ is
 at distance $O(1)$ from $e$. We may then apply Step 1 to $\q_1'$
 and $\q_2$ and deduce that $\dist (e,a_2)$ is $O(1)$. Then we apply
 Lemma \ref{ends} to $\q_1, \q_2$ and $A$, and deduce that $\dist (a_1,a_2)$
 is $O(1)$.

\medskip

 \textbf{Case 3.} Assume that $a_i\in \nn_M(A_i)$ for some $A_i\subset \Sat
 (\pgot_0 )$, $i=1,2$. Without loss of generality we may suppose that $A_2$ intersects $\pgot_0$ between
  its exit point from $\onn_M(A_1)$ and the end of $\pgot_0$. Let $e$ be the exit point of $\pgot_0$ from
 $\onn_M(A_1)$ and let $\pgot$ be the sub-quasi-geodesic of $\pgot_0$ between $e$ and its end. As above, $\q_1'=\q_1\sqcup [a_1,e]$ is a
 quasi-geodesic with the property that the intersection  $\q_1'\cap \Sat (\pgot )$ is
 at distance $O(1)$ from $e$.
  Step 2 for $\q_1'$, $\q_2$ and $\Sat (\pgot)$ implies that $\dist (e,
  a_2)$ is $O(1)$. Lemma \ref{ends} applied to $\q_1, \q_2$ and
  $A_1$ implies that $\dist (a_1,a_2)$ is $O(1)$.\endproof

\begin{lemma}\label{813g}
The statement in Lemma \ref{813} holds with $A$ replaced by the
saturation of a third $(L,C)$--quasi-geodesic $\pgot_0$ and $M\geq
\max (M(L,C), R)$, where $R$ is the constant given in Lemma
\ref{pl}, (1), for $(\lambda, \kappa)=(L,C)$, $\theta=\frac{1}{3}$
and $m=1$.
\end{lemma}

\begin{center}
\begin{figure}[!ht]
%TeXCAD Options
%\grade{\on}
%\emlines{\off}
%\epic{\off}
%\beziermacro{\on}
%\reduce{\on}
%\snapping{\off}
%\quality{8.00}
%\graddiff{0.01}
%\snapasp{1}
%\zoom{4.0000}
\unitlength .7 mm % = 2.85pt
\linethickness{0.4pt}
\ifx\plotpoint\undefined\newsavebox{\plotpoint}\fi % GNUPLOT compatibility
\begin{picture}(120.75,110.25)(-45,0)
\put(36.5,25.25){\circle*{2}} \put(81.5,26.75){\circle*{2}}
\qbezier(8,27.5)(9.38,37.13)(33.25,40.25)
\qbezier(33.25,40.25)(68.88,44.88)(94,41)
\qbezier(94,41)(116.63,38)(118.75,27)
\qbezier(118.75,27)(120.75,16.88)(57.75,7.25)
\qbezier(57.75,7.25)(19.88,1.5)(13.5,13.75)
\qbezier(13.5,13.75)(7.25,24.38)(8,27.5)
\put(31,21.5){\makebox(0,0)[cc]{$a_1$}}
\put(85,24.25){\makebox(0,0)[cc]{$a_2$}}
\put(55,15.5){\makebox(0,0)[cc]{$\Sat(\pgot_0)$}}
\qbezier(12,107)(15.75,103.13)(14.5,98.75)
\qbezier(14.5,98.75)(13.25,96.88)(16,93.5)
\qbezier(16,93.5)(18.88,90.13)(18.25,86.25)
\qbezier(18.25,86.25)(16.25,79.5)(21.25,76.75)
\qbezier(21.25,76.75)(23.25,75.38)(23.25,69.5)
\qbezier(23.25,69.5)(23,60.13)(25.75,61.25)
\qbezier(25.75,61.25)(27.13,62.63)(28,53.5)
\qbezier(28,53.5)(27.88,47.88)(30.25,46.75)
\qbezier(30.25,46.75)(34.13,43.5)(32.5,38.25)
\qbezier(32.5,38.25)(31.13,33.5)(36.25,26.75)
\qbezier(105.25,110.25)(107.25,103.13)(102.25,100.5)
\qbezier(102.25,100.5)(96.38,97.38)(100,91.75)
\qbezier(100,91.75)(103,87.25)(97,81.75)
\qbezier(97,81.75)(89.75,74.38)(93.5,70.5)
\qbezier(93.5,70.5)(94.88,69.25)(90.75,61)
\qbezier(90.75,61)(85.88,51.75)(87.5,48.5)
\qbezier(87.5,48.5)(88.5,45.13)(84.5,38.25)
\qbezier(84.5,38.25)(79.5,30.13)(81.5,27.5)
\put(22.75,81.25){\makebox(0,0)[cc]{$\q_1$}}
\put(99,78.75){\makebox(0,0)[cc]{$\q_2$}}
%\emline(36.25,25.5)(81.5,27)
\multiput(36.25,25.5)(1.0055556,.0333333){45}{\line(1,0){1.0055556}}
%\end
\put(53.25,31.25){\makebox(0,0)[cc]{$\ge D_0$}}
\end{picture}
%\centering
\caption{Lemma \ref{813g}.} \label{fig3}
\end{figure}
\end{center}

\proof \textbf{Step 1.} We prove that $\tq_1= \q_1\sqcup
[a_1,a_2]$ is a quasi-geodesic (the same argument implies that
$\tq_2=[a_1,a_2]\sqcup \q_2$ is a quasi-geodesic). Let $t\in
[0,\ell_1]$ and let $x\in [a_1,a_2]$. Consider a geodesic
$[\q_1(t),x]$. Its entrance point in $\onn_M(\Sat (\pgot ))$ is at
distance $O(1)$ from $a_1$, by Lemma \ref{endsg}. It follows that
$\dist (\q_1(t), x )= t+\dist (a_1,x) +O(1)$.

\medskip

\textbf{Step 2.} We now prove the statement in the Lemma. Let
$t\in [0, \ell_1]$ and $s\in [0,\ell_2 ]$. Let $\g $ be a geodesic
joining $\q_1(t)$ and $\q_2(s)$. By Lemma \ref{pl}, $\tq_1$ is
contained in the $\tau$-tubular neighborhood of $\Sat (\q_2) \cup
\Sat (\g)$. In particular the point $a_1$ is in the same tubular
neighborhood.

\medskip

\textbf{Step 2.a.} Assume that $a_1$ is in the $\tau$-tubular
neighborhood of $\Sat (\q_2)$. If $a_1$ is at distance less than
$\tau$ from $\q_2$, then by the hypothesis on $\q_2$ and Lemma
\ref{mm}, $a_1$ is at uniformly bounded distance from $a_2$. Thus
if $D_0$ is large enough this case cannot occur.

Suppose that $a_1$ is in the $\tau$-tubular neighborhood of
$A\subset \Sat (\q_2)$. Then the portion of $\tq_2$ between $a_1$
and a point in $\q_2\cap \nn_M(A)$ is contained in
$\nn_{\tau'}(A)$, hence the same holds for $[a_1,a_2]$. If $D_0$
is large enough this implies that $A\subset \Sat (\pgot_0 )$.

We may conclude in this case that $\q_1 \sqcup [a_1,a_2]\sqcup
\q_2$ is a quasi-geodesic by Lemma \ref{813}.

\medskip

\textbf{Step 2.b.} Suppose that $a_1$ is in $\nn_\tau (A)$ for
some $A\subset \Sat (\g )$. By Lemma \ref{ends} the exit point of
$\g$ from $\nn_M(A)$ is at distance $O(1)$ from a point $x$ in
$\tq_2$. If $x$ is not on $[a_1,a_2]$ then $x\in \q_2$, and it
follows that a part of $\tq_2$ containing $[a_1,a_2]$ is in some
$\nn_{\tau'} (A)$, whence $A\subset \Sat (\pgot_0)$ if $D_0$ is
large enough. Lemma \ref{813} then implies that $\q_1 \sqcup
[a_1,a_2]\sqcup \q_2$ is a quasi-geodesic.

Suppose then that $x\in [a_1,a_2]$.  Then $\dist (\q_1(t),
\q_2(s))=\dist (\q_1(t), x)+\dist (x, \q_2(s))+O(1) $. Since
$\tq_1$ is a quasi-geodesic, it follows that $\dist(\q_1(t),
x)\geq \frac{1}{L_1}(t+\dist(a_1,x))-C_1$ for some $L_1\geq 1$ and
$C_1\geq 0$. Also $\tq_2$ is a quasi-geodesic, therefore
$\dist(\q_2(s), x)\geq \frac{1}{L_1}(s+\dist(a_2,x))-C_1$. We
conclude that
$$
\dist (\q_1(t), \q_2(s)) +O(1)\geq
\frac{1}{L_1}(t+s+\dist(a_1,a_2))\, .
$$

Now assume that $a_1$ is at distance at most $\tau$ from a point
in $\g$.  Then the argument above can be repeated with $x=a_1$.

\endproof

\subsection{Relatively hyperbolic groups}\label{srhg}

The following result of the authors and Denis Osin \cite{DS} shows
that Cayley graphs of relatively hyperbolic groups are
asymptotically tree-graded metric spaces.

\begin{theorem} [Theorem 8.5 and Appendix of \cite{DS}] \label{tgr} A group $G$ is relatively
      hyperbolic with respect to its finitely generated subgroups $H_1,...,H_m$ if and only
      if every asymptotic cone of $G$ is tree-graded with respect to
      the collection of $\omega$--limits of sequences of cosets
      $\gamma_nH_i$ ($\gamma_n\in G, i=1,...,m$).
\end{theorem}

{\bf From now on we fix an infinite finitely generated group $G$
that is relatively hyperbolic with respect to a finite collection
of finitely generated peripheral subgroups $\mh=\{H_1,...,H_m\}$,
$G\ne H_i$. Let ${\cal G}$ be the set of all left cosets of
$H_i$.}

\begin{definition}
Recall that $H_i$ are called {\em peripheral subgroups} of $G$,
subgroups of conjugates of $H_i$ are called {\em parabolic
subgroups} of $G$, conjugates of $H_i$ are called {\em maximal
parabolic subgroups}.
\end{definition}

Note  that a maximal parabolic subgroup $\gamma H_i \gamma\iv$ is
the stabilizer of a left coset $\gamma H_i$.

\medskip

\begin{lemma}\label{stabc} Let
$\ck=\co{G;x,d}$ be an asymptotic cone of $G$ and let
$\lio{\gamma_n H}$ be an ultralimit, where $\gamma_n\in G$ and $H$
is a peripheral subgroup such that $\lio{\gamma_n H}$ is a
non-empty subset of $\ck$. The stabilizer in $x^\omega (\Pi_1
G/\omega) (x^\omega)\iv$ of $\lio{\gamma_n H}$ is inside $\Pi
\left( \gamma_n H \gamma_n\iv \right) /\omega $.
\end{lemma}

\proof Let $\sss $ be the stabilizer in $x^\omega (\Pi_1 G/\omega)
(x^\omega)\iv$ of $\lio{\gamma_n H}$. For every $(g_n)^\omega \in
\sss$, $\lio{g_n \gamma_n H }=\lio{\gamma_n
  H}$, therefore $\omega$--a.s. $g_n \gamma_n H=\gamma_n
  H$, and $g_n \in \gamma_n H \gamma_n\iv$.
\endproof

\begin{lemma}\label{conj} Let $H, H'$ be peripheral subgroups of $G$, and let $g\in G$
 be such that $gH'\neq H$.
Then the subgroup $H\cap gH'g\iv$ is a conjugate of a subgroup
inside the ball
  $B(1,R)$ for some uniform constant $R$. In particular, the size
  of this subgroup is uniformly bounded.
\end{lemma}

\proof Let $a\in H, b\in gH'$ be such that
$\dist(a,b)=\dist(H,gH')$. In particular, if $\dist(H,gH')=0$ then
one can take $a=b$.

Let $x\in H\cap gH'g\iv$. Then $\dist(xa, xb)=\dist(H,gH')$. By
Lemma \ref{813}, if $\dist(a,xa)> D_0$ then the union $[b,a]\sqcup
[a,xa]\sqcup [xa,xb]$ is a $(2,C_1)$--quasi-geodesic (a geodesic
if $a=b$). By Lemma \ref{815}, this quasi-geodesic is in the
$t$-neighborhood of $gH'$ for some uniform constant $t$. Hence
$[a,xa]\subseteq \nn_t (H) \cap \nn_t(gH')$. By $(\alpha_1)$ the
distance $\dist(a,xa)$ is uniformly bounded. The same is obviously
true in the case when $\dist(a,xa)\le D_0$.  Hence $a\iv (H\cap
gH'g\iv)a$ is in a ball of uniformly bounded radius.
\endproof

\begin{lemma} \label{rem1} Let $K$ be a subgroup of $G$.
If $K$ contains a central element of infinite
order, or if $K$ does not contain free non-Abelian subgroups, then
either $K$ is virtually cyclic or $K$ is a parabolic subgroup.
\end{lemma}

\proof Suppose that $K$ has a central element $z$ of infinite
order. If no power of $z$ is in a parabolic subgroup then the
normalizer of $\la z \ra$ is virtually cyclic by \cite{Osin}, so
$K$ is virtually cyclic. If $z^n$ is inside a maximal parabolic
subgroup $H$ then for every $k\in K$, $k Hk\iv$ intersects $H$ in
an infinite set, since the intersection contains $\langle z^n
\rangle$ and $z$ is of infinite order. Hence $k H=H$ by Lemma
\ref{conj}, and $k\in H$. Thus $K\subseteq H$.

%$%$%0.0.0.1-256049849

Suppose that $K$ contains no free non-Abelian subgroups. By
\cite{Tu} any subgroup of a relatively hyperbolic group
  either has a free non-Abelian subgroup or it is elementary
  (that is, either virtually cyclic or parabolic). In particular $K$
  is elementary.
\endproof

In the case of non-parabolic subgroups the following can be said.

\begin{lemma}\label{lxie}
There exist a constant $R$ and a positive integer $m_0$ depending
on $G, H_1,..., H_m\, ,$ such that for any finite subset $S$ of
$G$ one of the following three cases occurs:
\begin{itemize}
    \item $S$ is contained in a
parabolic subgroup of $G$;
    \item $S$ is conjugate to a subset inside the ball $B(1, R)$;
    \item $S^{m_0} =\{ s_1...s_k \; ;\; s_i \in S\cup S\iv,\, k\leq
m_0\}$ contains a hyperbolic element.
\end{itemize}

In particular, every finite non-parabolic subgroup of $G$ is
conjugate to a subgroup inside $B(1,R)$.
\end{lemma}

\proof The argument in the proof of \cite[Lemma 5.3]{Xie}, relying
on a result in \cite{Koubi}, proves in fact the statement of the
lemma.
\endproof

\begin{lemma}\label{subgroups}
\begin{itemize}
  \item[(1)] If each $H_i$ is hyperbolic relative to $H_{i,j}$ ($j=1,...,s_i$)
then $G$ is hyperbolic relative to $H_{i,j},\, (i=1,...,n,
j=1,...,s_i)$.

  \item[(2)] Each peripheral subgroup $H_i$ is undistorted in $G$.

  \item[(3)] If $H$ is an undistorted subgroup of $G$ then it is hyperbolic
with respect to finitely many parabolic subgroups of $H$.

  \item[(4)] Suppose that none of the peripheral subgroups $H_i$ is
hyperbolic relative to proper subgroups. Then every automorphism
of $G$ permutes maximal parabolic subgroups.
\end{itemize}
\end{lemma}

\proof Part (1) is proved in \cite[Corollary 1.14]{DS}. Part (2)
immediately follows from Lemma \ref{815} and Theorem \ref{tgr}.
Part (3) is proved in \cite[Theorem 1.8]{DS}.

To prove part (4) let $\psi$ be an automorphism of $G$. Suppose
that the image by $\psi$ of a maximal parabolic subgroup
$gH_ig\iv$ is not parabolic. It follows that $\psi (H_i)$ is
likewise non-parabolic. Since $H_i$ is undistorted by part (2),
the subgroup $\psi(H_i)$ is undistorted as well, and we can apply
part (3) and conclude that $\psi(H_i)$ is relatively hyperbolic
with respect to some parabolic subgroups of it. It follows that
$H_i$ is also relatively hyperbolic with respect to some proper
subgroups of it, a contradiction.

Thus for every maximal parabolic subgroup $H = g H_i g\iv $, $\psi
(H)$ is parabolic, hence contained in some maximal parabolic
subgroup $H' $. The same argument applied to $H'$ and $\psi\iv$
implies that $H < \psi\iv (H') < H''$, where $H''$ is maximal
parabolic. Since $H$ is not finite (otherwise it would be
relatively hyperbolic with respect to the trivial subgroup), it
follows by Lemma \ref{conj} that $H=H''$, thus $\psi (H)=H'$.
\endproof

\medskip

The following lemma shows how a maximal parabolic subgroup acts
outside the left coset that it stabilizes.

\begin{lemma} \label{l6c}
Let $g$ be an element in a maximal parabolic subgroup $\gamma H
\gamma \iv$ and let $x$ be a point in $G\setminus \gamma H$. Let
$x_1$ be a nearest point projection of $x$ onto $\gamma H$. Then
there exists a uniform constant $C$ such that one of the following
 two situations occurs:
\begin{itemize}
    \item[(1)] $\dist(x_1, gx_1) \leq C$;
    \item[(2)] $\dist (x, gx)\geq \dist (x, \gamma H)+\frac{1}{2}
    \dist (x_1,gx_1)-C$.
\end{itemize}
\end{lemma}

\proof Consider the constant $D_0$ provided by Lemma \ref{813} for
$(L,C)=(1,0)$ and
 $M=\delta=M(1,0)$.
 If $\dist
(x_1, g\, x_1)\geq D_0$ then $[x,x_1]\sqcup [x_1,
 g\, x_1]\sqcup [g\, x_1, g\, x]$ is a
 $(2,C_1)$-quasi-geodesic, where $C_1$ is a uniform constant. This implies that
 $\dist(x, g\, x)\geq \dist (x,x_1)+\frac{1}{2} \dist (x_1,g\, x_1) -C_1$.

\endproof

\begin{lemma} \label{l6}
Let $\ck =\co{G; x,d}$ be an asymptotic cone of $G$. Then any
subgroup $\sss< x^\omega (\Pi_1 G/\omega) (x^\omega)\iv$  which
stabilizes a piece in $\ck$ and a point outside the piece, is
conjugate to a subgroup in $\Pi B(1,R)/\omega$ for some universal
constant $R=R(G)$. In particular, its size is bounded by a
universal constant $D=D(G)$.
\end{lemma}

\proof Let $\lio{\delta_iH}$ and $b=\lio{b_i}$ be the piece and
respectively the point outside the piece, fixed by $\sss$. Let
$g=(g_i)^\omega$ be an
 element in $\sss$. Note that since
$g\cdot \lio{\delta_iH}=\lio{\delta_iH}$, we have $g_i\in
\delta_iH\delta_i\iv$ $\omega$--a.s.

The distance from $b$ to $\lio{\delta_iH}$ is positive. Therefore
the distance from $b_i$ to $\delta_iH$ is at least $O(d_i)$. The
distance from $b_i$ to $g_ib_i$ is $o(d_i)$ since $g\cdot b=b$.
Let $c_i$ be a nearest point projection of $b_i$ onto $\delta_i
H$. By Lemma \ref{l6c}, $c_i\iv g_ic_i$ is inside a ball of radius
$C$ $\omega$--a.s. Hence we can take $R(G)=C$. Let $D$ be the
number of elements in the ball of radius $C$ in $G$. The set
$\sss$ cannot contain more than $D$ elements, by Lemma \ref{D}.
\endproof

\medskip

\begin{lemma}\label{tripod}
Let $\ck =\co{G; x,d}$ be an asymptotic cone of $G$. Then any
subgroup $\sss< x^\omega (\Pi_1 G/\omega) (x^\omega)\iv$  which
fixes pointwise a non-degenerate tripod in a transversal tree of
$\ck$ is a conjugate of a subgroup in $\Pi B(1,R)/\omega$ for some
universal constant $R=R(G)$. In particular, the size of $\sss$ is
bounded by a universal constant $D=D(G)$.
\end{lemma}

\proof \textbf{Step 1.}\quad Let $u=\lio{u_i}, v=\lio{v_i}$ and
$w=\lio{w_i}$ be the endpoints of a tripod in a transversal tree.
All the statements about sequences $(a_i)$ below are to be
understood as holding for $\omega$--almost every $i$.

Let $R$ be the constant given by Lemma \ref{pl}, (1), for $m=1$,
$(\lambda , \kappa)=(L,C)=(1,0)$ and $\theta=\frac{1}{3}$. Consider
a point { $\pi_i$ } in $\onn_R\left(\Sat ([u_i,v_i]) \right)$
minimizing the distance to $w_i$. Then $\pi=\lio{\pi_i}\in \Sat
([u,v])$ and $\dist (w,\pi )\leq \dist (w, [u,v])$. If $\pi\not\in
[u,v]$ then $\pi\in A$ for some piece $A$ intersecting $[u,v]$ in a
point $p$. Since $\pi \neq p$, Lemma 2.28 in \cite{DS} implies that
$[w,\pi ]\sqcup [\pi , p]\sqcup [p,u]$ is a geodesic. This
contradicts the fact that $u$ and $w$ are in the same transversal
tree. Hence $\pi \in [u,v]$ and it is the center of the tripod
$(uvw)$. In particular it follows that the distances from $\pi_i$ to
$u_i$, $v_i$ and $w_i$ respectively are of order $d_i$.

Let $g=(g_i)^\omega$ be an element in $\sss$.

\medskip

\noindent \textit{Notation:}\quad In the sequel, for any
$a=\lio{a_i}\in \ck$ we use the notations $a_i'$ to denote $g_ia_i$
and respectively $a'$ to denote $\lio{a_i'}=g\cdot a$.

\medskip

Note that if $a$ is in the tripod $(uvw)$ then $a'=a$, that is
$\dist (a_i,a_i')=o(d_i)$.

With the above notation $\pi_i'$ is a point in $\onn_R\left(\Sat
([u_i',v_i']) \right)$ minimizing the distance to $w_i'$, and its
distances to $u_i',v_i'$ and respectively $w_i'$ are of order $d_i$.

It suffices to show that $\dist (\pi_i, \pi_i')=O(1)$. Indeed, if
this is true, then $(\pi_i)\iv g_i \pi_i$ is in a ball around $1$ of
uniformly bounded radius. Thus, in what follows we prove that $\dist
(\pi_i, \pi_i')=O(1)$.

Let $\hpi_i$ be a point in $\onn_R\left(\Sat ([u_i',v_i']) \right)$
minimizing the distance to $w_i$. Arguing as before we obtain that
$\lio{\hpi_i}$ is the center of the tripod $(uvw)$, therefore the
distances from $\hpi_i$ to $u_i,v_i,w_i$ and to their images by
$g_i$ are of order $d_i$.

According to Lemma \ref{813g} if $\dist (\hpi_i , \pi_i')\geq D_0$
then $[w_i, \hpi_i]\sqcup [\hpi_i ,\pi_i'] \sqcup [\pi_i', w_i']$ is
a quasi-geodesic. This contradicts the fact that $\dist
(w_i,w_i')=o(d_i)$, while $\dist (w_i , \hpi_i)$ and $\dist (\pi'_i,
w_i')$ are of order $d_i$. We conclude that $\dist (\hpi_i ,
\pi_i')\leq D_0$.

Thus, in order to finish the argument, we need to prove that $\dist
(\pi_i, \hpi_i)=O(1)$.

\medskip

\textbf{Step 2.}\quad For every $\varepsilon >0$ take $\bar{u}_i,
\bar{v}_i\in [u_i,v_i]$, at distance $\varepsilon d_i$ from $u_i$
and $v_i$ respectively. By Lemma \ref{pl}, (1), there exists a point
in $[u_i, \bar{u}_i]$ contained in $\nn_R(\Sat [u_i',v_i'])$ and a
similar point in $[\bar{v}_i, v_i]$. Lemma \ref{pl}, (2), implies
that $[\bar{u}_i,\bar{v}_i]\subset \nn_{\tau R} (\Sat [u_i',v_i'])$.

Take $y_i$ an arbitrary point on $[\bar{u}_i,\bar{v}_i]$ and
assume that $y_i\not\in \nn_{\tau R} ([u_i',v_i'])$. Then $y_i\in
\nn_{\tau R} (A_i)$, with $A_i\subset \Sat ([u'_i,v'_i])$. Let
$\bar{e}_i,\bar{f}_i$ be the entry and respectively exit point of
$[\bar{u}_i,\bar{v}_i]$ in $\onn_{\tau R} (A_i)$. Likewise let
$e_i',f_i'$ be the entry and respectively exit point of
$[u'_i,v_i']$ in $\onn_{M} (A_i)$. Lemmas \ref{mm} and \ref{813}
imply that if $\dist (\bar{e}_i,e_i')\geq D_0$ then $[u_i,
\bar{e}_i]\sqcup [\bar{e}_i,e_i'] \sqcup [e_i',u_i']$ is a
quasi-geodesic. This contradicts the fact that $\dist
(u_i,u_i')=o(d_i)$ while $\dist (u_i, \bar{e}_i)$ is of order
$d_i$. Thus $\dist (\bar{e}_i,e_i')\leq D_0$. Arguing similarly we
obtain that $\dist (\bar{f}_i,f_i')\leq D_0$.

If $\dist (\bar{e}_i,\bar{f}_i)$ is large enough then by Lemma
\ref{pl}, (3), $A_i\subset \Sat [u_i,v_i]$, hence $A_i\subset \Sat
[u_i,v_i]\cap \Sat [u_i',v_i']$. If not it follows that $y_i$ is
also at distance $O(1)$ from $\{e_i', f_i'\}$.

We have thus obtained that every $y_i\in [u_i,v_i]$ such that $\dist
(u_i,y_i)$ and $\dist (y_i, v_i)$ is of order $d_i$ is either at
distance $O(1)$ from $[u_i',v_i']$ or it is contained in
$\nn_R(A_i)$ for some $A_i\subset \Sat [u_i,v_i]\cap \Sat
[u_i',v_i']$, $y_i$ at large distance from the entrance and exit
point of $[u_i,v_i]$ in $\onn_M(A_i)$.

A similar statement can be formulated for points $y_i'\in
[u_i',v_i']$ and $[u_i,v_i]$.

\medskip

\noindent {\textbf{Step 3.}}\quad Suppose that $\pi_i\in \onn_R
([u_i,v_i])$. Likewise $\pi_i'\in \onn_R ([u_i',v_i'])$.  The
argument in Step 2 implies that $\pi_i'\in \nn_{\tau R+R}\left(\Sat
[u_i,v_i] \right)$. Hence $\hpi_i$ is in
$\onn_{R+D_0}([u_i',v_i'])\cap \onn_{(\tau+1)R+D_0} (\Sat
[u_i,v_i])$. Lemmas \ref{endsg} and \ref{mm} imply that the entrance
point $\tpi_i$ of $[w_i, \hpi_i]$ in $\onn_{(\tau+1)R+D_0} (\Sat
[u_i,v_i])$ and $\pi_i$ are at distance $O(1)$. Then $\dist (\tpi_i,
[u_i,v_i])=O(1)$. Also by Step 2 we can deduce that $\tpi_i$ is at
distance $O(1)$ from $\Sat ([u_i',v_i'])$. By Lemma \ref{pl}, (1),
$\tpi_i$ is at a distance from $\hpi_i$ which is at most thrice the
distance from $\Sat ([u_i',v_i'])$, otherwise $[\tpi_i ,\hpi_i]$
would intersect $\nn_R(\Sat ([u_i',v_i']))$, a contradiction. Thus
$\dist (\tpi_i, \hpi_i)=O(1)$ and $\dist (\tpi_i, \pi_i)=O(1)$,
which finishes the argument.

\medskip

\noindent {\textbf{Step 4.}}\quad  Suppose that $\pi_i\in \onn_R
(A_i)$ with $A_i$ a left coset in  $\Sat ([u_i,v_i])$, and that
$\pi_i$ is at large distance from $[u_i,v_i]$. Then also $\pi_i'\in
\onn_R (A_i')$ with $A_i'$ a left coset in $\Sat ([u_i',v_i'])$ and
$\pi_i'$ is at large distance from $[u_i',v_i']$.

Let $e_i, f_i$ be the entrance and, respectively, the exit point{s}
of $[u_i,v_i]$ { in } $\onn_M(A_i)$. By the argument in Step 2,
$e_i$ is at distance $O(1)$ from a point $y_i$ in $[u_i',v_i']$ or
it is contained in $\nn_M (B_i)$ for some $B_i\subset \Sat
([u_i,v_i])\cap \Sat ([u_i',v_i'])$, far from the extremities of
$[u_i,v_i]\cap \nn_M (B_i)$. Without loss of generality we may
suppose that either $y_i\in [f_i',v_i']$, or $\nn_M(B_i)$ intersects
$[f_i',v_i']$. In the second case let $y_i$ be the entrance point of
$[f_i',v_i']$ into $\onn_M(B_i).$

\begin{center}
\begin{figure}[!ht]
%TeXCAD Options
%\grade{\on}
%\emlines{\off}
%\epic{\off}
%\beziermacro{\on}
%\reduce{\on}
%\snapping{\off}
%\quality{8.00}
%\graddiff{0.01}
%\snapasp{1}
%\zoom{4.0000}
\unitlength .7mm % = 2.85pt
\linethickness{0.4pt}
\ifx\plotpoint\undefined\newsavebox{\plotpoint}\fi % GNUPLOT compatibility
\begin{picture}(171.5,83)(-20,0)
\qbezier(4.25,14.5)(38,18.63)(74.75,12.25)
\qbezier(86,19)(130.63,27.25)(170.75,15.5)
\put(4.25,10.5){\makebox(0,0)[cc]{$u_i$}}
\put(74,9.25){\makebox(0,0)[cc]{$v_i$}}
\put(86.25,15.5){\makebox(0,0)[cc]{$u_i'$}}
\put(171.5,12){\makebox(0,0)[cc]{$v_i'$}}
\qbezier(24.25,10.25)(24.5,35.88)(34.75,41)
\qbezier(122.25,15.5)(122.5,41.13)(132.75,46.25)
\qbezier(34.75,40.5)(38.75,42.5)(43.75,40.5)
\qbezier(132.75,45.75)(136.75,47.75)(141.75,45.75)
\qbezier(43.75,40.5)(52.38,36.88)(48.5,9.75)
\qbezier(141.75,45.75)(150.38,42.13)(146.5,15)
\qbezier(48.25,9.75)(47.75,4.88)(39.25,4.5)
\qbezier(146.25,15)(145.75,10.13)(137.25,9.75)
\qbezier(39.25,4.5)(24.25,2.13)(24.25,10.25)
\qbezier(137.25,9.75)(122.25,7.38)(122.25,15.5)
\put(37.5,26.75){\makebox(0,0)[cc]{$A_i$}}
\put(135.5,32){\makebox(0,0)[cc]{$A_i$}}
%\emline(21,20)(21.75,25.25)
\multiput(21,20)(.032609,.228261){23}{\line(0,1){.228261}}
%\end
%\emline(119,25.25)(119.75,30.5)
\multiput(119,25.25)(.032609,.228261){23}{\line(0,1){.228261}}
%\end
%\emline(22.75,28.75)(24,34)
\multiput(22.75,28.75)(.0328947,.1381579){38}{\line(0,1){.1381579}}
%\end
%\emline(120.75,34)(122,39.25)
\multiput(120.75,34)(.0328947,.1381579){38}{\line(0,1){.1381579}}
%\end
%\emline(26.25,38.5)(30,43)
\multiput(26.25,38.5)(.03348214,.04017857){112}{\line(0,1){.04017857}}
%\end
%\emline(124.25,43.75)(128,48.25)
\multiput(124.25,43.75)(.03348214,.04017857){112}{\line(0,1){.04017857}}
%\end
%\emline(33,44)(38.5,45)
\multiput(33,44)(.183333,.033333){30}{\line(1,0){.183333}}
%\end
%\emline(131,49.25)(136.5,50.25)
\multiput(131,49.25)(.183333,.033333){30}{\line(1,0){.183333}}
%\end
%\emline(41.5,44.5)(46.5,42.75)
\multiput(41.5,44.5)(.0961538,-.0336538){52}{\line(1,0){.0961538}}
%\end
%\emline(139.5,49.75)(144.5,48)
\multiput(139.5,49.75)(.0961538,-.0336538){52}{\line(1,0){.0961538}}
%\end
%\emline(47.75,41.25)(50.25,37)
\multiput(47.75,41.25)(.0333333,-.0566667){75}{\line(0,-1){.0566667}}
%\end
%\emline(145.75,46.5)(148.25,42.25)
\multiput(145.75,46.5)(.0333333,-.0566667){75}{\line(0,-1){.0566667}}
%\end
%\emline(51,34.75)(52,29)
\multiput(51,34.75)(.033333,-.191667){30}{\line(0,-1){.191667}}
%\end
%\emline(149,40)(150,34.25)
\multiput(149,40)(.033333,-.191667){30}{\line(0,-1){.191667}}
%\end
%\emline(52,26.5)(52,21.25)
\put(52,26.5){\line(0,-1){5.25}}
%\end
%\emline(150,31.75)(150,26.5)
\put(150,31.75){\line(0,-1){5.25}}
%\end
%\emline(52,19)(51.25,12.5)
\multiput(52,19)(-.032609,-.282609){23}{\line(0,-1){.282609}}
%\end
%\emline(150,24.25)(149.25,17.75)
\multiput(150,24.25)(-.032609,-.282609){23}{\line(0,-1){.282609}}
%\end
%\emline(50.5,10)(48.75,5.25)
\multiput(50.5,10)(-.0336538,-.0913462){52}{\line(0,-1){.0913462}}
%\end
%\emline(148.5,15.25)(146.75,10.5)
\multiput(148.5,15.25)(-.0336538,-.0913462){52}{\line(0,-1){.0913462}}
%\end
%\emline(46.75,3)(42,2.25)
\multiput(46.75,3)(-.206522,-.032609){23}{\line(-1,0){.206522}}
%\end
%\emline(144.75,8.25)(140,7.5)
\multiput(144.75,8.25)(-.206522,-.032609){23}{\line(-1,0){.206522}}
%\end
%\emline(39.75,2)(33.5,2)
\put(39.75,2){\line(-1,0){6.25}}
%\end
%\emline(137.75,7.25)(131.5,7.25)
\put(137.75,7.25){\line(-1,0){6.25}}
%\end
%\emline(30.5,2)(25,4.25)
\multiput(30.5,2)(-.0820896,.0335821){67}{\line(-1,0){.0820896}}
%\end
%\emline(128.5,7.25)(123,9.5)
\multiput(128.5,7.25)(-.0820896,.0335821){67}{\line(-1,0){.0820896}}
%\end
%\emline(23,5.75)(20.75,9.25)
\multiput(23,5.75)(-.0335821,.0522388){67}{\line(0,1){.0522388}}
%\end
%\emline(121,11)(118.75,14.5)
\multiput(121,11)(-.0335821,.0522388){67}{\line(0,1){.0522388}}
%\end
%\emline(19.75,12)(20.25,17.75)
\multiput(19.75,12)(.033333,.383333){15}{\line(0,1){.383333}}
%\end
%\emline(117.75,17.25)(118.25,23)
\multiput(117.75,17.25)(.033333,.383333){15}{\line(0,1){.383333}}
%\end
\put(17.5,19.25){\makebox(0,0)[cc]{$e_i$}}
\put(55,17.75){\makebox(0,0)[cc]{$f_i$}}
\put(114.25,24.5){\makebox(0,0)[cc]{$e_i'$}}
\put(153,22){\makebox(0,0)[cc]{$f_i'$}}
\qbezier(45.25,43.25)(59.5,52.38)(68.75,70)
\qbezier(69,70)(90.88,66)(109.25,53)
\qbezier(109.5,53.25)(120.13,45.88)(125.25,45)
\qbezier(125.25,45)(131.25,43.63)(149.5,20)
\qbezier(142,49)(153,65)(162,83)
\put(68,72.5){\makebox(0,0)[cc]{$w_i$}}
\put(116.25,52){\makebox(0,0)[cc]{$\tilde\pi_i$}}
\put(125.5,49){\makebox(0,0)[cc]{$\hat\pi_i$}}
\put(146.25,49.5){\makebox(0,0)[cc]{$\pi_i'$}}
\put(165.25,82.5){\makebox(0,0)[cc]{$w_i'$}}
\put(44.25,47.25){\makebox(0,0)[cc]{$\pi_i$}}
\put(20,15.75){\circle*{1}} \put(51.5,15.25){\circle*{1}}
\put(45.5,43){\circle*{1}} \put(69,69.75){\circle*{1}}
\put(117,48){\circle*{1}} \put(125.5,44.5){\circle*{1}}
\put(142.25,49){\circle*{1}} \put(162,82.75){\circle*{1}}
\put(149.75,20){\circle*{1}} \put(118.25,22.25){\circle*{1}}
\end{picture}

%\centering
\caption{Step 4 in Lemma \ref{tripod}.} \label{fig4}
\end{figure}
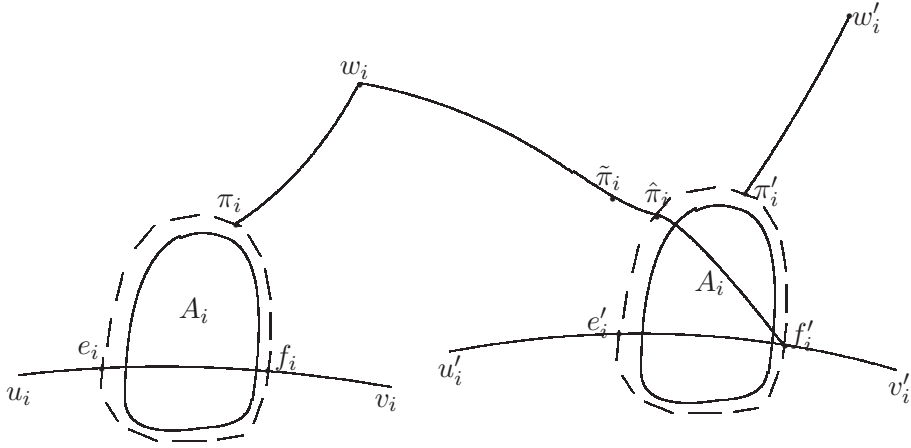
\end{center}

Since $\dist (\hpi_i, \pi_i')\leq D_0$, { we have } $\hpi_i\in
\onn_{R+D_0} (A_i')$ and { that $\hpi_i$ } is at { distance $> O(1)$
} from $[u_i',v_i']$. In particular { $\dist(\hpi_i, f_i')>O(1)$}.
Lemma \ref{813g} implies that $\q_i=[w_i, \hpi_i]\sqcup [\hpi_i,
f_i']$ is a quasi-geodesic. By Step 2, $f_i'$ is at distance $O(1)$
from $\Sat ([u_i,v_i])$. Let $D_1\geq R$ be such that $f_i'\in
\nn_{D_1}\left( \Sat ([u_i,v_i])\right)$. Lemmas \ref{endsg} and
\ref{mm} imply that the entrance point $\tpi_i$ of $\q_i$ in
$\onn_{D_1} (\Sat ([u_i,v_i]))$ and $\pi_i $ are at distance $O(1)$.
In particular $\tpi_i$ is at distance $O(1)$ from $A_i$ and at { $>
O(1)$ }distance from $[u_i,v_i]$.

\medskip

\noindent {\textbf{Step 4.a.}}\quad Assume that $\tpi_i \in [\hpi_i,
f_i']$. Then $\tpi_i \subset \nn_{\tau R} (A_i')$. If $\dist (\tpi_i
, [e_i', f_i'])$ is large, as both $\tpi_i$ and $[e_i', f_i']$ are
at distance $O(1)$ from both $A_i'$ and $\Sat ([u_i,v_i])$ it would
follow that $A_i'\subset \Sat ([u_i,v_i])$. This would contradict
the fact that $\tpi_i$ is the entrance point of $\q_i$ in
$\onn_{D_1} (\Sat ([u_i,v_i]))$.

Hence $\dist \left(\tpi_i , [e_i', f_i']\right)$ is $O(1)$.
Moreover, the fact that $A_i'$ cannot be in $\Sat ([u_i,v_i])$
implies that $\dist (e_i', f_i')$ is $O(1)$. Hence $\dist (e_i',
\tpi_i)=O(1)$, and since $\dist (\tpi_i, \pi_i)$ is likewise $O(1)$,
we deduce that $\dist (\pi_i, e_i')$ is $O(1)$.

By Step 2, either $e_i'$ is at distance $O(1)$ from $[u_i,v_i]$ or
it is in $\nn_M(B_i)\cap [u_i',v_i']$, far from the extremities of
the intersection, where $B_i\subset \Sat ([u_i,v_i]) \cap \Sat
([u_i',v_i'])$. In the first case, $\pi_i$ also is at distance
$O(1)$ from $[u_i,v_i]$. We are then back { in } the case of Step 3,
with a constant possibly larger than $R$, and we can argue
similarly.

In the second case, we have that $\pi_i$ is at distance $O(1)$ from
$B_i$. If $B_i\neq A_i$ then by Lemma \ref{pl}, (4), $\pi_i$ is at
distance $O(1)$ from $[u_i,v_i]$ and we argue as previously. Assume
therefore that $B_i=A_i$. By isometry $\dist (e_i, f_i)$ is also
$O(1)$. On the other hand, by the argument in Step 2, $e_i$ is at
distance $O(1)$ from the entrance point of $[u_i',v_i']$ into
$\onn_M(B_i)$ at the same for $f_i$ and the exit point. It follows
that $\diam \onn_M(B_i) \cap [u_i',v_i']$ is $O(1)$ and that $\dist
(e_i,e_i')$ is $O(1)$. Hence $\dist (\pi_i, e_i)=O(1)$ and we are
back again to Step 3.

\medskip

\noindent {\textbf{Step 4.b.}}\quad Assume that $\tpi_i \in [w_i,
\hpi_i]$. If $\nn_M(A_i)\cap [u_i,v_i]$ is too large, by Step 2 and
Lemma \ref{pl}, (4), it follows that $A_i\subset \Sat ([u_i',
v_i'])$. Together with the fact that $\tpi_i$ is in a tubular
neighborhood of $A_i$, with the choice of $\hpi_i$ and with Lemma
\ref{mm}, this implies that $\dist (\tpi_i,\hpi_i)=O(1)$. Hence we
may assume in what follows that $\nn_M(A_i)\cap [u_i,v_i]$ has small
diameter, and same for its image.

Let $D_2$ be the maximum between $R+D_0$ and $\dist (\tpi_i, A_i)$.
Let $[z_i, z_i']$ be either the sub-arc of $[\tpi_i, \hpi_i]$ of
extremities the exit point from $\onn_{D_2} (A_i)$ and the entry
point into $\onn_{D_2} (A_i')$ or a degenerate segment composed of
one point in $[\tpi_i, \hpi_i]\cap \onn_{D_2} (A_i)\cap \onn_{D_2}
(A_i')$.

Lemma \ref{pl}, (4), applied to the polygonal line ${\mathfrak
l}=[e_i, y_i]\cup [y_i,f_i']$, to the left cosets $A_i$ and $A_i'$
and the points $z_i,z_i'$ implies that $\{ z_i, z_i'\}$ is in
$\nn_\varkappa ({\mathfrak l})$.

By Lemma \ref{mm}, $\dist (z_i', \hpi_i )=O(1)$, which implies that
$\hpi_i$ is at distance $O(1)$ from $[e_i, y_i]\cup [y_i,f_i']$.

If $\dist (e_i,y_i)=O(1)$ then it follows that $\hpi_i$ is at
distance $O(1)$ from $[y_i,f_i']$. This implies that $\pi_i'$ is at
distance $O(1)$ from $[u_i', v_i']$, and we are back in the case of
Step 3.

The other possibility is that $\dist (e_i,y_i)$ is large, which by
Step 2, corresponds to the case when $e_i$ is in $\nn_M(B_i')\cap
[u_i,v_i]$, far from the extremities of this intersection, where
$B_i'\subset \Sat ([u_i,v_i])\cap \Sat ([u_i',v_i'])$.

Assume that $A_i'\neq B_i'$. Then one can apply Lemma \ref{813}
 {twice } and deduce that $[e_i, y_i]\sqcup [y_i,f_i']\sqcup
[f_i',\hpi_i]$ is a quasi-geodesic. Therefore $\dist (\hpi_i, [e_i,
y_i]\cup [y_i,f_i'])=O(1)$ implies that $\dist (\hpi_i, f_i')$ is
$O(1)$. But then it follows that $\dist (\pi_i', [u_i',v_i'])$ is
$O(1)$, and we can again argue as in Step 3.

Assume that $A_i'=B_i'$. By the argument in Step 2, $e_i'$ and
$f_i'$ are at distance $O(1)$ from the entrance and exit points of
$[u_i,v_i]$ into $\onn_M (A_i')$. By the argument in the beginning
of the current step, $\diam\, \onn_M (A_i')\cap [u_i',v_i']$ is
$O(1)$. It follows that $\diam\, \onn_M (A_i')\cap [u_i,v_i]$ is
$O(1)$, in particular $e_i$ is at distance $O(1)$ from $e_i'$. We
can then repeat the argument done above in the case when $\dist
(e_i,y_i)=O(1)$, with $y_i$ replaced by $e_i'$, deduce that $\hpi_i$
is at distance $O(1)$ from $[y_i,f_i']$, hence that $\pi_i'$ is at
distance $O(1)$ from $[u_i', v_i']$, and get back to Step 3.
\endproof

\begin{cor}\label{triple}
Let $\ck =\co{G; x,d}$ be an asymptotic cone of $G$. There exists
a constant $R=R(G)$ such that the following holds. Any subgroup
$\sss< x^\omega (\Pi_1 G/\omega) (x^\omega)\iv$  which fixes three
points not in the same piece nor on the same transversal geodesic
in $\ck$ is conjugate to a subgroup in $\Pi B(1,R)/\omega$.
\end{cor}

\proof Let $u,v,w$ be the three points. If any of the strict
saturations $\Sato \{a,b\}$ with $a\neq b$, $a,b\in \{u,v,w\}$,
contains a piece then $\sss$ stabilizes the piece and fixes a point
outside it, and we may apply Lemma \ref{l6}.

If not, then the three points are in the same transversal tree, as
vertices of a tripod. We apply Lemma \ref{tripod} in this
case.\endproof

\subsection{Homomorphisms into relatively hyperbolic groups}
\label{homs}

The following observation of Bestvina and Paulin is well known
(see for instance \cite{Best}).

\begin{lemma}\label{Pau} Let $\Lambda$ and $\Gamma$
be two finitely generated groups, let $S=S\iv$ be a finite set
generating $\Lambda$ and let $\dist$ be a word metric on $\Gamma$.
Given $\phi_n:\Lambda\to \Gamma$ an infinite sequence of
homomorphisms, one can associate to it a sequence of positive
integers defined by
\begin{equation}\label{dn}
d_n=\inf_{x\in \Gamma}\sup_{a\in S}\dist(\phi_n(a)x,x)\, .
\end{equation}

If $(\phi_n)$ are pairwise non-conjugate in $\Gamma$ then
$\lim_{n\to \infty}d_n=\infty $.
\end{lemma}

\begin{remark}\label{xa}
For every $n\in \N$, $d_n=\dist(\phi_n(a_n)x_n,x_n)$ for some
$x_n\in \Gamma$ and $a_n\in S$.
\end{remark}

Let $\Lambda =\la S\ra$, $(\phi_n)$ and $(d_n)$ be as in Lemma
\ref{Pau}, with $\Gamma =G$.

Consider an arbitrary ultrafilter $\omega$. According to Remarks
\ref{xa} and \ref{udisj}, there exists $a\in S$ and $x_n\in G$
such that $d_n=\dist(\phi_n(a)x_n,x_n)$ $\omega$--a.s.

\begin{lemma} \label{Pau1}
Under the assumptions of Lemma \ref{Pau}, the group $\Lambda$ acts
on the asymptotic cone $\ck_\omega=\co{G; (x_n), (d_n)}$ by
isometries, without a global fixed point, as follows:
\begin{equation}\label{act}
  g\cdot \lio{x_n}=\lio{\phi_n(g)x_n}\, .
\end{equation}

This defines a homomorphism $\phi_\omega$ from $\Lambda$ to the
group $x^\omega (\Pi_1 \Gamma/\omega) (x^\omega)\iv$ of isometries
of~$\ck_\omega$.
\end{lemma}

The action in Lemma \ref{Pau1} of a group $\Lambda$ on the
asymptotic cone $\ck_\omega$ (which  is tree-graded with respect
to limits of sequences of cosets from ${\cal G}$ by Theorem
\ref{tgr}) satisfies the hypotheses in Theorem \ref{split1} if one
more condition on $(\phi_n)$ holds.

\begin{definition}
A homomorphism $\phi : \Lambda \to G$ is called \textit{parabolic}
if its image is a parabolic group.
\end{definition}

\begin{proposition}\label{i}
Suppose that a finitely generated group $\Lambda =\la S\ra$ has
infinitely many non-parabolic homomorphisms $\phi_n$ into a
relatively hyperbolic group $G$, which are pairwise non-conjugate
in $G$.

Then the action of $\Lambda$ on an asymptotic cone $\ck_\omega$ of
$G$ defined by $(\phi_n)$ as in (\ref{act}) satisfies the
properties (i),(ii) and (iii) of Theorem \ref{split1}.
\end{proposition}

\medskip

The proof is done in several steps.

\medskip

As before (see the notation before Theorem \ref{1}), we use the
following notation:

\begin{itemize}
\item $\calc_1(\Lambda,\omega)$ is the set of stabilizers in $\Lambda$
of proper subsets in $\ck_\omega$ such that all their finitely
generated subgroups stabilize pairs of pieces in $\ck_\omega$;

\item $\calc_2(\Lambda,\omega)$ is the set of stabilizers in
$\Lambda$ of pairs of points in~$\ck_\omega$ that are not in the
same piece;

\end{itemize}

\begin{lemma}\label{A}
Let $K$ be a subgroup in $\Lambda$ such that all its finitely
generated subgroups stabilize pairs of pieces. Then $\phi_\omega(K)$
is a conjugate of a subgroup in the finite set $\Pi B(1,R)/ \omega $
for some uniform constant $R$.
\end{lemma}

\proof  Take an increasing sequence
\begin{equation}\label{seqk}
K_1\subset K_2\subset ...\subset K_i\subset ...
\end{equation}
of finitely generated subgroups of $K$ such that $K=\bigcup K_i$.
By hypothesis each $K_i$ stabilizes two distinct pieces, that is
two different nonempty limits $\lio{g_nH}$ and
  $\lio{h_nH'}$. By Lemma \ref{stabc}, for every $k\in K_i$
  $\omega$--a.s. $\phi_n(k) \in g_n Hg_n\iv \cap h_n H'h_n\iv$.
  In particular, given a finite generating set $S_i$ of $K_i$,
$\omega$--a.s. $\phi_n(S_i) \subset g_n Hg_n\iv \cap h_n H'h_n\iv$.
Then also $\phi_n(K_i)$ is contained in $g_n Hg_n\iv \cap h_n
H'h_n\iv$ $\omega$--a.s. By Lemma \ref{conj}, there exists an
element $a_n(i)$ in $g_n H$ such that $\phi_n(K_i)$ is conjugate by
$a_n(i)$ to a subgroup inside $B(1,R)$, for some uniform constant
$R>0$, $\omega$--a.s. Since there are finitely many subgroups inside
$B(1,R)$, Remark \ref{udisj} implies that there exists a finite
subgroup $U_i$ in $B(1,R)$ such that $\phi_n(K_i) = a_n(i)\iv U_i
a_n(i)$, $\omega$--a.s.

Again because there are finitely many subgroups inside $B(1,R)$,
there exists a subgroup $U$ inside $B(1,R)$ and a subsequence
$i_0<i_1<i_2<...$ such that $U_{i_m}=U$ for every $m\geq 0$. In
particular, since $\phi_n(K_{i_0}) < \phi_n(K_{i_m})$ for every
$m>0$ and both groups have the same cardinal as $U$ $\omega$--a.s.,
it follows that $\phi_n(K_{i_0}) = \phi_n(K_{i_m})$ $\omega$--a.s.
Hence for every $j\geq i_0$, $\phi_n(K_{j}) \phi_n(K_{i_0})$
$\omega$--a.s. Consequently, given $a_n=a_n(i_0)$, the group
$\phi_n(K_{j})$ is equal to $a_n\iv U a_n$, $\omega$--a.s.

Now an arbitrary element $k$ in $K$ is contained in some $K_j$ with
$j\geq i_0$, hence $\phi_n(k)\in \phi_n(K_{j})$ $\omega$--a.s. Then
$a_n$ conjugates $\phi_n(k)$
 to an element of $U$ $\omega$--a.s. Hence $(a_n)^\omega$ conjugates
 $\phi_\omega(k)$ to an element in the finite set $\Pi U/\omega$.
 \endproof

\begin{lemma}\label{trgeod}
Let $\g$ be a non-trivial geodesic segment in a transversal tree
of $\ck_\omega$, and let $L$ be the pointwise stabilizer of $\g$
in $\Lambda$. Then up to conjugacy, $\phi_\omega (L)$ is an
extension of a subgroup from
 $\Pi B(1,R)/\omega $ by an Abelian group.
\end{lemma}

\proof \textbf{Step 1.} Let $\gamma$ be an element in $L$. We
associate with it a sequence of translation numbers.

The geodesic $\g$ is the $\omega$--limit of geodesics $\g_n$ of
length $O(d_n)$. Let $a_n$ and $b_n$ be the initial and terminal
points of $\g_n$. The fact that $\gamma$ fixes $\g$ means that the
left action of $\phi_n(\gamma)$ on $\g_n$ must move it within
$o(d_n)$ distance from itself.

Note that $\g_n$ intersects each coset from $\cal G$ by a
sub-geodesic of length $o(d_n)$ $\omega$--a.s.

Let $\Sat(\g_n)$ be the saturation of $\g_n$ in the sense of
Definition \ref{defsat}.

Let $m_n$ be the midpoint of $\g_n$. Let $A_n$ be a coset from
$\cal G$ whose $M$--tubular neighborhood contains $m_n$ and such
that the length $\ell_n$ of the intersection of $\nn_M(A_n)$ with
$\g_n$ is maximal possible. Here $M=M(1,0)$.

If $\ell_n\geq \ell$, where $\ell=\ell(D)$ is a constant to be
defined later (see Step 2), then take $m_n'$ to be the entry point
of $\g_n$ into $\nn_M(A_n)$. We know that
$\dist(m_n',m_n)=o(d_n)$. If $\ell_n\leq \ell$, take $m_n'=m_n$.

Let $m_n(\gamma)$ be a projection of $\phi_n(\gamma)m_n'$ onto
$\g_n$.

We define the {\em $n$-th translation number} $\lambda_n(\gamma)$
as
$$\lambda_n(\gamma)=(-1)^{\epsilon}\dist(m_n',\phi_n(\gamma)m_n')$$ where $\epsilon=0$ if
$\dist(b_n,m_n(\gamma))\le \dist(b_n,m_n')$ and $\epsilon=1$
otherwise.

\medskip

\textbf{Step 2.} We are going to show that $\lambda_n$ satisfies
the quasi-homomorphism condition:
\begin{equation}\label{eq1}
\forall \gamma\, ,\, \zeta\mbox{ in }\Lambda\, ,\quad
\omega\mbox{-a.s.}\quad
|\lambda_n(\gamma\zeta)-(\lambda_n(\gamma)+\lambda_n(\zeta))|\le
\Delta \, ,\end{equation} where $\Delta$ is a universal constant.

Consider the geodesic $\phi_n(\gamma)\g_n$. By our assumption,
$\dist(\phi_n(\gamma)a_n,a_n)$ and $\dist(\phi_n(\gamma)b_n,b_n)$
are $o(d_n)$. Let $\pgot_n$ be the middle third of $\g_n$. Let
$a_n'$ and $b_n'$ be the initial and terminal points of $\pgot_n$.
Then by \cite[Lemmas 4.24 and 4.25]{DS}, $\phi_n(\gamma)\pgot_n$
is in the $D$--tubular neighborhood of $\Sat(\g_n)$ for a uniform
constant $D\ge M$.

According to \cite[Lemma 4.22]{DS}, for every $t>0$, every
geodesic $\cf$ and every $A\in \cg$ either the diameter of the
intersection $\nn_t(A) \cap \nn_t(\Sat(\cf))$ is uniformly bounded
by a constant $\ell(t)$ or $A\subset \Sat(\cf)$.

Since $\pgot_n$ is a geodesic, the intersection between $\nn_M
(\phi_n(\gamma)A_n)$ and the $D$--tubular neighborhood of
$\Sat(\pgot_n)$ has diameter at least $\ell_n$.

\medskip

\textbf{Step 2.a} Suppose that  $\ell_n\geq \ell(D)$. Then
$\phi_n(\gamma)A_n\subset \Sat(\pgot_n)$.

The endpoint $\phi_n(\gamma)a_n'$ is also in the $D$--tubular
neighborhood of $\Sat(\pgot_n)$.

Suppose that $\phi_n(\gamma)a_n'$ is in the $D$--tubular
neighborhood of a coset $B_n\subseteq\Sat(\g_n)$. Since the length
of the intersection of $\g_n$ with $\nn_D(B_n)$ is $o(d_n)$,
$B_n\ne A_n$ $\omega$--a.s. Then by \cite[Lemma 4.28]{DS},
$\dist(\phi_n(\gamma)m_n',m_n(\gamma))$ is uniformly bounded
$\omega$--a.s.

Suppose now that $\phi_n(\gamma)a_n'$ is in the $D$--tubular
neighborhood of $\g_n$. Then by \cite[Lemma 8.13]{DS}, again,
 $\dist(\phi_n(\gamma)m_n',m_n(\gamma))$ is uniformly bounded $\omega$--a.s.

\medskip

\textbf{Step 2.b} Now suppose that $\ell_n< \ell(D)$. Then
$m_n'=m_n$. Since $\phi_n(\gamma)m_n$ is in the $D$--tubular
neighborhood of $\Sat(\g_n)$, $\phi_n(\gamma)m_n$ is either in
$\nn_D(\g_n)$ or in $\nn_D(A_n')$ where $A_n'\subset\Sat(\g_n)$.
In the first case $\dist(\phi_n(\gamma)m_n,m_n(\gamma))$ is
uniformly bounded $\omega$--a.s.

Suppose therefore that $\phi_n(\gamma)m_n$ is in $\nn_D(A_n')$
where $A_n'\subset\Sat(\g_n)$.

Suppose moreover that $\phi_n(\gamma)\g_n$ does not intersect
$\nn_M(A_n')$. Then the length of the intersection of
$\phi_n(\gamma)\g_n$ with $\nn_D(A_n')$ is at most $3D+1$.
Otherwise property $(\alpha_2)$ and the choice of $M$ in
Definition \ref{defsat} would imply that $\phi_n(\gamma)\g_n$
intersects $\nn_M(A_n')$. In particular $\phi_n(\gamma)m_n$ is at
distance at most $3D+1$ from the entry point of
$\phi_n(\gamma)\g_n$ into $\nn_D(A_n')$. An argument as in Step
2.a implies that this entry point is at uniformly bounded distance
from $\g_n$, hence the same holds for $\phi_n(\gamma)m_n$.

Suppose that $\phi_n(\gamma)\g_n$ intersects $\nn_M(A_n')$. Let
$c_n$ be the entry point of $\phi_n(\gamma)\g_n$ into
$\nn_M(A_n')$. If $\phi_n(\gamma)m_n$ is not in $\nn_M(A_n')$,
then its distance to $c_n$ is at most $3D+1$, otherwise one
obtains a contradiction with the fact that $c_n$ is the entry
point into $\nn_M(A_n')$.

Suppose that $\phi_n(\gamma)m_n\in \nn_M(A_n')$. Note that the
intersection of $\phi_n(\gamma)\iv\nn_M(A_n')$ with $\g_n$
contains $m_n$ and has the same length as the intersection of
$\nn_M(A_n')$ and $\phi_n(\gamma)\g_n$. Therefore these lengths
are smaller than $\ell(D)$. In particular $\dist
(\phi_n(\gamma)m_n, c_n)\leq \ell (D)$.

An argument as in Step 2.a gives that $c_n$ is at uniformly
bounded distance from $\g_n$. Therefore this is also true for
$\phi_n(\gamma)m_n$.

We conclude that $\dist(\phi_n(\gamma)m_n', m_n(\gamma))$ is
uniformly bounded $\omega$--a.s.

Thus in all cases, for some constant $D''$,
\begin{equation}\label{eq2} \dist(\phi_n(\gamma)m_n',
m_n(\gamma))<D''\quad \omega\mbox{-a.s.}\end{equation}

\medskip

\textbf{Step 2.c}\quad Now we are ready to prove (\ref{eq1}). For
simplicity, assume that $\lambda_n(\gamma), \lambda_n(\zeta)\ge 0$
(the other cases are similar). All the equalities in the proof
below are true $\omega$--a.s. By (\ref{eq2}),
$\dist(m_n',m_n(\gamma))=\lambda_n(\gamma)+O(1)$,
$\dist(m_n',m_n(\gamma\zeta))=\lambda_n(\gamma\zeta)+O(1)$,
$\dist(\phi_n(\gamma)m'_n,\phi_n(\gamma)m_n(\zeta))=\lambda_n(\zeta)+O(1)$.
The last equality implies that $\dist(m_n(\gamma),
m_n(\gamma\zeta))=\lambda_n(\zeta)+O(1)$. Combining these
equalities together, we obtain (\ref{eq1}).

\medskip

\textbf{Step 3.} Let $\gamma,\zeta$ be two elements of $L$,
$[\gamma,\zeta]=\gamma\zeta\gamma\iv\zeta\iv$ be their commutator.
Then by (\ref{eq1}), $\lambda_n(\phi_n([\gamma,\zeta]))=O(1)$
$\omega$--a.s. Therefore
$$
\dist(\phi_n([\gamma,\zeta])m_n',m_n')=O(1)\; \;
\omega\mbox{--a.s.}
$$
Therefore $|(m_n')\iv \phi_n([\gamma,\zeta])m_n'| = O(1)$
$\omega$--a.s. Hence up to conjugacy $\phi_n([\gamma,\zeta])$ is
in the ball $B(1, R)$ $\omega$--a.s. for some $R$.

This implies that up to conjugacy the set of commutators of
$\phi_\omega(L)$ is contained in the set $\Pi B(1, R)/\omega $
which is of bounded cardinality by Lemma \ref{D}. Therefore every
finitely generated subgroup $L_1\le\phi_\omega(L)$ has conjugacy
classes of bounded size, i.e. $L_1$ is an $FC$-group \cite{N}. By
\cite{N}, the set of all torsion elements of $L_1$ is finite, and
the derived subgroup of $L_1$ is finite and is generated, up to
conjugacy, by a subset of $\Pi B(1, R)/\omega$ for some $R$. There
exists only finite number of finite subgroups generated by subsets
of $\Pi B(1, R)/\omega$. Since elements of $\Pi B(1, R)$ are
sequences $(g_n)^\omega$ that are constants $\omega$--a.s., every
finite subgroup generated by a subset of $\Pi B(1,R)/\omega$ is
inside $\Pi B(1,R')/\omega$ for some $R'>R$. Hence the derived
subgroup of $\phi_\omega(L)$ is conjugate to a subgroup of $\Pi
B(1,R')/\omega$.
\endproof

\begin{lemma} \label{B} Let $L\in \calc_2(\Lambda, \omega)$.
Then $\phi_\omega (L)$ is inside a conjugate of an extension of a
subgroup in $\Pi B(1,R)/\omega $ by an Abelian group.
\end{lemma}

\noindent Let $x,y$ be two points in $\ck_\omega$ that are not in
the same piece and let $L$ be the stabilizer in $\Lambda$ of the
pair $(x,y)$.

\textbf{Case 1.} Suppose that $\Sato\{x,y\}$ contains a piece $A$.
Then the stabilizer of $x,y$ coincides with the stabilizer of $A\cup
\{x,y\}$. Either $x$ or $y$ is not in $A$. Lemma \ref{l6} implies
that $\phi_\omega (L)$ is conjugated to a subgroup in $\Pi
B(1,O(1))/ \omega$.
\medskip

\textbf{Case 2.} Suppose now that $\Sato \{x,y\}$ contains no
piece. Then $x,y$ are contained in the same transversal tree, and
so they are joined by a unique geodesic $\g$ in that transversal
tree and the stabilizer of $x,y$ coincides with the stabilizer of
$\g$. It remains to use Lemma \ref{trgeod}.
\endproof

\noindent \textit{Proof of Proposition \ref{i}.} \quad
\textbf{(i)} Obviously property (i) in Theorem \ref{split1} is
satisfied: the action of $\Lambda$ permutes pieces of
$\ck_\omega$.

\me

\textbf{(ii)} According to Lemma \ref{Pau1} there is no point in
$\ck_\omega$ fixed by the whole $\Lambda$.

If $\Lambda \cdot A= A$ for some piece $A\in \pp$ then, by Lemma
\ref{stabc}, $ \phi_\omega (\Lambda ) \subset \Pi \left( g_n Hg_n\iv
\right)/\omega$ for some sequence $(g_n)$ in $G$ and some peripheral
subgroup $H\in \mh$. In particular $\omega$--a.s. $\phi_n(S)\subset
g_n Hg_n\iv$, hence the image of $\phi_n$ is a parabolic group. This
contradicts the hypothesis that $\phi_n$ are non-parabolic
homomorphisms.

\medskip

\textbf{(iii)} We prove (iii) in two steps.

\medskip

\noindent{\textbf{(iii.a)}}  Let us prove that $\calc_1(\Lambda
,\omega )$ satisfies ACC.

Let
\begin{equation}\label{seqk2}
K_1\subset K_2\subset ...\subset K_i\subset ...
\end{equation} be an increasing sequence of subgroups from
$\calc_1(\Lambda ,\omega )$. By Lemma \ref{A}, for every $i\in \N$,
the cardinality $\card\phi_\omega (K_i)$ is bounded by a constant
$D$. The sequence
$$
\card \phi_\omega (K_1)\leq \card \phi_\omega(K_2)\leq ...\leq
\card \phi_\omega (K_i)\leq ...
$$ must stabilize. Thus we may assume
that all $\card \phi_\omega (K_i)$ are the same. It follows that for
all $i>1$, $\phi_\omega (K_i)=\phi_\omega (K_1)$.

Now since each $K_i$ is in $\calc_1(\Lambda ,\omega )$, it is the
stabilizer of some proper subset $\mathcal{M}_i$ in $\ck$. The
equality $\phi_\omega (K_i)=\phi_\omega (K_1)$ implies that for
every $k_i\in K_i$ there exists $k_1\in K_1$ such that
$\phi_n(k_i)=\phi_n(k_1)$ $\omega$--a.s. In particular $k_i$ also
stabilizes $\mathcal{M}_1$, thus $K_i\subset \Stab
(\mathcal{M}_1)=K_1$. We obtain that $K_i=K_1$ for every $i>1$.

\medskip

\noindent \textbf{(iii.b)} Let us prove that $\calc_2(\Lambda
,\omega)$ satisfies ACC. Let $L$ be an arbitrary subgroup in
$\calc_2(\Lambda ,\omega)$, that is $L$ is the stabilizer of two
points $x,y$ in $\ck_\omega$ not in the same piece.

If $\Sato \{ x,y\}$ contains a piece then $\phi_\omega (L)$ is
conjugate to a subgroup in $\Pi \, B(1,R)/\omega $, by Lemma
\ref{l6}. In this case ACC is proved with an argument as in (iii.a).

If $\Sato \{ x,y\}$ contains no piece, then $x$ and $y$ are the
endpoints of a transversal geodesic $\pgot$. Recall that the action
of $\Lambda$ on $\ck_\omega$ induces an action of $\phi_\omega
(\Lambda)$ on the $\R$-tree $T= \ck_\omega /\approx \, $. Let $\g$
be the projection of $\pgot$ onto $T$. By Lemma \ref{inv},
$\phi_\omega (\Lambda)$ is the stabilizer of $\g$. Note that the
action of $\phi_\omega (\Lambda) $ on $T$ has finite of bounded size
tripod stabilizers, by Lemma \ref{inv} and Corollary \ref{triple},
and it has (finite of bounded size)-by-Abelian arc stabilizers, by
Lemma \ref{trgeod}. It follows by Lemma \ref{stable} that for every
ascending sequence of subgroups from $\calc_2(\Lambda ,\omega)$
\begin{equation}\label{seqk3}
L_1\subset L_2\subset ...\subset L_i\subset ...
\end{equation}
the ascending sequence of images
\begin{equation}\label{iseqk3}
\phi_\omega (L_1)\subset \phi_\omega (L_2)\subset ...\subset
\phi_\omega (L_i)\subset ...
\end{equation}
 must stabilize, because it is a sequence of arc stabilizers for the action of
$\phi_\omega (\Lambda)$ on $T$. Assume that $\phi_\omega
(L_i)=\phi_\omega (L_1)$, for every $i>1$. Since $L_1$ is defined
as the stabilizer of two points $x_1,y_1$ in $\ck_\omega$, not in
the same piece, it follows that every element $l$ in $L_i$ must
also stabilize $x_1,y_1$, as $\phi_\omega (l)\in \phi_\omega
(L_1)$. Hence $L_i\subset L_1$, therefore $L_i=L_1$.
\hspace*{\fill} $\Box$

\begin{proposition}\label{iplus}
If, in addition to the assumptions of Proposition \ref{i}, the
kernel of the homomorphism $\phi_\omega$ is finite, then the
action of $\Lambda$ on $\ck_\omega$ satisfies property (iii) of
Theorem \ref{splitg}.
\end{proposition}

\proof This follows immediately from Lemmas \ref{A} and \ref{B},
and Corollary \ref{triple}.\endproof

Here are the main applications of Theorem \ref{splitg} together
with Propositions \ref{i} and \ref{iplus} to relatively hyperbolic
groups. We assume that Theorem \ref{gui} is correct (recall that
the proof is not published yet). If one does not want to assume
that, one can add the assumption that the group $\Lambda$ in these
applications is either torsion-free or finitely presented and use
Theorem \ref{sela} or \ref{bf}, respectively.

{First we show that a ``generic" finitely generated group has only
finitely many non-parabolic homomorphisms into a given relatively
hyperbolic group $G$, up to conjugacy by elements of $G$.}

As an immediate {consequence of Theorem \ref{split1} and Proposition
\ref{i}} we obtain the following statement. Recall that a group
$\Lambda$ satisfies {\em property FA} of Serre \cite{Ser} if every
action of $\Lambda$ on a simplicial tree has a global fixed point.

\begin{cor}\label{corcor} If a finitely generated group $\Lambda$
satisfies property FA then for every relatively hyperbolic group
$G$ there are only finitely many pairwise non-conjugate
non-parabolic homomorphisms $\Lambda\to G$.
\end{cor}

If $\Lambda$ is an arbitrary group, we can still obtain a lot of
information about its homomorphisms into relatively hyperbolic
groups.

\begin{definition} Let $G$ be a relatively hyperbolic group, and let $K<A$ be
two subgroups of $G$. The subgroup $K$ is called {\em locally
parabolic in} $A$ if for every finitely generated subgroup $K_1$
of $K$ there exists an embedding $\phi\colon A\to G$ such that
$\phi(K_1)$ is parabolic in $G$.
\end{definition}

\begin{remarks}\label{rem888} 1. We do not know any examples where
$K$ is locally parabolic in $A\le G$ but no embedding $A\to G$
maps $K$ into a parabolic subgroup of $G$.

2. If $K$ is finitely generated and locally parabolic in $A$ then
(obviously) some embedding $A\to G$ maps $K$ inside a parabolic
subgroup of $G$.

3. By Lemma \ref{rem1}, if parabolic subgroups do not contain
non-Abelian free subgroups then one can drop the ``finitely
generated" assumption in 2.
\end{remarks}

\begin{theorem} \label{prop888} Let $\Lambda$ be a finitely generated subgroup
in $G$ which is neither virtually cyclic nor parabolic. Assume
moreover that $\Lambda$ does not split over any subgroup $K$ of it
such that $K$ is virtually cyclic or locally parabolic in
$\Lambda$. Then the number of conjugacy classes in $G$ of
injective non-parabolic homomorphisms $\Lambda\to G$ is finite.
\end{theorem}

\proof Suppose, by contradiction, that there exists a sequence of
injective non-parabolic homomorphisms $\phi_n\colon \Lambda \to G$
pairwise non-conjugate in $G$.

By Proposition \ref{i}, the sequence $(\phi_n)$ defines an action
of $\Lambda$ on an asymptotic cone $\ck_\omega$ of  $G$ satisfying
properties (i), (ii) and (iii) of Theorem \ref{split1}. Moreover
the homomorphism $\phi_\omega$ is injective, hence by Proposition
\ref{iplus} the assumptions of Theorem \ref{splitg} also hold for
the action of $\Lambda$ on $\ck_\omega$. Consequently one of the
cases (1), (2) or (3) of Theorem \ref{splitg} must occur for
$\Lambda$.

Suppose that (1) holds. Then $\Lambda$ splits over a (finite of
uniformly bounded size)-by-Abelian-by-(virtually cyclic) subgroup
$K$. By Lemma \ref{rem1}, the subgroup $K$ is either virtually
cyclic or it is a parabolic subgroup of $G$. This contradicts the
assumption that $\Lambda$ does not split over a virtually cyclic
or locally parabolic subgroup.

Suppose (2) holds. Then $\Lambda$ splits over the pre-image $K$ by
$\phi_\omega$ of a stabilizer $\Stab(B,p)$ where $B\in \pp, p\in B$.
By Lemma \ref{stabc}, $\Stab(B,p)$ is inside the ultraproduct
$\Pi_{n\in \N} P_n/\omega$ of some maximal parabolic subgroups
$P_n<G$. This means that for every element $k\in K$ $\omega$--almost
all $\phi_n$ map $k$ to $P_n$. Therefore $K$ is locally parabolic in
$\Lambda$, a contradiction.

If (3) holds then by Lemma \ref{rem1} the group $\Lambda$ is
either virtually cyclic or parabolic. This again contradicts the
choice of $\Lambda$.
\endproof

Theorem \ref{prop888} immediately implies the following corollary.
For every subgroup $\Lambda<G$ let $N_G(\Lambda)$ (resp.
$C_G(\Lambda)$ ) be the normalizer (resp. the centralizer) of
$\Lambda$ in $G$. Clearly there exists a natural embedding
$\varepsilon$ of $N_G(\Lambda)/C_G(\Lambda)$ into the group of
automorphisms $\Aut(\Lambda)$.

\begin{cor}\label{cor646}
Suppose that $\Lambda\le G$ is neither virtually cyclic nor
parabolic, and that it does not split over a locally parabolic or
virtually cyclic subgroup. Then
$\varepsilon(N_G(\Lambda)/C_G(\Lambda))$ has finite index in
$\Aut(\Lambda)$. In particular, if $\Out(\Lambda)$ is infinite then
$\Lambda$ has infinite index in its normalizer.
\end{cor}

\proof Indeed, every automorphism of $\Lambda$ is an embedding of
$\Lambda$ into $G$.\endproof

\begin{lemma} \label{out}
Suppose that the peripheral subgroups of $G$ are not relatively
hyperbolic with respect to proper subgroups.

Suppose that $\Out(G)$ is infinite, and for some sequence
$(\phi_n)$ of coset representatives of $\Aut(G)/\Inn(G)$ define a
non-trivial action of $G$ on an asymptotic cone $\ck_\omega$ of
$G$, as in Lemma \ref{Pau1}. Then

\begin{itemize}

\item[(1)] the stabilizers of pieces of $\ck_\omega$ in $G$ are
either conjugates of subgroups in a fixed ball $B(1,R)$ or maximal
parabolic subgroups;

\item[(2)] the stabilizer of a point $y=\lim^\omega(y_n)$ of
$\ck_\omega$ is the subgroup $$G_y=\left\{g\in G\mid \lio
{\frac{|y_n\iv\phi_n(g)y_n|}{d_n}}=0\right\}\, .$$
\end{itemize}
\end{lemma}

\proof (1) Let $L$ be a stabilizer of a piece $\lim^\omega
(\gamma_n H)$ of $\ck_\omega$, for some $\gamma_n\in G$ and a
peripheral subgroup $H$. By Lemma \ref{stabc}, $\phi_\omega(L)$ is
inside $\Pi \left( \gamma_n H \gamma_n\iv \right)/\omega$. This
means that for every $a\in L$, $\phi_n(a)$ is in $\gamma_n
H\gamma_n\iv$ $\omega$--a.s. Hence $a\in P_n=\phi_n\iv(\gamma_n
H\gamma_n\iv)$ $\omega$--a.s. Since peripheral subgroups are not
relatively hyperbolic with respect to proper subgroups, by Lemma
\ref{subgroups}, (4), the groups $P_n=\phi_n\iv(\gamma_n
H\gamma_n\iv)$ are maximal parabolic.

Suppose that $L$ contains $N+1$ distinct elements, where $N$ is
the cardinal of the ball $B(1,R)$ appearing in Lemma \ref{conj}.
Let $a_1,..., a_{N+1}$ be these elements. There exists $I\subset
\N$ with $\omega (I)=1$ and such that $\{a_1,..., a_{N+1}\}
\subset P_n$ for $n\in I$. Lemma \ref{conj} implies that for any
$n\in I$ the group $P_n$ coincides with a fixed maximal parabolic
group $P$.

Suppose that $L$ has at most $N$ elements but that it is not
conjugate to a subgroup in $B(1,R)$. A similar argument allows to
conclude that there exists $I\subset \N$ with $\omega (I)=1$ such
that for all $n\in I$,  $P_n$ coincides with some fixed $P$.

Thus $L<P$. On the other hand, $P$ clearly stabilizes the piece
$\lim^\omega (\gamma_n H)$. Hence $L=P$. This implies (1).

Statement (2) can be obtained by a direct computation.
\endproof

\begin{theorem}\label{thout} Suppose that the peripheral
subgroups of $G$ are not relatively hyperbolic with respect to
proper subgroups. If $\Out(G)$ is infinite then one of the
followings cases occurs:
\begin{itemize}
    \item[(1)] $G$ splits over a virtually cyclic subgroup;
    \item[(2)] $G$ splits over a parabolic
  [(finite of uniformly bounded size)-by-Abelian]-by-(virtually cyclic)
  subgroup;
    \item[(3)] $G$ can be represented as a fundamental group
    of a graph of groups such that each vertex group is either maximal parabolic
    or of the form $G_y$ in Lemma \ref{out}, (2), and the edge groups are parabolic; thus
    $G$ splits as an amalgamated product or an HNN extension with a
maximal parabolic subgroup $H$ as a vertex group and a proper
subgroup of $H$ as an edge group.
\end{itemize}
\end{theorem}

\proof Given a sequence $(\phi_n)$ of automorphisms as in Lemma
\ref{out}, consider the corresponding action of $G$ on an
asymptotic cone $\ck_\omega$ of $G$.

Conditions (i) and (ii) of Theorem \ref{splitg} obviously hold,
condition (iii) holds by Lemmas \ref{A} and \ref{B}, and Corollary
\ref{triple}. Thus $G$ is either in Case (1) or in Case (2) of
Theorem \ref{splitg}. Case (3) cannot occur since $G$ contains a
non-Abelian free subgroup by Lemma \ref{rem1}.

Suppose that Case (1) of Theorem \ref{splitg} occurs. By Lemma
\ref{rem1}, a (finite of uniformly bounded
size)-by-Abelian-by-(virtually cyclic)
 subgroup of $G$ is either virtually cyclic, and then Case (1) in
 the present theorem
 holds, or it is parabolic and Case (2) of the theorem
 holds.

If Case (2) of Theorem \ref{splitg} occurs then by Lemma \ref{out}
we have Case (3) of this theorem.

\endproof

The next theorem describes co-Hopfian relatively hyperbolic
groups. The following lemma answers a question in \cite{BS}.

\begin{lemma}\label{cor22} For any monomorphism $\phi\colon G \to G$ such that $\phi^k(G)$ is not
parabolic for any $k$, let $Z_k$ be the (finite) centralizer of
$\phi^k(G)$. Then the increasing union $Z$ of $Z_k$ is finite.
\end{lemma}
\proof Suppose that $Z$ is infinite. It is clear that $Z$ is
locally finite. Hence $Z$ is a parabolic subgroup by Lemma
\ref{rem1}, since it does not contain free non-Abelian subgroups.
Let $R$ be the constant from Lemma \ref{conj}. For some $k>>1$,
$Z_k$ contains more elements than the ball $B(1,R)$. Note that
conjugation by any element $g\in \phi^k(G)$ fixes elements of
$Z_k$. Let $H$ be the maximal parabolic subgroup containing $Z$.
Then $gHg\iv\cap H\ge Z_k$. By Lemma \ref{conj}, $g\in H$. Hence
$H$ contains $\phi^k(G)$, a contradiction.
\endproof

\begin{theorem} \label{cohopf} Suppose that $G$ is not co-Hopfian.  Let $\phi$ be an injective
but not surjective homomorphism $G\to G$. Then one of the
following holds:

\begin{itemize}
\item $\phi^k(G)$ is parabolic for some $k$;

\item $G$ splits over a parabolic or virtually cyclic subgroup.
\end{itemize}
\end{theorem}

\proof By Lemma \ref{cor22} and \cite[Theorem 3.1]{RS}, we can
assume that there are infinitely many pairwise non-conjugate
powers of $\phi$. Then $G$ acts on an asymptotic cone of it
$\ck_\omega$. The rest of the proof is similar to the proof of
Theorem \ref{thout} and it is left to the reader.
\endproof

\addtocontents{toc}{\contentsline {section}{\numberline { }
References \hbox {}}{\pageref{bibbb}}}
\begin{minipage}[t]{2.9 in}
\noindent Cornelia Dru\c tu\\ Department of Mathematics\\
University of Lille 1\\ and UMR CNRS 8524 \\ Cornelia.Drutu@math.univ-lille1.fr\\
\end{minipage}
\begin{minipage}[t]{2.6 in}
\noindent Mark V. Sapir\\ Department of Mathematics\\
Vanderbilt University\\
m.sapir@vanderbilt.edu\\
%http://www.math.vanderbilt.edu/$\sim$msapir\\
\end{minipage}

\end{document}